\def\no{\if01}
\def\iftwelvept{\no}

\def\ifusepdf{\no}
\def\ifpsfont{\no}

\iftwelvept
\documentclass[leqno,12pt]{amsart}
\else
\documentclass[leqno,10pt]{amsart}
\fi
\usepackage{amssymb}
\usepackage{amscd}
\usepackage{latexsym}
\usepackage{verbatim}
\usepackage[all]{xy}

\ifusepdf
\usepackage{hyperref}
\else\fi
\ifpsfont
\usepackage[T1]{fontenc}
\usepackage{times}
\else\fi

\iftwelvept
\else\fi

\theoremstyle{plain}
\newtheorem{Theorem}{Theorem}[section]

\newtheorem{Proposition}[Theorem]{Proposition}
\newtheorem{Lemma}[Theorem]{Lemma}
\newtheorem{Corollary}[Theorem]{Corollary}

\newtheorem{Claim}{Claim}[Theorem]

\theoremstyle{definition}

\newtheorem{Definition}[Theorem]{Definition}
\newtheorem{Remark}[Theorem]{Remark}
\newtheorem{Construction}[Theorem]{Construction}
\newtheorem{Example}[Theorem]{Example}

\newcommand{\ZZ}{\mathbf{Z}}
\newcommand{\TT}{\mathbb{T}}

\newcommand{\RR}{\mathbb{R}}
\newcommand{\RRR}{\mathbf{R}}
\newcommand{\CC}{\mathbb{C}}

\newcommand{\DD}{\mathbb{D}}

\newcommand{\NNNN}{\operatorname{N}}
\newcommand{\GG}{\mathcal{G}}
\newcommand{\HH}{\operatorname{\mathcal{HH}}}

\newcommand{\FF}{\mathcal{F}}

\newcommand{\AAA}{\mathcal{A}}

\newcommand{\MM}{\mathcal{M}}

\newcommand{\DDD}{\mathcal{D}}
\newcommand{\uni}{\mathbf{1}}
\newcommand{\CCC}{\mathcal{C}}

\newcommand{\PR}{\operatorname{Pr}^{\textup{L}}}

\newcommand{\SSSS}{\mathbb{S}}
\newcommand{\Rep}{\operatorname{Rep}}

\newcommand{\OO}{{\mathcal{O}}}

\newcommand{\LL}{{\mathbb{L}}}

\newcommand{\MMM}{\mathcal{M}}

\newcommand{\PPP}{\mathcal{P}}

\newcommand{\Hom}{\operatorname{Hom}}

\newcommand{\Ker}{\operatorname{Ker}}

\newcommand{\Spec}{\operatorname{Spec}}

\newcommand{\Perf}{\operatorname{Perf}}

\newcommand{\SP}{\operatorname{Sp}}
\newcommand{\SPS}{\operatorname{Sp}^{\Sigma}}

\newcommand{\Mod}{\operatorname{Mod}}

\newcommand{\SSS}{\mathcal{S}}

\newcommand{\colim}{\operatorname{colim}}
\newcommand{\CAT}{\operatorname{Cat}}
\newcommand{\Cat}{\textup{Cat}_{\infty}}

\newcommand{\Map}{\operatorname{Map}}

\newcommand{\Fun}{\operatorname{Fun}}
\newcommand{\Alg}{\operatorname{Alg}}

\newcommand{\End}{\operatorname{End}}

\newcommand{\Aff}{\operatorname{Aff}}

\newcommand{\FIN}{\operatorname{Fin}_\ast}

\newcommand{\wCat}{\widehat{\textup{Cat}}_{\infty}}
\newcommand{\wCatone}{\widehat{\textup{Cat}}}

\newcommand{\CAlg}{\operatorname{CAlg}}
\newcommand{\CCAlg}{\operatorname{CAlg}^{\le0}}
\newcommand{\CCAlgft}{\operatorname{CAlg}^{\le0,\diamondsuit}}
\newcommand{\CCAlgftec}{\operatorname{CAlg}^{\le0,\square}}
\newcommand{\DDT}{\DD_{\etwo}}
\newcommand{\AAAA}{\mathbb{A}}
\newcommand{\BBBB}{\mathbb{B}}
\newcommand{\KK}{\mathbb{K}}

\newcommand{\BAR}{\operatorname{Bar}}

\newcommand{\QC}{\operatorname{QC}}

\newcommand{\EXT}{\operatorname{Art}^{\textup{tsz}}}
\newcommand{\TSZ}{\mathsf{TSZ}}
\newcommand{\NIL}{\mathsf{Nil}}
\newcommand{\BEXT}{\operatorname{Art}^{\textup{nil}}}

\newcommand{\CST}{\mathsf{St}}
\newcommand{\FST}{\widehat{\mathsf{St}}^{\ast}}
\newcommand{\GFST}{\widehat{\mathsf{St}}^{!}}

\newcommand{\LM}{\operatorname{\mathcal{LM}}}

\newcommand{\Ind}{\operatorname{Ind}}

\newcommand{\Coh}{\operatorname{Coh}}

\newcommand{\PRST}{\PR_{\textup{St}}}
\newcommand{\Free}{\operatorname{Free}}

\newcommand{\wSSS}{\widehat{\mathcal{S}}}

\newcommand{\assoc}{\operatorname{As}}
\newcommand{\Sect}{\operatorname{Sect}}

\newcommand{\eone}{\mathbf{E}_1}
\newcommand{\etwo}{\mathbf{E}_2}
\newcommand{\eenu}{\mathbf{E}_n}
\newcommand{\einf}{\mathbf{E}_\infty}

\newcommand{\KS}{\mathbf{KS}}
\newcommand{\Lie}{\mathbf{Lie}}

\newcommand{\LMod}{\operatorname{LMod}}
\newcommand{\RMod}{\operatorname{RMod}}
\newcommand{\ST}{\operatorname{\mathcal{S}t}}

\newcommand{\Mfldf}{\mathsf{Mfld}^{\textup{fr}}}
\newcommand{\Disk}{\mathsf{Disk}}
\newcommand{\Diskf}{\mathsf{Disk}^{\textup{fr}}}

\newcommand{\DR}{\textup{dR}}
\newcommand{\PST}{\Fun(\CCAlgftec_k,\SSS)}

\newcommand{\RM}{\operatorname{\mathcal{RM}}}

\newcommand{\Def}{\operatorname{Def}}

\newcommand{\Sym}{\operatorname{Sym}}
\newcommand{\Proof}{{\sl Proof.}\quad}
\newcommand{\QED}{{\unskip\nobreak\hfil\penalty50\quad\null\nobreak\hfil
{$\Box$}\parfillskip0pt\finalhyphendemerits0\par\medskip}}

\begin{document}

\title{On D-modules of categories I}

\author{Isamu Iwanari}


\address{Mathematical Institute, Tohoku University, Sendai, Miyagi, 980-8578,
Japan}

\email{isamu.iwanari.a2@tohoku.ac.jp}

\maketitle

\section{Introduction}

A theorem of Goodwillie \cite{Good} shows that the periodic cyclic homology has homotopy invariance
in characteristic zero.
Getzler \cite{Get} constructs, under a suitable condition, a flat connection on the relative periodic cyclic homology of  an associative algebra over a base commutative ring of characteristic zero.
The connection can be thought of as a noncommutative analogue of Gauss-Manin connection.

Morita invariant property is of fundamental importance in noncommutative algebraic geometry.
Unfortunately, to the knowledge of the author, Morita invariance of the original connection remains open.
Nowadays, in view of several motivations including homological mirror symmetry,
the significance of Morita invariance of the associated connection becomes clear.
Equivariant contexts (that is, categories with group actions)
naturally arise from interesting situations: for example  matrix factorizations \cite{PMF}, schemes with group actions, etc.
However, explicit formulas often destroy symmetries,
and it is ill-suited for equivariant contexts.

Let us move to the world of higher categories.
In the present paper, we give constructions of $D$-module structures on the periodic cyclic 
homology of a stable $\infty$-category.
The theory is applicable to equivariant
contexts and has favorable features
such as Morita invariance.
Also, a generalization from the setting of algebras to stable $\infty$-categories or the likes
is useful. For example, it enables us to use exact sequences associated to
 exact sequences of categories. 
Many operations for $\infty$-categories are intolerant of explicit methods so that
they often require conceptual or geometric understandings.
We need to give a conceptual way of constructing connections.

We propose two different methods. Here is a brief overview of them.

\vspace{2mm}

(I) {\it Canonical extensions of factorization homology to mapping stacks.}
The first approach is based on a simple observation on factorization homology.
Let $\MM$ be a symmetric monoidal $\infty$-category satisfying some condition on colimits.
Factorization homology has been extensively developed 
in last decade, see \cite{And}, \cite{AF}, \cite{AFT}, \cite{Ho}, \cite{HA} for the foundation.
We also refer to the survey \cite{Ginot}.
The factorization homology theory takes a manifold $M$ and a structured algebra
$B$ such as an $\eenu$-algebra in $\MM$ as input
and produces an object of $\MM$ as output.
We would like to explain the observation in a simple setting.
Let $k\to A$ be a map of connective commutative ring spectra
and let $B$ be an $\eenu$-algebra over $A$.
Given a framed $n$-manifold $M$, we have the factorization homology $\int_MB/A$,
that is an $A$-module (spectrum).
Now let us regard $B$ as an $\eenu$-algebra over $k$ through the restriction along
$k\to A$. It gives rise to the factorization homology $\int_MB/k$, that is a $k$-module.
The point is that $\int_MB/k$ is promoted to a module spectrum over $A\otimes_kM$ and there is a canonical equivalence
\[
\big(\int_MB/k\big)\otimes_{A\otimes_kM}A\simeq \int_MB/A.
\]
The symbol $A\otimes_kM$ means the tensor of $A$ by the underlying topological space
of $M$ as a commutative algebra over $k$. 
Since $A\otimes_kM$ is the ring of functions of the mapping stack
$\Map(M,\Spec A)$ in the sense of derived geometry, we may interpret $\int_MB/k$
as an extension/deformation of $\int_MB/A$ on $S=\Spec A$ to $\Map(M,S)$.
By the gluing property (see Lemma~\ref{gluingexample}),
this extension can be globalized: $\Spec A$ is generalized to
an arbitrary derived scheme $S$ (or algebraic space) (see Thereom~\ref{geometricthm1}).
The idea is simple and is based on the local-to-global principle which is bulit in the very definition of the factrization homology.

Let us consider the special case when $M=S^1$ and $B$ is an associative algebra ($\eone$-algebra)
over 
$A$. In this case, $\int_{S^1}B/A$ is the Hochschild homology spectrum $\HH_\bullet(B/A)$.
Let us replace $B$ with an $A$-linear stable $\infty$-category $\CCC$.
In Section~\ref{categoricalextension}, we generalize the above observation
to the case of $\CCC$ (see Theorem~\ref{generalizedeasy}, Theorem~\ref{preeasy}).
The work of Ben-Zvi and Nadler \cite{BN} and Preygel \cite{P}
relates equivariant complexes on the free loop space $LS=\Map(S^1,S)$
and $2$-periodic $D$-modules in characteristic zero.
Using this bridge, we associate a $D$-module structure on the periodic cyclic
homology/complex $\mathcal{HP}_\bullet(\CCC/A)$.

\vspace{3mm}

(II) {\it Hochschild pair and moduli theory.}
Suppose that $k$ is a field of characteristic zero and $A$ is a smooth algebra over $k$.
Let $\CCC$ be an $A$-linear stable $\infty$-category.
The second approach uses the Hochschild pair, that is, the pair of the Hochschild cochain complex $\HH^\bullet(\CCC/A)$ and the Hochschild chain complex $\HH_\bullet(\CCC/A)$.
Besides, the construction makes use of the local moduli (deformation) theory related to the Hochschild pair
\cite{IMA} and the theory of correspondence between pointed formal stacks and dg Lie algebras
developed in \cite{Gai2}, \cite{H}, \cite[X]{DAG}.
The interplay of deformation functors, formal stacks, and dg Lie algebras
plays a central role.
In \cite{I} for an $A$-linear stable $\infty$-category we give a conceptual construction of an algebra structure on the Hochschild pair
$(\HH^\bullet(\CCC/A),\HH_\bullet(\CCC/A))$, encoded by means of so-called
Kontsevich-Soibelman operad $\KS$.
Moreover, in \cite{IMA} we give a moduli-theoretic interpretation of  $(\HH^\bullet(\CCC/A),\HH_\bullet(\CCC/A))$. 
Using them we extend $\HH_\bullet(\CCC/A)$ to
an $S^1$-equivariant (Ind-coherent) complex on the loop space $LS$
(so that as above we can apply the relation of loop spaces and $D$-modules to it).
In order to achieve this, we also construct a Kodaira-Spencer morphism for $\CCC$
(where $\CCC$ is regared as a family of stable $\infty$-categories over $\Spec A$).
Our construction reveals
how $(\HH^\bullet(\CCC/A),\HH_\bullet(\CCC/A))$
together with the Kodaira-Spencer morphism yields a promotion of the Hochschild homology $\HH_\bullet(\CCC/A)$
to the free loop space $LS$. 
Consider $LS$ to be a pointed formal stack
$S\to LS\to S$ over $S=\Spec A$, where $S\to LS$ is induced by constant loops.
According to the correspondence between 
pointed formal stacks and dg Lie algebras, there is an essentially unique dg Lie algebra $L$
corresponding to $LS$. Moreover, a Lie algebra module over $L$ amounts to an Ind-coherent 
complex on $LS$.
At a first glimpse, it seems like a good idea to have a Lie algebra action of $L$ on $\HH_\bullet(\CCC/A)$.
However, $LS\to S$ does not commute with the natural $S^1$-action on $LS$.
As a result, there is no way to describe an $S^1$-equivariant Ind-coherent complex on $LS$
as a Lie algebra module over $L$.
Instead of $LS$, we consider a pointed formal stack $S\to \widehat{S\times LS}\stackrel{\textup{pr}_1}{\to} S$
which is a completion along the graph of $S\to LS$ 
and admits an $S^1$-action induced by that of $LS$
(actually, there are two versions of completions
$\widehat{S\times LS}$ and $(S\times LS)^\wedge_S$, but in this introduction
we do not distinguish them).
We use a method of constructing complexes
on the formal stacks in the diagram $S \stackrel{\textup{pr}_2}{\leftarrow} \widehat{S\times S}\to \widehat{S\times LS}$
in a compatible way.

\vspace{2mm}

We introduced two different approaches. Both have their own features.
As for $\eenu$-algberas,
the extensions to the mapping stacks can be viewed as a far-reaching generalization of the extensions of Hochschild homology to the loop spaces.
In other words, the extension to the free loop space may be understood as a simple case of
the elementary observation about factorization homology.
The approach (I) has a good functoriality with respect to inputs,
that is, stable $\infty$-categories. In particular, one can use
the localization sequence.
Moreover, using the extended object defined on the free loop space and the Koszul duality 
we can study the resulting $D$-module, but we defer the study to subsequent papers. 
The approach (I) is based on the simple idea.
One may hope for various generalizations and variants.
For example, it is reasonable to consider its generalization to 
canonical extensions of factorization homology of $(\infty,n)$-categories.
The approach (II) is not so simple as the approach (I).
However, the associated complex on $LS$ has a direct link to the structure arising from
$(\HH^\bullet(\CCC/A),\HH_\bullet(\CCC/A))$ which are described as 
Lie algebra actions. It is useful to study the resulting objects.
In fact, in a future paper we apply the approach (II) to prove a generalization of Griffiths transversality theorem in this situation (cf. \cite{DC3}).
At this point, the question naturally arises: what is the relationship between the two approaches?
In a subsequent paper \cite{DC2}
we will prove that the associated two $D$-modules coincide when $S$ is affine.
This means that one can take advantages of both approaches.

\vspace{2mm}

{\it Related works.}
Getzler's method on connections on periodic cyclic homology uses 
Cartan homotopy/magic formula arising from 
the algebraic structure on the Hochschild cochain complex and the Hochschild chain complex.
Our second approach also relies on Cartan homotopy/magic formula.
In that sense, one may say shortly that the second approach is a generalization of
Getzler's original approach.
On the other hand, closer inspection reveals that
the strategy of the second approach differs quite fundamentally
from Getzler's one.
Notwithstanding,
both possess the advantage of having access to the algebraic strucuture on
 the Hochchild pair.
 Our construction uses methods which Getzler's one does not use.
The algebra $(\HH^\bullet(\CCC/A),\HH_\bullet(\CCC/A))$ over $\KS$
which is constructed in \cite{I} has Morita invariance
in the sense that an equivalence $\CCC\simeq \CCC'$
induces an equivalence of Hochschild pairs with algebraic structures. 
As mentioned above, 
we apply the moduli-theoretic description of $(\HH^\bullet(\CCC/A),\HH_\bullet(\CCC/A))$
established in \cite{IMA} and the relation between free loop spaces and $D$-modules.
These machineries allow us to obtain the bridge between $(\HH^\bullet(\CCC/A),\HH_\bullet(\CCC/A))$
and the geometry of the free loop space, and the $D$-module of the periodic cyclic homology/complex.

Recently, Hoyois, Safronov, Scherotzke and Sibilla develope the theory of
categorified Chern character and categorified Grothendieck-Riemann-Roch theorem
and apply it to produce a $D$-module structure on the periodic cyclic homology/complex
\cite[Section 6]{HSSS}.
While they also use the free loop spaces,
their work is based on different ideas and techniques so that
it would be interesting to compare it with ours.

\vspace{2mm}

{\it Organization.} Let us give some instructions
to the reader.
Section 2 collects conventions and some of the notations that we will
use. In Section 3, we will review
some of background material relevant to this paper.
In Section 4, we carry out the first construction (I):
we make an observation about factorizations homology. 
The results are presented in Theorem~\ref{restrictionfactorization}, Theorem~\ref{diagramfactorization}
and Theorem~\ref{easymain}.
Most of Sections 5--7 is devoted to the second construction (II).
Section 5 contains the construction of Kodaira-Spencer morphisms.
Section 6 contains a main procedure in (II). 
The main conclusion in Section 7 is Theorem~\ref{cyclicextension}.
Section 4 and Sections 5--7 are logically independent from one another.
The reader can read a preferred one.
In Section 8, we apply the relation between free loop spaces and 2-periodic $D$-modules
to show how to obtain $D$-module of periodic cyclic homology from equivariant 
complexes constructed in Section 4 and Section 5--7 (see Theorem~\ref{maindmodule}).

\vspace{2mm}

{\it Acknowledgements.}
The author would like to thank Takuo Matsuoka for valuable comments and advice on
Section 4 and many discussions since the summer of 2017.
The topic of this paper is included in the contents of the graduate course offered in fall term 2020
at Tohoku University.
He also thanks students for their feedbacks.
This work is supported by JSPS KAKENHI grant.

\section{Notation and Convention}

\label{NCsection}

{\it $(\infty,1)$-categories.}
Throughout this paper we use 
the language of $(\infty,1)$-categories.
We use the theory of {\it quasi-categories} as a model of $(\infty,1)$-categories.
We assume that the reader is familiar with this theory.
We will use the notation similar to that used in \cite{I}, \cite{IMA}.
A quasi-category is a simplicial set which
satisfies the weak Kan condition of Boardman-Vogt.
Following \cite{HTT}, we shall refer to quasi-categories
as {\it $\infty$-categories}.
Our main references are \cite{HTT}
 and \cite{HA}.
 To an ordinary category, we can assign an $\infty$-category by taking
its nerve, and therefore
when we treat ordinary categories we often omit the nerve $\NNNN(-)$
and directly regard them as $\infty$-categories.
 
 Here is a list of (some) of the conventions and notation that we will use:

\begin{itemize}



\item $\Delta^n$: the standard $n$-simplex

\item $\textup{N}$: the simplicial nerve functor (cf. \cite[1.1.5]{HTT})

\item $\SSS$: $\infty$-category of small spaces/$\infty$-groupoids. We denote by $\widehat{\SSS}$
the $\infty$-category of large spaces (cf. \cite[1.2.16]{HTT}).

\item $\CCC^\simeq$: the largest Kan subcomplex of an $\infty$-category $\CCC$. Namely, $\CCC^\simeq$ is the largest $\infty$-groupoid contained in $\CCC$.

\item $\CCC^{op}$: the opposite $\infty$-category of an $\infty$-category. We also use the superscript ``op" to indicate the opposite category for ordinary categories and enriched categories.


\item $\SP$: the stable $\infty$-category of spectra.

\item $\Fun(A,B)$: the function complex for simplicial sets $A$ and $B$. If $A$ and $B$ are $\infty$-categories, we regard $\Fun(A,B)$ as the functor category. 

\item $\Map_{\mathcal{C}}(C,C')$: the mapping space from an object $C\in\mathcal{C}$ to $C'\in \mathcal{C}$ where $\mathcal{C}$ is an $\infty$-category.
We usually view it as an object in $\mathcal{S}$ (cf. \cite[1.2.2]{HTT}).

\end{itemize}

{\it Operads and Algebras.}
We will use operads.
We employ
the theory of {\it $\infty$-operads} which is
thoroughly developed in \cite{HA}.
The notion of $\infty$-operads gives
one of the models of colored operads.
Here is a list of (some) of the notation about $\infty$-operads and algebras over them that we will use:

\begin{itemize}

\item $\FIN$: the category of pointed finite
sets $\langle 0 \rangle, \langle 1 \rangle,\ldots \langle n \rangle,...$
where $\langle n \rangle=\{*,1,\ldots n \}$ with the base point $*$.
We write $\Gamma$ for $\NNNN(\FIN)$. $\langle n\rangle^\circ=\langle n\rangle\backslash*$.
Notice that the (nerve of) Segal's gamma category is the opposite category
of our $\Gamma$.

\item Let $\mathcal{M}^\otimes\to \mathcal{O}^\otimes$ be a fibration of $\infty$-operads. We denote by
$\Alg_{/\mathcal{O}^\otimes}(\mathcal{M}^\otimes)$ the $\infty$-category of algebra objects (cf. \cite[2.1.3.1]{HA}).  We often write
$\Alg_{\mathcal{O}^\otimes}(\MMM)$ for $\Alg_{/\mathcal{O}^\otimes}(\MMM^\otimes)$.

\item $\CAlg(\mathcal{M}^\otimes)$: $\infty$-category of commutative
algebra objects in a symmetric
monoidal $\infty$-category $\mathcal{M}^\otimes\to \NNNN(\FIN)=\Gamma$.
When the symmetric monoidal structure is clear,
we usually write $\CAlg(\mathcal{M})$ for $\CAlg(\mathcal{M}^\otimes)$.

\item $\Mod_R^\otimes(\mathcal{M}^\otimes)$: the symmetric monoidal
$\infty$-category of
$R$-module objects
where $\mathcal{M}^\otimes$
is a symmetric monoidal $\infty$-category.
Here $R$ belongs to $\CAlg(\mathcal{M}^\otimes)$
cf. \cite[3.3.3, 4.5.2]{HA}.
We write $\Mod_R(\mathcal{M}^\otimes)$ for the underlying $\infty$-category.

\item $\Mod_R$: Suppose that $R$ is a commutative ring spectrum, i.e., a commutative algebra
object in $\SP$.
Unless stated otherwise, we write
$\Mod^\otimes_R$ for $\Mod^\otimes_R(\SP)$, that is,
the symmetric monoidal $\infty$-category of $R$-module spectra.
By $\Mod_R$ we mean the underlying $\infty$-category.
When $R$ is the
Eilenberg-MacLane spectrum $HC$ of an ordinary commutative ring $C$, we write $\Mod_C$
for $\Mod_R$ (thus
$\Mod_C$ is not the category of usual $C$-modules). 
If $\mathcal{D}^\otimes(C)$ denotes the symmetric monoidal stable $\infty$-category obtained from
the category of (possibly unbounded) chain complexes of $C$-modules by inverting quasi-isomorphisms,
there is a canonical equivalence $\Mod^\otimes_C \simeq  \mathcal{D}^\otimes(C)$ of symmetric monoidal $\infty$-categories, see \cite[7.1.2,7.1.2.13]{HA} for more details.
Let $A$ be an obejct of $\CAlg(\Mod_R)$ (see below),
and let $A'\in \CAlg(\SP)$ be the image of $A$ under the forgetful functor $\Mod_R\to \SP$. 
Then the induced functor $\Mod_A(\Mod_R)\to \Mod_{A'}(\SP)=\Mod_{A'}$ is an equivalence.
By abuse of notation, we usually write $\Mod_A$ for $\Mod_A(\Mod_R)$.

\item $\CAlg_R$: the $\infty$-category $\CAlg(\Mod_R^\otimes)$ of commutative
algebra objects in the symmetric monoidal $\infty$-category $\Mod_R^\otimes$
where $R$ is a commutative ring spectrum.
If $R$ is the sphere spectrum, we write $\CAlg$ for $\CAlg_R$.
We write $\CAlg^+_R$ for $(\CAlg_R)_{/R}\simeq \CAlg(\SP)_{R//R}$.
When $R$ is the Eilenberg-MacLane spectrum $HC$ with a commutative ring $C$,
then we write $\CAlg_C$ for $\CAlg_{HC}$.
If $C$ is an ordinary commutative ring over a field $k$ of characteristic zero,
the $\infty$-category $\CAlg_C=\CAlg(\Mod_C^\otimes)\simeq \CAlg(\mathcal{D^\otimes}(C))$ is
equivalent to the $\infty$-category
obtained from the model category of commutative
differential graded $C$-algebras by inverting quasi-isomorphisms
(cf. \cite[7.1.4.11]{HA}). Therefore,
we often regard an object of  $\CAlg_C$ as a commutative differential graded (dg) algebra
over $C$.
We shall refer to an object of $\CAlg_C$ as a commutative dg algebra over $C$.

Suppose that $A\in \CAlg_k$ is a connective commutative dg algebra over a field $k$ of characteristic zero
(a connective commutative dg algebra is a commutative dg algebra $A$ such that
$H^i(A)=\pi_{-i}(A)=0$ for $i>0$).
Let $\Mod_k^{\le0}$ be the full subcategory of $\Mod_k$ spanned by 
connective objects (those objects $M$ such that $H^i(M)=0$ for $i>0$).
We regard $A$ as an object of $\CAlg(\Mod_k^{\le0})$
and set $\CCAlg_{A}=\CAlg(\Mod_A(\Mod_k^{\le0}))$.
The symmetric monoidal fully faithful functor $\Mod_k^{\le0} \hookrightarrow \Mod_k$
exhibits $\CCAlg_{A}=\CAlg(\Mod_A(\Mod_k^{\le0}))$ as a full subcategory of $\CAlg(\Mod_A)=\CAlg_A$.

\item $\mathbf{E}^\otimes_n$: the $\infty$-operad of little $n$-cubes (cf. \cite[5.1]{HA}).
For a symmetric monoidal $\infty$-category $\CCC^\otimes$,
we write $\Alg_{n}(\CCC)$ or $\Alg_{\eenu}(\CCC)$ for
the $\infty$-category of algebra objects over $\mathbf{E}^\otimes_n$ in $\CCC^\otimes$.
We refer to an object of $\Alg_{n}(\CCC)$ as an $\eenu$-algebra in $\CCC$.
If we denote by $\assoc^\otimes$ the associative operad (\cite[4.1.1]{HA}),
there is the standard equivalence $\assoc^\otimes\simeq \eone^\otimes$
of $\infty$-operads. We usually identify $\Alg_{1}(\CCC)$
with the $\infty$-category $\Alg_{\assoc^\otimes}(\CCC)$, that is, the 
$\infty$-category of associative algebras in $\CCC$.
We write $\Alg_{n}^+(\CCC)$ for the $\infty$-category $\Alg_n(\CCC)_{/\uni_{\CCC}}$
of augmeneted objects where $\uni_{\CCC}$ denotes the unit algebra.

\item $\LM^\otimes$: the $\infty$-operad defined in \cite[4.2.1.7]{HA}.
Roughly, an algebra over $\LM^\otimes$ is a pair $(A,M)$
such that $A$ is an unital associative algebra
and $M$ is a left $A$-module.
For a symmetric monoidal $\infty$-category $\CCC^\otimes \to \Gamma$,
we write $\LMod(\CCC^\otimes)$ or $\LMod(\CCC)$
for $\Alg_{\LM^\otimes}(\CCC^\otimes)$.
There is the natural inclusion of $\infty$-operads $\assoc^\otimes\to \LM^\otimes$.
This inclusion determines $\LMod(\CCC)\to \Alg_{\assoc^\otimes}(\CCC)\simeq \Alg_1(\CCC)$
which sends $(A,M)$ to $A$.
For $A\in \Alg_1(\CCC)$, we define $\LMod_A(\CCC)$
to be the fiber of $\LMod(\CCC)\to \Alg_1(\CCC)$ over $A$ in $\Cat$.

\item $\RM^\otimes$: these $\infty$-operads are variants of $\LM^\otimes$
which are used to define structures of right modules over associative algebras
 \cite[4.2.1.36]{HA}.
Informally, an algebra over $\RM^\otimes$ is a pair $(A,M)$
such that $A$ is an unital associative algebra and
$M$ is a right $A$-module.
For a symmetric monoidal $\infty$-category $\CCC^\otimes \to \Gamma$,
we write $\RMod(\CCC^\otimes)$ (or simply $\RMod(\CCC)$)
for $\Alg_{\RM^\otimes}(\CCC^\otimes)$.

\item Unless otherwise stated, $k$ is a base field of characteristic zero.

\end{itemize}

{\it Group actions.}
Let $G$ be a group object in $\SSS$ (see e.g. \cite[7.2.2.1]{HTT} for the notion of group objects).
The main example in this paper is the circle $S^1$.
Let $\mathcal{C}$ be an $\infty$-category.
For an object $C\in \CCC$,
a $G$-action on $C$ means a lift of $C\in \CCC$
to $\Fun(BG, \CCC)$, where $BG$ is the classifying space of $G$.
A $G$-equivariant morphism means a morphism in $\Fun(BG, \CCC)$.
We often identify $\Fun(BG, \CCC)$ as the limit (``$G$-invariants''/homotopy fixed pionts) of the trivial $G$-action on $\CCC$
and write $\CCC^G$ for $\Fun(BG, \CCC)$ (e.g., $\Mod_A^{S^1}=\Fun(BS^1,\Mod_A)$).
We remark that when we regard $G$ as a group object, $\CCC^G$ is not the cotensor with the space
$G$.

In the text, $\GG_{\CCC}^{S^1}$ and $\mathbb{T}_{A/k}[-1]^{S^1}$ are cotensors by $S^1\in \SSS$
in $Lie_A$,
while $\End^L(\HH_\bullet(\CCC/A))^{S^1}$ means the homotopy fixed points.
We hope that these symbols may not be confusing.

\vspace{3mm}

{\it Stable $\infty$-categories.}
We recall some formulations of $\infty$-categories of stable $\infty$-categoires,
see \cite[Section 3.1]{IMA}, \cite[Section 3]{BGT1}, \cite[4.8]{HA} for details.
Let $\ST$ be the $\infty$-category of small
stable idempotent-complete $\infty$-categories
whose morphisms are exact functors.
This $\infty$-category is compactly generated.
Let $\CCC$ be a small
stable idempotent-complete $\infty$-category
and let 
$\Ind(\CCC)$ denote the $\infty$-category of Ind-objects \cite[5.3.5]{HTT}.
Then $\Ind(\CCC)$ is a compactly generated stable $\infty$-category.
The inclusion $\CCC\to \Ind(\CCC)$ identifies
the essential image with the full subcategory $\Ind(\CCC)^\omega$ spanned by
compact obejcts in $\Ind(\CCC)$.
We let $\PR$ denote the $\infty$-category of presentable $\infty$-categories 
such that mapping spaces are spaces of functors which preserve small colimits
(i.e., left adjoint functors) \cite[5.5.3]{HTT}. It has a closed symmetric monoidal structure
whose internal Hom/mapping objects are given by the functor category
$\Fun^{\textup{L}}(-,-)$ of left adjoint functors,
see \cite[4.8.1.15]{HA}.
The $\infty$-category $\ST$ also admits a closed symmetric monoidal structure
whose internal Hom/mapping objects are given by the functor category
$\Fun^{\textup{ex}}(-,-)$ of exact functors.
The stable $\infty$-category of compact spectra is a unit object in $\ST$.
We set $\PRST=\Mod_{\SP^\otimes}(\PR)$,
which can be regarded as the full subcategory of $\PR$
that consists of stable presentable $\infty$-categories.
The $\Ind$-construction
$\CCC\mapsto \Ind(\CCC)$ determines a symmetric monoidal functor
\[
\ST\longrightarrow \PR_{\textup{St}}.
\]

Let $A\in \CAlg(\SP)$.
Let $\Mod^\otimes_A\in \CAlg(\PRST)$ be the symmetric monoidal $\infty$-category of $A$-modules in $\SP$.
Let $\Perf^\otimes_A\in \CAlg(\ST)$ be the symmetric monoidal $\infty$-category of compact $A$-modules in $\SP$.
Let $\Mod^\otimes_{\Perf_A^\otimes}(\ST)$ be
the symmetric monoidal $\infty$-category of $\Perf_A^\otimes$-modules in $\ST$.
This symmetric monoidal $\infty$-category is presentable, and
the tensor product functor preserves small colimits separately in each variable. 
We refer to an object of the underlying $\infty$-category $\Mod_{\Perf_A^\otimes}(\ST)$
as an $A$-linear small stable $\infty$-category.
For ease of notation, put $\ST_A^\otimes=\Mod^\otimes_{\Perf_A^\otimes}(\ST)$
and $\ST_A:=\Mod_{\Perf_A^\otimes}(\ST)$.
We refer the reader to \cite{I} for the description of $\ST_A^\otimes$ by means of spectral
categories.
We write $(\PR_A)^{\otimes}$ for $\Mod_{\Mod_A^\otimes}^\otimes(\PRST)$.
We refer to an object of the underlying $\infty$-category $\PR_A$
as an $A$-linear stable presentable $\infty$-categories.

For $B\in \Alg_{1}(\Mod_A)$, we denote by $\LMod_B(\Mod_A)$ (resp. $\RMod_A(\Mod_A)$)
(or simply $\LMod_B$ and (resp. $\RMod_B$)) the $\infty$-category of left $B$-modules
(resp. right $B$-module spectra) (thare is a canonical equivalence $\LMod_B(\Mod_A)\simeq \LMod_B(\SP)$ induced by the forgetful functor $\Mod_A\to \SP$).

\section{Preliminaries}

In this section, we review some of the theories we will use from the next section.
While
 the prerequisite for Section~\ref{factorizationsection}
 is Section~\ref{hochschildhomologysection} and~\ref{derivedscheme},
Section 5--7 needs also Section~\ref{stablederived} and Section~\ref{formalstack}.

\subsection{Hochschild homology and Hochschild cohomology}
\label{hochschildhomologysection}

Let $\ST_A$ denote $\Mod_{\Perf_A^\otimes}(\ST)$.
Let
\[
\HH_\bullet(-/A):\ST_A \longrightarrow \Mod_A^{S^1}=\Fun(BS^1,\Mod_A)
\]
be the symmetric monoidal functor which carries $\mathcal{C} \in \ST_A$ to the Hochschild homology $A$-module spectrum $\HH_\bullet(\mathcal{C}/A)$
over $A$.
We refer the reader to \cite[Section 6, 6.14]{I} for the construction
of $\HH_\bullet(-/A)$
(in {\it loc. cit.} we use the symbol $\HH_\bullet(\mathcal{C})$ instead of  $\HH_\bullet(\mathcal{C}/A)$).
Let $B\in \Alg_1(\Mod_A)$.
We denote by $\LMod_B=\LMod_{B}(\Mod_A)$ the stable $\infty$-category of left $B$-modules in $\Mod_A$. Let $\textup{Perf}_B$ be the smallest stable idempotent-complete subcategory of $\LMod_B(\Mod_A)$
which contains $B$.
We will write  $\HH_\bullet(B/A)$ for $\HH_\bullet(\Perf_B/A)$.
If $B$ is an $\etwo$-algebra in $\Mod_A$, then $\LMod_B$ has an asscoiative monoidal
structure. 
More precisely, $\LMod_B$ lies in $\Alg_1(\PR_A)$. 
Moreover, $\Perf_B$ is promoted to an object of $\Alg_1(\ST_A)$.
Thus, the symmetric monoidal functor $\HH_\bullet(-/A)$ gives $\HH_\bullet(B/A)$
which belongs to $\Alg_1(\Mod^{S^1}_A)$.

We review the Hochschild cohomology of $\CCC\in \ST_A$.
Let $\Ind(\CCC)$ be the Ind-category that belongs to
$\PR_A$. Moreover, it is compactly generated.
We denote by $\theta_A:\Alg_2(\Mod_A)\to \Alg_1(\PR_A)$
the functor informally given by $B\mapsto \LMod^\otimes_B$.
By definition, the endomorphism algebra object $\mathcal{E}nd_A(\Ind(\CCC))\in \Alg_1(\PR_A)$
endowed with a tautological action on $\Ind(\CCC)$ is a final object of
$\LMod(\PR_A)\times_{\PR_A}\{\Ind(\CCC)\}$.
There exists a final object of $\Alg_2(\Mod_A)\times_{\Alg_1(\PR_A)}\LMod(\PR_A)\times_{\PR_A}\{\Ind(\CCC)\}$.
To see this, note first that
there exists a right adjoint $\Alg_{1}(\PR_A)\to \Alg_{2}(\Mod_A)$ of
$\theta_A$, see \cite[4.8.5.11, 4.8.5.16]{HA}. Namely, there exists a adjoint pair
\[
\xymatrix{
\theta_A:\Alg_{2}(\Mod_A)  \ar@<0.5ex>[r] &    \Alg_{1}(\PR_A):E_A  \ar@<0.5ex>[l]  
}
\]
such that $E_A$ sends $\MM^\otimes$ to the endomorphism algebra of the unit object $\uni_{\MM}$.
Then its final object is given by $E_A(\mathcal{E}nd_A(\Ind(\CCC)))\in \Alg_2(\Mod_A)$
with the left module action of $\LMod_{E_A(\mathcal{E}nd_A(\Ind(\CCC)))}$ on $\Ind(\CCC)$
determined by the counit map $\theta_A(E_A(\mathcal{E}nd_A(\Ind(\CCC))))\to \mathcal{E}nd_A(\Ind(\CCC))$.
The Hochschild cohomology $\HH^\bullet(\CCC/A)$ 
is defined to be $E_A(\mathcal{E}nd_A(\Ind(\CCC)))$.
We refer to  $\HH^\bullet(\CCC/A)$ as the Hochschild cohomology of $\CCC$ (over $A$).

Though we use the word ``homology'' and ``cohomology'',
$\HH_\bullet(\CCC/A)$ and $\HH^\bullet(\CCC/A)$
are not graded modules obtained by passing to (co)homology but spectra (or chain complexes) with algebraic structures.

\subsection{Derived schemes}

\label{derivedscheme}

Let $A$ be a connective commutative dg algebra over $k$.
We fix our convention on derived schemes.
The references are \cite{TV2}, \cite{Gai2}, \cite{DAG}.
We set $\Aff_A=(\CCAlg_A)^{op}$.
We refer to an object of the functor category $\Fun(\CCAlg_A,\wSSS)$
as a derived prestack or (simply) a functor.
The $\infty$-category $\CST_A$ of derived stacks
is defined to be the full subcategory of $\Fun(\CCAlg_A,\wSSS)$
which consists of derived prestacks satisfying the \'etale descent property.
By the \'etale descent we mean the descent condition with respect to Cech nerves of \'etale coverings
\cite[I, Ch. 2, 2.3]{Gai2}. See e.g. \cite{TV2} or \cite[I, Ch. 2,  2.1]{Gai2}
for flat, \'etale morphisms, open immersions, and related properties.
By abuse of notation, we write $\Aff_A$ also for the essential image 
of the Yoneda embedding $\Aff_A\to \Fun(\CCAlg_A,\wSSS)$.
We denote by $\Spec B$ the image of $B\in (\CCAlg_A)^{op}$.
Any representable prestack $\Spec B$ satisfies the descent condition so that
the essential image $\Aff_A$ is contained in $\CST_A$.
If a stack $F$ is representable, we call $F$ a derived affine scheme.
A derived scheme over $A$ is a stack $F$ which satisfies the following conditions:
\begin{itemize}
\item The diagonal morphism $F\to F\times F$ is affine schematic. That is, for any morphism
$\Spec B \to F\times F$, the fiber product $\Spec B\times_{F\times F}F$ is represented by
a derived affine scheme,

\item There exists a set of derived affine schemes $\{T_i\}_{i\in I}$
and morphisms $\{ \phi_i:T_i\to F \}_{i\in I}$ such that (i) each $\phi_i$ is an open immersion
(i.e., for any $\Spec B\to F$, the base change $T_i\times_{F}\Spec B\to \Spec B$ is 
an open immersion of derived affine schemes), and (ii) for any $\Spec B\to F$ the set of base changes
$\{T_i\times_{F}\Spec B \to \Spec B\}_{i \in I}$ is a family of open immersions which
covers $\Spec B$.

\end{itemize}

More precisely, this definition gives what we call derived schemes with affine diagonal.

{\it Mapping stacks.}
Let $X$ be a derived scheme over a derived scheme $Y$.
Let $S$ be a space, that is, an object of $\SSS$.
We define the mapping stack $\Map_Y(S,X)$ (also called the mapping scheme).
It is defined as a functor $\bigl((\Aff)_{/Y}\bigr)^{op} \to \SSS$.
Consider the functor $e_S:(\Aff)_{/Y}\to \Fun(\CAlg,\SSS)_{/Y}$
given by the product with $S$, that is, it carries
$\Spec R\to Y$ to $S\times \Spec R\stackrel{\textup{pr}_2}{\to} \Spec R\to Y$
where $S$ is regarded as the contant functor with value $S$. 
The functor $\Map_Y(S,X):\bigl((\Aff)_{/Y}\bigr)^{op} \to \SSS$ is defined to be the compoite
\[
\bigl((\Aff)_{/Y}\bigr)^{op}\stackrel{e_S^{op}}{\longrightarrow} \bigl(\Fun(\CAlg,\SSS)_{/Y}\bigr)^{op}\stackrel{h_X}{\longrightarrow} \SSS
\]
where $h_X$ is rerpesented by $X$.
Namely, it sends $\Spec R\to Y$ to $\Map_{\Fun(\CAlg,\SSS)_{/Y}}(S\times \Spec R, X)$.

There are two typical cases such that
$\Map_Y(S,X)$ is a derived scheme over $Y$.

(i)
Suppose that $S$ has the homotopy type of a finite CW complex.
Note that $S$ belongs to the full subcategory
which contains the contractible space and is closed under finite colimits.
Thus, $\Map_Y(S,X)$ is contained in the full subcategory of $\Fun(\CAlg,\SSS)_{/Y}$ which containes $X$ and is closed under finite limits.
It follows that $\Map_Y(S,X)$ is a derived scheme.

(ii)
Suppose that $X\to Y$ is affine (but $S$ is arbitrary).
In this case, $\Map_Y(S,X)$ is a derived scheme. The problem
is local on $Y$. Assume that $Y$ is affine.
When $X=\Spec B\to \Spec A=Y$, Then 
$\Map_Y(S,X)=\Spec B\otimes_AS$ where $\otimes_AS$ means
the tensor with $S$ in $\CAlg_A$.

We define the (free) loop space $L_YX$ of $X$ to be $\Map_Y(S^1,X)$.
By definition, it is easy to see that $\Map_Y(S^1, X)\simeq X\times_{X\times_YX}X$. When $Y$ is clear, we write $LX$ for $L_YX$.

\subsection{}
\label{stablederived}
Let $\AAA$ be a small $\infty$-category.
Let $\PPP(\AAA)$ denote the functor category $\Fun(\AAA^{op},\SSS)$ where the $\SSS$ is the $\infty$-category of spaces/$\infty$-groupoids. There
is the Yoneda embedding $\mathfrak{h}_{\AAA}:\AAA\to \PPP(\AAA)$.
Let $\PPP_{\Sigma}(\AAA)\subset \PPP(\AAA)$ be the full subcategory spanned by
those functors $\AAA^{op}\to \SSS$ which preserve finite products \cite[5.5.8.8]{HTT}.
The $\infty$-category  $\PPP_{\Sigma}(\AAA)$ is compactly generated,
and $\PPP_{\Sigma}(\AAA)\subset \PPP(\AAA)$ is characterized as
the smallest full subcategory which contains
the essential image of the Yoneda embedding and is closed under sifted colimits.

Suppose that $\AAA$ admits finite coproducts and a zero object $0$.
We recall the $\infty$-category $\mathcal{P}^{st}_{\Sigma}(\AAA)$
associated to $\AAA$ (cf. \cite{H}, \cite{IMA}).
Consider the set of morphisms $S=\{\mathfrak{h}_{\AAA}(0)\sqcup_{\mathfrak{h}_{\AAA}(M)}\mathfrak{h}
_{\AAA}(0)\to \mathfrak{h}_{\AAA}(0\sqcup_{M}0)\}$
where $M\in \AAA$
such that the pushout $0\sqcup_{M}0$ exists in $\AAA$.
Let $\mathcal{P}^{st}_{\Sigma}(\AAA)$ be
the presentable $\infty$-category obtained from $\mathcal{P}_{\Sigma}(\AAA)$
by inverting morphisms in $S$ (see e.g. \cite[5.5.4]{HTT} for the localization).
The $\infty$-category
$\mathcal{P}^{st}_{\Sigma}(\AAA)$ can be regarded as the full subcategory of $\PPP_\Sigma(\AAA)$
spanned by $S$-local objects. In other words, $\mathcal{P}^{st}_{\Sigma}(\AAA)$ is the full subcategory
which consists of functors $F$ such that the canonical morphism
$F(0\sqcup_M0)\to \ast\times_{F(M)}\ast$ is an equivalence for any $M\in \AAA$ such that
$0\sqcup_M0$ exists in $\AAA$.
Any object of the essential image of the Yoneda embedding
$\AAA\to \PPP_\Sigma(\AAA)$ is $S$-local.

The presentable $\infty$-category $\mathcal{P}^{st}_{\Sigma}(\AAA)$
can be characterized by a universal property. 
Let $\DDD$ be a presentable $\infty$-category.
Let $\Fun^{\textup{L}}(\mathcal{P}^{st}_{\Sigma}(\AAA),\DDD)$ be the full subcategory
of $\Fun(\mathcal{P}^{st}_{\Sigma}(\AAA),\DDD)$ which consists of colimit-preserving functors (i.e.,
left adjoint functors).
Let $\Fun^{st}(\AAA,\DDD)$ be the full subcategory of $\Fun(\AAA,\DDD)$
spanned by those functors $f$ which preserve finite coproducts 
and carry pushouts of the form $0\sqcup_M0$ ($M\in \AAA$)
to $f(0)\sqcup_{f(M)}f(0)$.
Taking into account the universal properties of $\PPP_\Sigma$
and the localization  \cite[5.5.8.15, 5.5.4.20]{HA},
we see that the composition with the fully faithful functor $\AAA\hookrightarrow \PPP^{st}_\Sigma(\AAA)$ induced by the
Yoneda embedding determines an equivalence
\[
\Fun^{\textup{L}}(\PPP^{st}_\Sigma(\AAA),\DDD)\stackrel{\sim}{\to} \Fun^{st}(\AAA,\DDD).
\]

\subsection{Formal stacks}
\label{formalstack}
Let $A$ be a connective commutative dg algebra $A$ over a field $k$ of characteristic zero.
We use a correspondence between pointed formal stacks and dg Lie algebras over $A$, which is proved by Hennion \cite{H} which generalizes the correspondence between dg Lie algebras and formal moduli problems
over a field of characteristic zero (see Lurie \cite[X]{DAG}). 

Let $Lie_A$ be the $\infty$-category of dg Lie algebras.
The $\infty$-category $Lie_A$ is obtained from the model category of dg Lie algebras 
(whose fibrations are termwise surjective maps) by inverting
quasi-isomorphisms
(another equivalent approach is to define it as
the $\infty$-category of algebras over the Lie operad $\Lie$).
Let
$\EXT_A$ be the full subcategory of $\CCAlg_{A//A}:=(\CCAlg_A)_{/A}$,
which is spanned by  trivial square zero extensions $A=A\oplus 0\hookrightarrow A\oplus M\stackrel{p_1}{\to} A$ such that
$M$ is a connective $A$-module of the form $\oplus_{1\le i\le n}A^{\oplus r_i}[d_i]$ ($r_i\ge0,\ d_i\ge0)$.
By abuse of notation, we often write $R$ for an object $A\to R\to A$ of $\CCAlg_{A//A}$.
Similarly, we often omit the augmentations from the notation.
Let $\TSZ_A$ denote the opposite category of $\EXT_A$.
We define the $\infty$-category $\FST_A$ of pointed formal stacks over $A$ to be
$\PPP_\Sigma^{st}(\TSZ)$. We often regard $\TSZ_A$ as a full subcategory of $\FST_A$.
We refer to an object of $\Fun(\EXT_A,\SSS)$ as a pointed formal prestack (or simply a formal stack).
By definition, $\FST_A$ is the full subcategory of $\Fun(\EXT_A,\SSS)$
so that we usually think of a pointed formal stack as a functor $\EXT_A\to \SSS$.
Unfolding the definition, a pointed formal stack over $A$ is nothing but a functor $F:\EXT_A \to \SSS$
which satisfies the following conditions
\begin{itemize}

\item it preserves finite products,

\item $F(A\times_{R}A)\simeq \ast\times_{F(R)}\ast$ for any $R \in \EXT_A$
such that $A\times_RA\in \EXT$. 

\end{itemize}

Let $\Free_{Lie}:\Mod_A\to Lie_A$ be the free Lie algebra functor
which is a left adjoint to the forgetful functor $Lie_A\to \Mod_A$. 
Let $\Mod_{A}^{f}\subset \Mod_A$ be the full subcategory that consists of
objects of the form $\oplus_{1\le i\le n}A^{\oplus r_i}[d_i]$ ($d_i\le -1$).
Let $Lie_A^f$ be the full subcategoy of $Lie_A$, which is the essential image of the restriction of
the free Lie algebra functor  $\Mod_A^{f}\to Lie_A$.
According to \cite[1.2.2]{H}, the inclusions $\Mod_A^f\hookrightarrow \Mod_A$ and 
$Lie_A^f\hookrightarrow Lie_A$ are extended to equivalences
$\PPP_\Sigma^{st}(\Mod_A^f) \stackrel{\sim}{\to} \Mod_A$ and $\PPP_\Sigma^{st}(Lie_A^f) \stackrel{\sim}{\to}Lie_A$ in an essentially unique way (cf. Section~\ref{stablederived}).
Let
\[
\xymatrix{
Ch^\bullet:Lie_A  \ar@<0.5ex>[r] &   (\CAlg_{A//A})^{op}:\DD_\infty   \ar@<0.5ex>[l]
}
\]
be the adjoint pair where the left adjoint $Ch^\bullet$
is the ``Chevalley-Eilenberg cochain functor''
which carries $L\in Lie_A$ to the Chevalley-Eilenberg cochain complex $Ch^\bullet(L)$
(i.e., the $A$-linear dual of Chevalley-Eilenberg chain complex), see e.g. \cite[1.4]{H}, \cite[X, 2.2]{DAG}.
Thanks to \cite[1.5.6]{H}, this adjoint pair induces an adjoint pair
\[
\xymatrix{
\mathcal{F}:Lie_A  \ar@<0.5ex>[r] &   \FST_A :\mathcal{L}.\ar@<0.5ex>[l]  
}
\]
Moreover, if $A$ is noetherian, both $\mathcal{F}$ and $\mathcal{L}$ are categorical equivalences. 
Recall that $A$ is noetherian if $H^0(A)$
is
an ordinary noetherian ring, and $H^n(A)$ is trivial when $n$ is big enough and of finite type over $H^0(A)$
for any $n$.
The left adjoint $\mathcal{F}$ is defined as follows. The restriction of the functor $Ch^\bullet$ to $Lie_A^f$
induces $Lie_A^f\to \TSZ\subset (\CAlg_{A//A})^{op}$ such that $Lie_A^f\stackrel{Ch^\bullet}{\to} \TSZ_A\hookrightarrow \FST_A$ belongs to $\Fun^{st}(Lie_A^f,\FST_A)$.
There exists an essentially unique left adjoint functor $\mathcal{F}:Lie_A\simeq \PPP_\Sigma^{st}(Lie_A^f)\to \FST_A$
which extends $Lie_A^f\to \FST_A$ (cf. Section~\ref{stablederived}).
The right adjoint $\mathcal{L}$ is given by the restriction of 
$\Fun(\EXT_A,\SSS)\to \Fun((Lie_A^f)^{op},\SSS)$ determined by
the composition with $Lie_A^f\to \TSZ_A$.
If $A$ is noetherian, $\mathcal{F}$ and $\mathcal{L}$ are
reduced to a pair of mutually inverse functors $Ch^\bullet:Lie_A^f\simeq \TSZ_A:\DD_\infty$.
For $L\in Lie_A$, we usually write $\mathcal{F}_L$ for the associated formal stack $\mathcal{F}(L)\in \FST_A$.

\vspace{2mm}

{\it Completions.}
Let $F:\CCAlg_k\to \SSS$ be a functor, that is, a derived prestack over $k$.
Let $\Spec A\stackrel{i}{\to} F\stackrel{p}{\to} \Spec A$
be a lift to an object of $\Fun(\CCAlg_k,\SSS)_{\Spec A//\Spec A}$.
We think of it as a derived prestack over $\Spec A$ with a section from $\Spec A$.
It is equivalent to giving a functor $F':(\CCAlg_A)_{/A}\to \SSS$
such that $F'(A\stackrel{id}{\to} A)$ is a contractible space.
We briefly recall the fomal completion of $F$ (see \cite{H} for a general treatment).
Since there exists a canonical fully faithful functor $\EXT_A\subset (\CCAlg_A)_{/A}$,
the composite functor $\EXT_A\hookrightarrow (\CCAlg_A)_{/A}\stackrel{}{\to} \SSS$
gives rise to an object $\widehat{F}$ of $\Fun(\EXT_A,\SSS)$ such that $\widehat{F}(A)=\widehat{F}(A\stackrel{id}{\to} A)$ is a contractible space.
We refer to $\widehat{F}$ as the formal completion of $F$ (along $i:\Spec A\to F$).
Suppose that $F$ is a derived scheme (or, more generally, a derived Artin stack \cite[1.3.3]{TV2}, \cite[I, Ch. 2, Section 4]{Gai2}).
Then the formal completion $\widehat{F}$ belongs to
$\FST_A$. Namely, $\widehat{F}$ is a pointed formal stack over $A$
(see \cite[Lemma 2.2]{H}).
In this case, there exists an essentially unique dg Lie algebra $L$
such that $\mathcal{F}_L\simeq \widehat{F}$.
Suppose furthermore that $F$ is locally of finite presentation over $A$:
this condition implies that the cotangent complex $\LL_{F/A}$ is a perfect complex.
The underlying complex of $L$ is equivalent to the pullback $i^*(\mathbb{T}_{F/A}[-1])\in \Mod_A$
where $\mathbb{T}_{F/A}$ is the tangent complex
that is a dual of $\LL_{F/A}$ (see \cite[1.5.6, 2.2.9]{H}).

\section{Canonical Extensions of factorization homology to Mapping stacks}

\label{factorizationsection0}

In this section, we will give a canonical extension of factorization homology.
This extension provides an extension of factorization homology to the mapping
stack which is canonically defined.
The main results are Theorem~\ref{restrictionfactorization}
and Theorem~\ref{generalizedeasy}. 
In Section~\ref{factorizationsection}, we work with a general setting of $\eenu$-algebras.
Section~\ref{factorizationsection} is much more general than Section~\ref{categoricalextension} in the sense that
it holds for arbitrary (structured) manifolds and arbitrary $\eenu$-algebras.
The idea is based on the local-to-global principle which is manifested in 
the very definition of factorization algebras.
 The argument in Section~\ref{factorizationsection} provides a blueprint for that in Section~\ref{categoricalextension}.
The gluing property is also
important (cf. Lemma~\ref{gluingexample}) in applications.
In Section~\ref{categoricalextension}, we give an elaborated version of Section~\ref{factorizationsection}:
an algebra is replaced with a stable $\infty$-category (and the manifold is the circle).
In Section~\ref{dmodule},
we will apply Example~\ref{geometricex2} (and Theorem~\ref{preeasy}).

\subsection{}
\label{factorizationsection}

We briefly recall the definition of factorization homology. 
We refer the reader to \cite{AF}, \cite{AFT}, \cite{Ho} for the theory of factorization homology
we need (in \cite{HA},
the theory is developed under the name of topological chiral homology).
Let $(\Mfldf_n)^{\otimes}$ be the symmetric monoidal $\infty$-category of framed smooth $n$-manifolds
such that the mapping spaces are the space
of smooth embeddings endowed the data of compatibility of framings
(see e.g. \cite[2.1]{AF}, \cite{AFT} or \cite{Ho} for the detailed account).
The symmetric monoidal structure is given by disjoint union.
We write $\Mfldf_n$ for the underlying $\infty$-category.
Let $\Diskf_n$ be the full subcategory of  $\Mfldf_n$ spanned by those manifolds
which are diffeomorphic to a (possibly empty) finite disjoint union of $\RRR^n$.
Let $p:(\Diskf_n)^{\otimes}\to (\Mfldf_n)^{\otimes}$ denote the fully faithful symmetric monoidal functor.

By default, in this paper $\Mod_R$ means $\Mod_{R}(\SP)$.
However, results in this Section~\ref{factorizationsection} holds in a more general setting.
Let $\mathcal{P}^\otimes$ be a symmetric monoidal presentable $\infty$-category
whose tensor product functor $\otimes: \mathcal{P}\times \mathcal{P}\to \mathcal{P}$
preserves small colimits separately in each variable. 
Let $R$ be a commutative algebra object in $\mathcal{P}$, that is, $R\in \CAlg(\mathcal{P})$.
In this Section~\ref{factorizationsection}, we write $\Mod_R^\otimes$ for
the symmetric monoidal $\infty$-category $\Mod_R^\otimes(\mathcal{P}^\otimes)$.
As usual,  $\Mod_R$ means the underlying $\infty$-category.

Let $\Fun^\otimes((\Diskf_n)^{\otimes},\Mod_R^\otimes)$ be the $\infty$-category
of symmetric monoidal functors.
There exists a canonical equivalence $\Fun^\otimes((\Diskf_n)^{\otimes},\Mod_R^\otimes)\simeq \Alg_n(\Mod_R)$
(obtained by identifying $(\Diskf_n)^{\otimes}$ with a symmetric monoidal envelope of the operad $\eenu$,
see \cite[Section 2.2.4]{HA} for monoidal evelopes).
Consider the following adjoint pair
\[
\xymatrix{
p_!: \Fun^\otimes((\Diskf_n)^{\otimes},\Mod_R^\otimes)  \ar@<0.5ex>[r] &  \Fun^\otimes((\Mfldf_n)^{\otimes},\Mod_R^\otimes)  : p^*   \ar@<0.5ex>[l]  
}
\]
where $p^*$ is determined by composition with $p$.
The left adjoint $p_!$ sends $\beta:(\Diskf_n)^{\otimes}\to \Mod_R^\otimes$
to a symmetric monoidal (operadic) left Kan extension $(\Mfldf_n)^{\otimes}\to \Mod_R^\otimes$ of $\beta$  along $p$.
Let us regard $\beta$ as an $\eenu$-algebra $B$ in $\Mod_R$.
For any $M\in \Mfldf_n$, 
we write $\int_{M}B/R$ for the image of $M$ under $p_!(\beta)$.
We refer to $\int_MB/R$ as the factorization homology
of $M$ in coefficients in $B$ over $\Mod_R$.
The factorization homology $\int_{M}B/R$ can naturally be identified with a colimit of 
$(\Diskf_n)_{/M}\to \Diskf_n\to \Mod_R$
where the second functor
is the {\it underlying functor}
of the symmetric monoidal functor $\beta:(\Diskf_n)^\otimes\to \Mod_R^\otimes$.

The functor $p_!:\Alg_n(\Mod_R)\simeq \Fun^\otimes((\Diskf_n)^{\otimes},\Mod_R^\otimes) \to \Fun^\otimes((\Mfldf_n)^{\otimes},\Mod_R^\otimes)$
given by $B\mapsto \int_{M}B/R$ is extended to a symmetric monoidal functor.
To see this, consider the symmetric monoidal structure on 
$\Fun^\otimes((\Mfldf_n)^{\otimes},\Mod_R^\otimes)$ which is induced by that of $\Mod_R$.
Namely, given two symmetric monoidal functors $F,G:(\Mfldf_n)^{\otimes}\to \Mod_R^\otimes$,
the tensor product $F\otimes G$ is informally defined by $(F\otimes G)(M)=F(M)\otimes_R G(M)$ for any $M\in \Mfldf_n$.
For  the precise construction of this symmetric monoidal structure, we refer to \cite[Construction 7.9]{I}.
The $\infty$-category $\Fun^\otimes((\Diskf_n)^{\otimes},\Mod_R^\otimes)\simeq \Alg_n(\Mod_R)$ 
is promoted to a symmetric monoidal $\infty$-category in a similar way,
and the restriction functor $p^*$ is extended to a symmetric monoidal functor.
We observe:

\begin{Lemma}
\label{symmetricfact}
The functor $p_!$ is a symmetric monoidal functor.
\end{Lemma}

\Proof
By the relative adjoint functor theorem \cite[7.3.2.11]{HA}, the left adjoint $p_!$ 
of the symmetric monoidal functor $p^*$
is a oplax symmetric monoidal functor.
Let $F,G:(\Diskf_n)^{\otimes}\to \Mod_R^\otimes$
be two symmetric monoidal functors.
It will suffice to prove that the canonical morphism $p_!(F\otimes G)\to p_!(F)\otimes p_!(G)$
is an equivalence.
To this end, it is enough to show that the evaluation 
$t:p_!(F\otimes G)(M) \to p_!(F)(M)\otimes p_!(G)(M)$
is an equivalence in $\Mod_R$ for any $M\in \Mfldf_n$.
The evaluation $p_!(F\otimes G)(M)$ can be identified with a colimit of
\[
\phi:(\Diskf_n)_{/M}\to \Diskf_n\stackrel{F\otimes G}{\longrightarrow} \Mod_R,
\]
while $p_!(F)(M)\otimes p_!(G)(M)$ is a colimit of
\[
\psi:(\Diskf_n)_{/M}\times (\Diskf_n)_{/M} \to \Diskf_n\times \Diskf_n \stackrel{F\times G}{\longrightarrow} \Mod_R\times \Mod_R\stackrel{\otimes}{\to} \Mod_R.
\]
The functor $F\otimes G:\Disk_n\to \Mod_R$ is naturally equivalent to the composite $\Diskf_n\stackrel{\textup{diagonal}}{\longrightarrow} \Diskf_n\times \Diskf_n \stackrel{\otimes \circ (F\times G)}{\longrightarrow} \Mod_R$ so that
the composite of $\psi$ and the diagonal map $\delta:(\Diskf_n)_{/M} \to (\Diskf_n)_{/M}\times (\Diskf_n)_{/M}$
is equivalent to $\phi$. The map $t$ is induced by the universality of colimits.
Thus, it is enough to show that $\delta$ is cofinal.
According to \cite[3.22]{AF}, $(\Diskf_n)_{/M}$ is sifted so that $\delta$ is cofinal.
\QED

Let us consider  the composite of symmetric monoidal functors
\[
\int_{M}(-)/R: \Alg_n(\Mod_R)\simeq \Fun^\otimes((\Diskf_n)^{\otimes},\Mod_R^\otimes)  \stackrel{p_!}{\longrightarrow}  \Fun^\otimes((\Mfldf_n)^{\otimes},\Mod_R^\otimes)\stackrel{\textup{ev}_M}{\longrightarrow} \Mod_R^\otimes
\]
where $\textup{ev}_M$ is a symmetric monoidal functor given by the evaluation at $M\in \Mfldf_n$.

We first make the following observation. 
Recall that given $R\in \CAlg(\mathcal{P}^\otimes)$,
we write $\Mod_R$ for $\Mod_{R}(\mathcal{P}^\otimes)$.
We can think of $R$ as an object of $\CAlg(\Alg_n(\Mod_R))\simeq \CAlg(\Mod_R)$ where the
equivalence is induced by the equivalence $\Gamma\otimes \eenu\simeq \Gamma$ of $\infty$-operads
coming from Dunn additivity theorem \cite[5.1.2.2]{HA}
(we can also observe this from the facts that $\CAlg(\Mod_A)$ is a coCartesian symmetric monoidal 
$\infty$-category, and $\eenu$ is a unital $\infty$-operad, see \cite[2.4.3.9]{HA}).
Let $R'\to R$ be a morphism in $\CAlg(\mathcal{P}^\otimes)$.  
The restriction functor $\Mod_R\to \Mod_{R'}$ along $R'\to R$ sends
the unit algebra $R$ to an object of $\CAlg(\Alg_n(\Mod_{R'}))\simeq \CAlg(\Mod_{R'})$ which we denote also by $R$.
Since  $\int_M(-)/R':\Alg_n(\Mod_{R'})\to \Mod_{R'}$ is symmetric monoidal
(cf. Lemma~\ref{symmetricfact}),
$\int_MR/R'$ is a commutative algebra object in $\Mod_{R'}$. 
It gives rise to
\[
\int_M(-)/R': \Mod_{R}(\Alg_n(\Mod_{R'}))\to  \Mod_{\int_MR/R'}(\Mod_{R'}).
\]

\begin{Theorem}
\label{restrictionfactorization}
Let $\Delta^1\times \Delta^1\to \CAlg(\mathcal{P}^\otimes)$
be a diagram which is depicted as
\[
\xymatrix{
k \ar[r]^{s'} \ar[d]_s & B^\circ \ar[d]^g \\
A \ar[r]^f & B.
}
\]
Let $C$ be an object of $\Alg_n(\Mod_A)$.
Then there  exists a canonical
equivalence
\[
\int_M C/k\otimes_{\int_MA/k}\int_MB/B^\circ \stackrel{\sim}{\longrightarrow} \int_MC\otimes_AB/B^\circ
\]
in $\Mod_{B^\circ}$. Moreover, the equivalence is natural with respect to $C\in \Alg_n(\Mod_A)$ and $M\in \Mfldf_n$
(see Construction~\ref{canonicalmapconstruction} for the formulation).
\end{Theorem}

\begin{Remark}
From the observation before Theorem~\ref{restrictionfactorization}, $\int_M C/k$ is a module over $\int_MA/k$.
According to Construction~\ref{canonicalmapconstruction} (ii), there is a canonically defined morphism
$\int_MA/k\to \int_MB/B^\circ$ in $\CAlg(\Mod_k)$.
It can be identified with the canonical
morphism $A\otimes_kM\to B\otimes_{B^\circ}M$.
\end{Remark}

\begin{Construction}
\label{canonicalmapconstruction0}
For a morphism $\phi:R'\to R$ in $\CAlg(\mathcal{P}^\otimes)$,
we let $l_\phi:\Mod_{R'}\rightleftarrows \Mod_R:r_\phi$ denote
the adjunction such that $r_{\phi}$ is given by the restriction along $\phi$, and $l_{\phi}$
is given by the base change $\otimes_{R'}R$.
By abuse of notation we write $l_\phi:\Alg_n(\Mod_{R'})\rightleftarrows \Alg_n(\Mod_R):r_\phi$ for the induced adjunction.
The diagram in the statement of Theorem~\ref{restrictionfactorization}
induces the (noncommutative) diagram
\[
\xymatrix{
\Alg_n(\Mod_A) \ar[r]^{l_f} \ar[rrd]_{r_s} & \Alg_n(\Mod_B) \ar[r]^{r_g} & \Alg_n(\Mod_{B^\circ}) \ar[d]^{r_{s'}} \\
  &   &  \Alg_n(\Mod_k).
}
\]
Since $l_f:\Mod_A\rightleftarrows \Mod_B:r_f$ can be regarded
the relative adjunction of symmetric monoidal $\infty$-categories over
$\Gamma$ in the sense of \cite[7.3.2]{HA},
it determines the adjunction $l_f:\Alg_n(\Mod_A)\rightleftarrows \Alg_n(\Mod_B):r_f$ (cf. \cite[7.3.2.6, 7.3.2.7]{HA}).
If we consider the unit map $\textup{id}\to r_f\circ l_f$ of
the adjunction $(l_f,r_f)$, it gives rise to the natural tranformation
\[
\sigma:r_s\to r_s\circ r_f\circ l_f \simeq r_{f\circ s}\circ l_f\simeq r_{s'}\circ r_g\circ l_f.
\]

Next, we consider module categories of these symmetric monoidal $\infty$-categories.
To this end, we briefly recall $\Mod(\MM^\otimes)$ for a symmetric monoidal $\infty$-categorie $\MM^\otimes$.
Let $\Gamma^+\to \Gamma$ denote the undercategory $\Gamma^+=\Gamma_{\langle1\rangle/}$ with the forgetful functor to $\Gamma$
($\Gamma$ is the (nerve of) pointed finite sets, see Section~\ref{NCsection}).
We say that a map in $\Gamma^+$ is inert if the image in $\Gamma$ is inert (see \cite{HA} for inert maps).
Note that there exists a section $\Gamma\simeq \Gamma_{\langle 0\rangle/}\to \Gamma^+$ induced
by $\langle1\rangle\to \langle0\rangle$. 
Given a symmetric monoidal $\infty$-category $u:\MM^\otimes\to \Gamma$,
$\Mod(\MM^\otimes)$ is the full subcategory of $\Fun_{\Gamma}(\Gamma^+,\MM^\otimes)$
spanned by maps $\Gamma^+\to \MM^\otimes$ which carries inert maps to $u$-coCartesian morphisms.
Composition with the section $\Gamma\to \Gamma^+$ 
determines the forgetful functor $\Mod(\MM^\otimes)\to \CAlg(\MM^\otimes)$.

Note that $l_f$ is (promoted to) a symmetric monoidal functor,
and $r_s$, $r_{s'}$, $r_f$ and $r_g$ are lax symmetric monoidal functors.
Thus, applying $\Mod(-)$ to the above diagram we have the induced
diagram 
\[
\xymatrix{
\Mod(\Alg_n(\Mod_A)) \ar[r]^{\Mod(l_f)} \ar[rrd]_{\Mod(r_s)} & \Mod(\Alg_n(\Mod_B)) \ar[r]^{\Mod(r_g)} & \Mod(\Alg_n(\Mod_{B^\circ})) \ar[d]^{\Mod(r_{s'})} \\
  &   &  \Mod(\Alg_n(\Mod_k))
}
\]
As above, 
$\Mod(l_f):\Mod(\Alg_n(\Mod_A))\rightleftarrows \Mod(\Alg_n(\Mod_B)):\Mod(r_f)$
constitutes an adjunction so that its unit map determines
a natural transformation $\sigma':\Mod(r_s)\to \Mod(r_{s'})\circ \Mod(r_g)\circ \Mod(l_f)$ which extends $\sigma$.
\end{Construction}

\begin{Construction}
\label{canonicalmapconstruction}
We will construct
$\int_MC/k\otimes_{\int_MA/k}\int_M B/B^\circ \to \int_MC\otimes_AB/B^\circ$.

\vspace{1mm}

(i)
We first observe the exchange of 
restriction functors
and
functors which define factorization
homology.
Consider the diagram
\[
\xymatrix{
\Alg_{n}(\Mod_{B^\circ})\simeq \Fun^\otimes((\Diskf_n)^{\otimes},\Mod_{B^\circ}^\otimes)  \ar[r]^(0.6){p_!} \ar@<0.5ex>[d]^{r_{s'}} &   \Fun^\otimes((\Mfldf_n)^{\otimes},\Mod_{B^\circ}^\otimes)  \ar@<0.5ex>[d]^{r_{s'}^M}\\
\Alg_n(\Mod_k)\simeq \Fun^\otimes((\Diskf_n)^{\otimes},\Mod_k^\otimes)  \ar@<0.5ex>[u]^{l_{s'}} \ar[r]^(0.6){q_!} & \Fun^\otimes((\Mfldf_n)^{\otimes},\Mod_k^\otimes) \ar@<0.5ex>[u]^{l_{s'}^M}
}
\]
where $l_{s'}$ and $l^M_{s'}$ are symmetric monoidal functors induced by compositions with
$B^\circ\otimes_k:\Mod^\otimes_k\to \Mod_{B^\circ}^\otimes$, and $p_!$ is a left adjoint functor
of the restriction functor $p^*:\Fun^\otimes((\Mfldf_n)^{\otimes},\Mod_{B^\circ}^\otimes)\to \Fun^\otimes((\Diskf_n)^{\otimes},\Mod_{B^\circ}^\otimes)$.
Likewise, $q_!$ is defined to be a left adjoint of the restriction functor $q^*:\Fun^\otimes((\Mfldf_n)^{\otimes},\Mod_k^\otimes)\to \Fun^\otimes((\Diskf_n)^{\otimes},\Mod_k^\otimes)$.
The arrows are symmetric monoidal functors except that
the down arrows are lax symmetric monoidal functors.
The adjoint pairs $(l_{s'},r_{s'})$, $(p_!,p^*)$, and $(q_!,q^*)$ are relative adjunction over $\Gamma$.

The diagram except down arrows commutes up to canonical homotopy.
To see this, note that $l_{s'}\circ q^*\simeq p^*\circ l_{s'}^M$ and $\textup{id}\simeq q^*\circ q_!$.
Using the unit map $\textup{id}\to q^*\circ q_!$ and the counit map
$p_!\circ p^*\to \textup{id}$ we have the natural
transformation $\xi:p_!\circ l_{s'}\to p_!\circ l_{s'}\circ q^*\circ q_!\simeq p_!\circ p^*\circ l_{s'}^M\circ q_!\to l^M_{s'}\circ q_!$.
Since $B^\circ\otimes_k:\Mod_k\to\Mod_{B^\circ}$ preserves small colimits, we see that
$\xi$ is a natural equivalence.
Consequently,
$p_!\circ l_{s'}\simeq l_{s'}^M\circ q_!$ and $l_{s'}\circ r_{s'}\to \textup{id}$ determine $t: l^M_{s'}\circ q_!\circ r_{s'}\simeq  p_!\circ l_{s'}\circ r_{s'}\to p_!$.

By applying $\Mod(-)$ to
the unit $l_{s'}\circ r_{s'}\to \textup{id}$ of the relative adjunction, it is promoted to
a natural transformation between functors $\Mod(\Alg_{n}(\Mod_{B^\circ}))\to \Mod(\Alg_n(\Mod_{B^\circ}))$.
Thus, $t$ induces the natural transformation
$t':\Mod(l^M_{s'}\circ q_!\circ r_{s'})\to \Mod(p_!)$.
We have
\begin{eqnarray*}
\Delta^1\times \Mod(\Alg_{n}(\Mod_{B^\circ})) \stackrel{T'}{\to} \Mod( \Fun^\otimes((\Mfldf_n)^{\otimes},\Mod_{B^\circ}^\otimes))
\end{eqnarray*}
where $T'$ corresponds to $t'$.

\vspace{1mm}

(ii)
Next, we will combine (i) and Construction~\ref{canonicalmapconstruction0}.
The natural tranformations
$t'$
and $\sigma'$ (see Construction~\ref{canonicalmapconstruction0}) naturally determine the sequence of natural transformations
\begin{eqnarray*}
\Mod(l^M_{s'}) \circ \Mod( q_!) \circ\Mod(r_s) &\to& \Mod(l^M_{s'}) \circ \Mod( q_!) \circ\Mod(r_{s'})\circ \Mod(r_g)\circ \Mod(l_f)   \\
&\to& \Mod(p_!)\circ \Mod(r_g)\circ \Mod(l_f)
\end{eqnarray*}
between functors $\Mod(\Alg_n(\Mod_A))\to \Mod(\Fun^\otimes((\Mfldf_n)^{\otimes},\Mod_{B^\circ}^\otimes)))$.
For ease of notation, we write $U\to V\to W$ for the sequence.
Consider the composition with $\Alg_n(\Mod_A)\simeq \Mod_A(\Alg_n(\Mod_A))\to \Mod(\Alg_n(\Mod_A))$.
We regard $A$ as a unit algebra in $\Alg_{n}(\Mod_{A}))$. 
Let us consider the map $\Delta^1\times\{A\}\to \Mod( \Fun^\otimes((\Mfldf_n)^{\otimes},\Mod_{B^\circ}^\otimes))$
induced by $V\to W$.
The image of $\{1\}\times \{A\}$ is the functor $\int_{(-)}B/B^\circ$
informally given by $M\mapsto \int_MB/B^\circ$, which is regarded as a module
over itself $\int_{(-)}B/B^\circ$.
If we write $\int_{(-)}B/k$ for $q_!(r_{s'}(B))$,
the image of $\{0\}\times \{A\}$ is
\[
l_{s'}^M(q_!(r_{s'}(B)))\simeq B^\circ\otimes_k\bigl(\int_{(-)}B/k\bigr)
\]
which is defined by the formula $M\mapsto B^\circ\otimes_k\bigl(\int_{M}B/k\bigr)$.
The counit $B^\circ\otimes_kB=l_{s'}\circ r_{s'}(B)\to B$ induces
\[
B^\circ \otimes_k\bigl(\int_{M}B/k\bigr)\simeq \int_M B^\circ\otimes_kB/B^\circ \to \int_MB/B^\circ
\]
which corresponds to 
the composition of $\textup{ev}_M:\CAlg( \Fun^\otimes((\Mfldf_n)^{\otimes},\Mod_{B^\circ}^\otimes))\stackrel{}{\to} \CAlg(\Mod_{B^\circ})$.
Next, let us consider
the map $\Delta^1\times\{A\}\to \Mod( \Fun^\otimes((\Mfldf_n)^{\otimes},\Mod_{B^\circ}^\otimes))$
induced by $U\to V$. The image of
$\{0\}\times \{A\}$ is the canoncal module $B^\circ\otimes_k\int_{(-)}A/k$ over $B^\circ\otimes_k\int_{(-)}A/k\in \CAlg(\Fun^\otimes((\Mfldf_n)^{\otimes},\Mod_{B^\circ}^\otimes))$ where $\int_{(-)}A/k$ is defined by $M\mapsto \int_MA/k$. 
The morphism $B^\circ\otimes_k\int_{(-)}A/k\to B^\circ\otimes_k\int_{(-)}B/k$ 
induced by $U\to V$ is determined by $A\to B$ in $\CAlg(\Mod_A)$.
Next, we consider the image of $C\in \Alg_n(\Mod_A)$ under $U\to V\to W$.
As above, by unfolding the definition, it turns out that the image is
\[
B^\circ\otimes_k\bigl(\int_{(-)}C/k\bigr)\to 
B^\circ\otimes_k\bigl(\int_{(-)}C\otimes_AB/k\bigr)\to \int_{(-)}C\otimes_AB/B^\circ
\]
lying over the sequence $B^\circ\otimes_k\bigl(\int_{(-)}A/k\bigr)\to 
B^\circ\otimes_k\bigl(\int_{(-)}B/k\bigr)\to \int_{(-)}B/B^\circ$
in $\CAlg(\Fun^\otimes((\Mfldf_n)^{\otimes},\Mod_{B^\circ}^\otimes))$.

\vspace{1mm}

(iii)
We will define $\chi$ below. 
Let $\overline{\theta}:\Delta^1\times \Alg_n(\Mod_A)\to \Mod(\Fun^\otimes((\Mfldf_n)^{\otimes},\Mod_{B^\circ}^\otimes))$ be the functor determined by
the composite $U\to W$.
Let 
\[
\theta:\Delta^1\times \{A\}\to \Delta^1\times \Alg_n(\Mod_A)\stackrel{\overline{\theta}}{\to} \Mod(\Fun^\otimes((\Mfldf_n)^{\otimes},\Mod_{B^\circ}^\otimes))\stackrel{\textup{forget}}{\to}  \CAlg( \Fun^\otimes((\Mfldf_n)^{\otimes},\Mod_{B^\circ}^\otimes))
\]
be the composite.
Consider the functor
\[
f:\Alg_n(\Mod_{A})\to \mathcal{A}:=\Fun(\Delta^1,\Mod( \Fun^\otimes((\Mfldf_n)^{\otimes},\Mod_{B^\circ}^\otimes)))\times_{\Fun(\Delta^1,\CAlg( \Fun^\otimes((\Mfldf_n)^{\otimes},\Mod_{B^\circ}^\otimes)))}\{\theta\}
\] 
induced by $\overline{\theta}$.
Since $\Fun^\otimes((\Mfldf_n)^{\otimes},\Mod_{B^\circ}^\otimes)$ has sifted colimits
which commutes with those in $\Mod_{B^\circ}$, it follows from \cite[4.5.3.6, 4.5.1.4]{HA}
that 
  $\Mod(\Fun^\otimes((\Mfldf_n)^{\otimes},\Mod_{B^\circ}^\otimes))\to \CAlg(\Fun^\otimes((\Mfldf_n)^{\otimes},\Mod_{B^\circ}^\otimes))$ is a coCartesian fibration. It follows that the evident inclusion
\begin{eqnarray*}
\Mod_{\int_{(-)} B/B^\circ}( \Fun^\otimes((\Mfldf_n)^{\otimes},\Mod_{B^\circ}^\otimes))\ \ \ \ \ \ \ \ \ \ \ \ \ \ \ \  \ \ \ \ \ \ \ \ \ \ \ \ \ \ \ \ \ \ \ \ \ \ \ \ \ \ \ \ \ \ \ \ \ \ \ \ \ \ \ \ \ \ \ \ \ \ \ \ \ \ \ \ \ \ \ \ \ \ \ \ \ \ \ \ \ \ \ \ \ \ \ \ \ \ \ \  \ \\
\to \mathcal{B}:=\Mod( \Fun^\otimes((\Mfldf_n)^{\otimes},\Mod_{B^\circ}^\otimes))\times_{\CAlg( \Fun^\otimes((\Mfldf_n)^{\otimes},\Mod_{B^\circ}^\otimes))} \CAlg(\Fun^\otimes((\Mfldf_n)^{\otimes},\Mod_{B^\circ}^\otimes)))_{/\int B/B^\circ}
\end{eqnarray*} 
admits a left adjoint $L$ (see \cite[Lemma 3.2.6]{AFF}).
If $(P, D\to \int_{(-)}B/B^\circ)$ is an object of $\mathcal{B}$
such that $P$ lies over
$D\in \CAlg( \Fun^\otimes((\Mfldf_n)^{\otimes},\Mod_{B^\circ}^\otimes))$,
$L$ sends it
to $P\otimes_D\int_{(-)}B/B^\circ$.
The composition of $L^{\Delta^1}$, the evident functor $\mathcal{A}\to \Fun(\Delta^1,\mathcal{B})$, and $f$ gives
\[
\chi:\Alg_n(\Mod_{A})\to \Fun(\Delta^1,\Mod_{\int_{(-)} B/B^\circ}( \Fun^\otimes((\Mfldf_n)^{\otimes},\Mod_{B^\circ}^\otimes))).
\] 

\vspace{1mm}

(iv)
Finally, we will observe the evaluation of $\chi$ at $C \in \Alg_n(\Mod_{A})$ and $M\in \Mfldf_n$.
Since $\Mod_{B^\circ}$ has (sifted) colimits whose tensor product functor
$\Mod_{B^\circ}\times \Mod_{B^\circ} \to \Mod_{B^\circ}$ preserves (sifted) colimits in each variable,
it follows from \cite[3.2.3.1 (1), (2)]{HA} that 
$\Fun^\otimes((\Mfldf_n)^{\otimes},\Mod_{B^\circ}^\otimes )$ has sifted colimits 
such that 
for any $M\in \Mfldf_n$, the symmetric monoidal functor $\textup{ev}_M:\Fun^\otimes((\Mfldf_n)^{\otimes},\Mod_{B^\circ}^\otimes )\to \Mod_{B^\circ}^\otimes$ determined by the evaluation at $M$ preserves sifted colimits.
Consequently,
the base change functor 
$\Mod_{B^\circ \otimes_k(\int_{(-)} A/k)}(\Fun((\Mfldf_n)^{\otimes},\Mod_{B^\circ}^\otimes ))\to \Mod_{\int_{(-)} B/B^\circ}(\Fun((\Mfldf_n)^{\otimes},\Mod_{B^\circ}^\otimes ))$ commutes with the base change functor  
$\Mod_{B^\circ\otimes_k(\int_{M}A/k)}(\Mod_{B^\circ}^\otimes )\to \Mod_{\int_MB/B^\circ}(\Mod_{B^\circ}^\otimes)$
under the evaluation functors, up to canonical homotopy.
It follows that $\chi$ carries $C\in \Alg_n(\Mod_{A})$
to a morphism in $\Fun^\otimes((\Mfldf_n)^{\otimes},\Mod_{B^\circ}^\otimes)$
such that for any $M\in \Mfldf_n$, the evaluation at $M$ is
\[
\chi_M(C):B^\circ\otimes_k\bigl(\int_{M}C/k\bigr)\otimes_{B^\circ\otimes_k(\int_{M}A/k)}\int_MB/B^\circ\to \int_MC\otimes_AB/B^\circ
\]
in $\Mod_{B^\circ}$.
The left-hand side is naturally equivalent to $\int_{M}C/k\otimes_{\int_{M}A/k}\int_MB/B^\circ$.

\end{Construction}

{\it Proof of Theorem~\ref{restrictionfactorization}.}
Recall that  $\int_MC\otimes_AB/B^\circ$ is a (sifted) colimit of
$(\Diskf_n)_{/M}\to \Diskf_n\stackrel{C\otimes_AB}{\to} \Mod_{B^\circ}$ where $C\otimes_AB$ denotes the underlying functor,
that is, $\Diskf_n\stackrel{C\otimes_AB}{\to} \Mod_{B^\circ}$ carries $(\RRR^n)^{\sqcup m}$ (i.e., $n$-disks with $m$ connected components)
to the $m$-fold tensor product $(C\otimes_AB)\otimes_{B^\circ}\cdots \otimes_{B^\circ}(C\otimes_AB)$ in $\Mod_{B^\circ}$.
On the other hand, since the geometric realization of a simplicial diagram
(appeared in the bar construction) naturally commutes
with a colimit (appeared as in the definition of factorization homology), 
$B^\circ \otimes_k\bigl(\int_{M}C/k\bigr)\otimes_{B^\circ\otimes_k(\int_{M}A/k)}\int_MB/B^\circ$
is a sifted colimit of the functor 
$(\Diskf_n)_{/M}\to \Diskf_n\stackrel{(B^\circ \otimes_kC)\otimes_{B^\circ\otimes_kA}B}{\longrightarrow}  \Mod_{B^\circ}$
such that the second functor is given on objects by
\begin{eqnarray*}
(\RRR^n)^{\sqcup m} &\mapsto& \bigl((B^\circ\otimes_kC)\otimes_{B^\circ}\cdots \otimes_{B^\circ}(B^\circ\otimes_kC)\bigr) \otimes_{ ((B^\circ\otimes_kA)\otimes_{B^\circ}\cdots \otimes_{B^\circ}(B^\circ\otimes_kA))}\bigl(B\otimes_{B^\circ}\cdots \otimes_{B^\circ}B\bigr) 
\end{eqnarray*}
where 
$B^\circ\otimes_kC$ and $B^\circ\otimes_kA$ are images of $C=r_s(C) \in \Alg_n(\Mod_k)$ and $A=r_s(A) \in \CAlg(\Mod_k)$
under the base change $B^\circ\otimes_k$,
respectively (each $(\bullet)\otimes\cdots \otimes(\bullet)$ means the $m$-fold tensor product).
Note that the factorization homology is based on the local-to-global principle:
the factorization homology functor is defined as a left Kan extension of 
 a functor $\Diskf_n\to \Mod_{B^\circ}$.
That is, by taking the composition with $(\Diskf_n)_{/M}\to \Diskf_n$ 
and the colimit of the composite, the natural transformation
$(B^\circ\otimes_kC)\otimes_{(B^\circ\otimes_kA)}B\to C\otimes_AB$ between functors $\Diskf_n\to \Mod_{B^\circ}$
induces 
$\chi_M(C):B^\circ\otimes_k\bigl(\int_{M}C/k\bigr)\otimes_{B^\circ\otimes_k(\int_{M}A/k)}\int_MB/B^\circ \to \int_MC\otimes_AB/B^\circ$.
Thus, it is enough to show that for any $m\ge0$, the composite of canonical maps
\begin{eqnarray*}
&\ \ & \bigl((B^\circ\otimes_kC)\otimes_{B^\circ}\cdots \otimes_{B^\circ}(B^\circ\otimes_kC)\bigr) \otimes_{ ((B^\circ\otimes_kA)\otimes_{B^\circ}\cdots \otimes_{B^\circ}(B^\circ\otimes_kA))}\bigl(B\otimes_{B^\circ}\cdots \otimes_{B^\circ}B\bigr) \\
&\to& \bigl((B^\circ\otimes_kC)\otimes_{(B^\circ\otimes_kA)}B\bigr)\otimes_{B^\circ} \cdots \otimes_{B^\circ}\bigl((B^\circ\otimes_kC)\otimes_{(B^\circ\otimes_kA)}B\bigr) \\
&\to& \bigl(C\otimes_{A}B\bigr)\otimes_{B^\circ} \cdots \otimes_{B^\circ}\bigl(C\otimes_{A}B\bigr)
\end{eqnarray*}
is an equivalence. These morphisms are obviously equivalences.
\QED

\begin{Remark}
\label{valuefunctorial}
Let $\mathcal{P}^\otimes$ be the symmetric monoidal 
$\infty$-category in Section~\ref{factorizationsection} and the statement of Theorem~\ref{restrictionfactorization}.
Let $\mathcal{Q}^\otimes$ be another symmetric monoidal presentable $\infty$-category
whose tensor product $\otimes :\mathcal{Q}\times \mathcal{Q}\to \mathcal{Q}$
preserves small colimits separately in each variable.
Let $\phi:\mathcal{P}^\otimes\to \mathcal{Q}^\otimes$ be 
a symmetric monoidal functor which preserves small colimits.
Let $\phi(k)$, $\phi(A)$, $\phi(B)$, $\phi(B^\circ)$
denote the induced commutative algebra objects in $\mathcal{Q}^\otimes$.
Let $\phi(C)$ be the induced object in $\Alg_n(\Mod_{\phi(A)})$.
Note that bar contruction (which defines the relative tensor product)
and the factorization homology are defined as sifted colimits of diagrams.
Since $\phi$ preserves sifted colimits, the image of the equivalence under $\phi$  in Theorem~\ref{generalizedeasy}
can naturally be identified with
\[
\int_M \phi(C)/\phi(k)\otimes_{\int_M\phi(A)/\phi(k)}\int_M\phi(B)/\phi(B^\circ) \stackrel{\sim}{\longrightarrow} \int_M\phi(C)\otimes_{\phi(A)}\phi(B)/\phi(B^\circ).
\]
\end{Remark}

\begin{Remark}
We outline a generalization of Theorem~\ref{restrictionfactorization} without proof\footnote{We learn it from Takuo
Matsuoka with the sketch of proofs. But the proof needs some more setups so that here we will not pursue this direction}.
Let $A_1\leftarrow A_0\to A_2$ be a diagram in $\CAlg(\mathcal{P})$.
Let $B_0$ be an object of $\Alg_{n+1}(\Mod_{A_0})$.
For $i=1,2$, let $B_i$ be an object of $\Alg_n(\LMod_{B_0\otimes_{A_0}A_i})$.
Let $M$ be a framed $n$-manifold.
Consider $B_0\in \Alg_{n+1}(\Mod_{A_0})$
as a symmetric monoidal functor $(\Diskf_n)^{\otimes} \to \Alg_1^\otimes(\Mod_{A_0})$.
We define $\int_{M}B_0/A_0$ to be a colimit of $(\Diskf_n)_{/M}\to  \Alg_1(\Mod_{A_0})$.
Then there exists an equivalence 
\[
\int_MB_1/A_1\otimes_{\int_MB_0/A_0}\int_MB_2/A_2 \simeq \int_M (B_1\otimes_{B_0}B_2)/(A_1\otimes_{A_0}A_2)
\]
in $\Mod_{A_1\otimes_{A_0}A_2}$.

\end{Remark}

\begin{Example}
\label{functorcategoryexample}
Let $I$ be a small $\infty$-catergory. Let $k$ be a commutative ring spectrum
and let $\Mod_k:=\Mod_k(\SP)$.
The functor category
$\Fun(I,\Mod_k)$ admits a symmetric monoidal structure arising from that of $\Mod_k$, which is defined by
the coCartesian fibration
\[
\Fun(I,\Mod_k)^\otimes:=\Fun(I,\Mod^\otimes_k)\times_{\Fun(I,\Gamma)}\Gamma\to \Gamma.
\] 
Set $\mathcal{P}^\otimes=\Fun(I,\Mod_k)^\otimes$.
Let $a_I, b_I, b_I^\circ:I\to \CAlg(\Mod_k)$ be objects of 
$\Fun(I,\CAlg(\Mod_k))\simeq \CAlg(\mathcal{P})$.
Let $\textup{cons}(k):I\to \CAlg(\Mod_k)$ be the 
constant functor with value $k$.
Suppose that $\textup{cons}(k)$, $a_I$, $b_I$ and $b_I^\circ$
is extended to a diagram
\[
\xymatrix{
\textup{cons}(k) \ar[r] \ar[d] & b_I^\circ \ar[d] \\
a_I \ar[r] & b_I.
}
\]
Let $c_I$ be an object of $\Alg_n(\Mod_{a_I}(\Fun(I,\Mod_k)))$.
By Theorem~\ref{restrictionfactorization}, we have
\[
\chi_M(c_I):\int_M c_I/\textup{cons}(k)\otimes_{\int_Ma_I/\textup{cons}(k)}\int_Mb_I/b_I^\circ \stackrel{\sim}{\longrightarrow} \int_Mc_I\otimes_{a_I}b_I/b_I^\circ
\]
in $\Mod_{b_I^\circ}(\Fun(I,\Mod_k))$.

Consider the symmetric monoidal functor $\textup{ev}_i:\Fun(I,\Mod_k)^\otimes\to \Mod_k^\otimes$ given by the evaluation at $i\in I$. Then $\textup{ev}_i$ preserves colimits.
Set $A_i=a_I(i)$, $B_i=b_I(i)$, $B_i^\circ=b_I^\circ(i)$.
Let $C_i=c_I(i) \in \Alg_n(\Mod_{A_i})$.
Then by Remark~\ref{valuefunctorial},
the image of $\chi_M(c_I)$ under $\textup{ev}_i$
can be identified with
\[
\int_M C_i/k\otimes_{\int_M A_i/k}\int_MB_i/B_i^\circ \stackrel{\sim}{\longrightarrow} \int_MC_i\otimes_{A_i}B_i/B_i^\circ
\]
in $\Mod_{B^\circ_i}$

To formulate the gluing property (cf. Lemma~\ref{gluingexample}), we introduce the notion of coCartesian over $a_I$.
Let $r$ be a functor $I\to \CAlg(\Mod_k)$.
Let $m$ be a map $I\to \Mod(\Mod_k)$ 
such that $\pi\circ m=r$ where $\pi$ is the projection
$\Mod(\Mod_k)\to \CAlg(\Mod_k)$.
If $m$ sends every morphism in $I$ to a coCartesian morphism in $\Mod(\Mod_k)$ (over $\CAlg(\Mod_k)$),
we say that $m$ is a coCartesian section over $r$.
Let $\Mod_{r}(\Fun(I,\Mod_k))$ be the fiber of $\Mod(\Fun(I,\Mod_k))\to \CAlg(\Fun(I,\Mod_k))$
over $r$.
Unwinding the definition,
we observe that $m$ is equivalent to giving an object of $\Mod_{r}(\Fun(I,\Mod_k))$.
If an object $m'$ of $\Mod_{r}(\Fun(I,\Mod_k))$ corresponds to a coCartesian section $m$ over $r$,
we say that $m'$ is coCartesian over $r$.
Let $\Mod_{r}(\Fun(I,\Mod_k))^{\textup{coCart}}$ denote the full subcategory of coCartesian objects.
For any morphism $A\to A'$ in $\CAlg(\Mod_k)$, the base change functor $\otimes_AA':\Mod_A\to \Mod_{A'}$
is a symmetric monoidal functor which preserves small colimits. It follows that $\Mod_{r}(\Fun(I,\Mod_k))^{\textup{coCart}}$
is closed under tensor products and small colimits in $\Mod_{r}(\Fun(I,\Mod_k))$.
Let $\Mod_{r}^\otimes(\Fun(I,\Mod_k))$ denote $\Mod_{r}(\Fun(I,\Mod_k))$
endowed with the standard symmetric monoidal structure.  
We denote by $\Mod_{r}^\otimes(\Fun(I,\Mod_k))^{\textup{coCart}}$
the symmetric monoidal $\infty$-category whose symmetric monoidal structure is induced by 
$\Mod_{r}^\otimes(\Fun(I,\Mod_k))$.
Let $c$ be an object of $\Alg_n(\Mod_{r}(\Fun(I,\Mod_k)))$.
We say that $c$ is coCartesian over $r$ if the underlying object of $c$ belongs to $\Mod_{r}(\Fun(I,\Mod_k))^{\textup{coCart}}$.
\end{Example}
 
\begin{Lemma}[Gluing property]
\label{gluingexample}
In the setting of Example~\ref{functorcategoryexample},
if $c_I$ is coCartesian over $a_I$, then $\int_{(-)} c_I/\textup{cons}(k)$, defined as an object of $\Mod_{\int_{(-)} a_I/\textup{cons}(k)}(\Fun(I,\Mod_k))$, is coCartesian over $\int_{(-)} a_I/\textup{cons}(k)$.
\end{Lemma} 

\Proof
The assumption means that the canonical morphism
$C_i\otimes_{A_i}A_j\to C_j$ is an equivalence for any morphism
$i\to j$
in $I$.
Thus, it will suffice to prove that the canonical morphism
\[
\int_MC_i/k\otimes_{\int_MA_i/k}\int_MA_j/k\to  \int_M (C_i\otimes_{A_i}A_j)/k
\]
is an equivalence.
This equivalence follows from the special case of Theorem~\ref{restrictionfactorization}
where $k=B^\circ$.
\QED

\begin{Example}
\label{geometricex1}
Let $X$ be a derived scheme over a base connective commutative ring spectrum $k$.
 Let $J$ be the category (poset) of affine open sets in $X$.
Set $I=J^{op}$. Let $a_I:I\to \CAlg(\Mod_k) $ be the functor defined by
$I\ni U_i\simeq \Spec A_i\mapsto A_i\in \CAlg(\Mod_k)$, where $\mathcal{P}^\otimes=\Fun(I,\Mod_k)$ (see Example~\ref{functorcategoryexample}).
We use the notation in Example~\ref{functorcategoryexample}. 
Consider $\Mod_{a_I}(\Fun(I,\Mod_k))^{\textup{coCart}}$.
By \cite[3.3.3.2]{HTT}, 
the canonical symmetric monoidal functor
\[
\Mod_{a_I}(\Fun(I,\Mod_k))^{\textup{coCart}} \to \lim_{U_i\in I}\Mod_{A_i}(\Mod_k)
\]
is an equivalence.
 By the descent theory, the right-hand side is naturally equivalent to the symmetric monoidal stable $\infty$-category
 $\QC^\otimes(X)$ of quasi-coherent complexes on $X$
 (the reader may adopt $\Mod_{a_I}(\Fun(I,\Mod_k))^{\textup{coCart}}$
 as a definition of the category of quasi-coherent complexes/sheaves on $X$).
We consider 
\begin{eqnarray*}
p_!: \Alg_n(\Mod_{a_I}(\Fun(I,\Mod_k))) \to \Fun^\otimes((\Mfldf_n)^{\otimes},\Mod_{a_I}(\Fun(I,\Mod_k)))  
\end{eqnarray*}
(see Example~\ref{functorcategoryexample}).

By the gluing property (Lemma~\ref{gluingexample}),
for $\mathcal{B}\in \Alg_n(\Mod_{a_I}(\Fun(I,\Mod_k)^{\textup{coCart}}))\simeq \Alg_n(\QC(X))$, the image $p_!(\mathcal{B})$
lies in
$\Fun^\otimes((\Mfldf_n)^{\otimes},\Mod_{a_I}(\Fun(I,\Mod_k))^{\textup{coCart}})  \simeq \Fun^\otimes((\Mfldf_n)^{\otimes},\QC(X)^\otimes)$.
Let $\Alg_n(\QC(X)) \to \Fun^\otimes((\Mfldf_n)^{\otimes},\QC(X)^\otimes)$ be
the functor induced by the restriction.
Fix $M\in \Mfldf_n$.
If we write $\int_M(-)/X$ for the composite
$\Alg_n(\QC(X)) \to \Fun^\otimes((\Mfldf_n)^{\otimes},\QC(X)^\otimes) \stackrel{\textup{ev}_M}{\to} \QC(X)$ with the evaluation at $M$,
according to Example~\ref{functorcategoryexample},
 there exists a commutative diagram
\[
\xymatrix{
\Alg_n(\QC(X))  \ar[r]^{\int_M(-)/X} \ar[rd] & \QC(X)  & \\
   &    \Mod_{\int_Ma_I/k}(\Fun(I,\Mod_k))^{\textup{coCart}} \ar[u] \ar@{-}[r]^(0.65){\sim} & \QC(\Map(M,X))
}
\]
in $\Cat$. 
Note that
$\int_MA_I/k$ is 
the composite
$I\stackrel{a_I}{\to} \CAlg(\Mod_k)\stackrel{\otimes M}{\to} \CAlg(\Mod_k)$
which belongs to
$\CAlg(\Fun(I,\Mod_k))\simeq \Fun(I,\CAlg(\Mod_k))$ (cf. Remark~\ref{commutativecase}).
The vertical functor is induced by the base change along $\int_Ma_I/\textup{cons}(k)\to \int_M a_I/a_I\simeq a_I$ determined by
the trivial map $M\to \ast$.
The mapping stack $\Map(M,X)$ is a derived scheme
obtained by gluing affine schemes $\Spec (A_i\otimes M)$ (cf. Section~\ref{derivedscheme}).
Moreover,   
$\Mod_{\int_MA_I/k}(\Fun(I,\Mod_k))^{\textup{coCart}}$
 can be identified with $\QC(\Map(M,X))$, and the vertical functor
 is equivalent to $*$-pullback functor along the morphism $X\to \Map(M,X)$
induced by the constant maps.
Let $\int_M\widetilde{\mathcal{B}/X}$ denote the image of $\int_M\mathcal{B}/k$ in $\QC(\Map(M,X))$.
Then $\int_M\mathcal{B}/X\in \QC(X)$
is ``extended to'' an object $\int_M\widetilde{\mathcal{B}/X}$ in $\QC(\Map(M,X))$.

\end{Example}

To summarize:

\begin{Theorem}
\label{geometricthm1}
Let $\mathcal{B}$ be an object of $\Alg_n(\QC(X))$. Let $M$ be a framed smooth $n$-manifold.
Then $\int_M\mathcal{B}/X \in \QC(X)$ is naturally promoted to
$\int_M\widetilde{\mathcal{B}/X}\in \QC(\Map(M,X))$ endowed with
\[
\iota^* \int_M\widetilde{\mathcal{B}/X}\simeq \int_M\mathcal{B}/X
\]
where $\iota^*:\QC(\Map(M,X))\to \QC(X)$ is the pullback along the morphism $\iota:X\to \Map(M,X)$ induced by constant maps.
Moreover, if $\End(M)$ is the monoid space ($\eone$-algebra in $\SSS$) of endomorphisms of $M$
in $\Mfldf_n$, then $\int_M\mathcal{B}/X\in \Fun(B\End(M), \QC(X))$ is naturally promoted to an object $\int_M\widetilde{B/X}$ of $\Fun(B\End(M),\QC(\Map(M,X))$. Here $B\End(M)$ is thought of as the full subcategory
 of $\Mfldf_n$ spanned by the single object $M$.
\end{Theorem}

\begin{Remark}
\label{commutativecase}
For $A\in \CAlg(\mathcal{P}^\otimes)$,
$\int_{M}A$ 
can be identified with $A\otimes M \in \CAlg(\mathcal{P}^\otimes)$
where $\otimes M$ indicates
the tensor with the underlying sapce of $M$ in $\SSS$.
From \cite[Prop. 5.1]{AF}, there exists a commutative diagram of symmetric monoidal $\infty$-categories
\[
\xymatrix{
\CAlg(\PPP^\otimes) \ar[r]^(0.4)\xi \ar[d]_{\textup{forget}} & \Fun^\otimes((\Mfldf_n)^\otimes,\PPP^\otimes) \\
\Alg_n(\PPP^\otimes) \ar[ur]_{q_!} &  
}
\]
in $\CAlg(\Cat)$. The arrow $\xi$ is the composition of the functor
$\CAlg(\PPP^\otimes) \to \Fun^\otimes((\Mfldf_n)^\otimes, \CAlg(\PPP^\otimes))$
given by the tensor with the underlying spaces of manifolds,
and the symmetric monoidal fogetful functor $\CAlg(\PPP^\otimes)\to\PPP^\otimes$.
Apply $\CAlg(-)$ to this diagram, that is, take the $\infty$-categories of commutative algebra objects.
Since the tensor products $\Gamma\otimes \Gamma$ and $\Gamma\otimes \eenu$ of $\infty$-operads
are naturally equivalent to the commutative $\infty$-operad $\Gamma$
(cf. \cite[3.2.4.5]{HA}, Dunn addtivity theorem \cite[5.1.2.2]{HA}, \cite[5.1.1.5]{HA}),
the induced vertical functor $\CAlg(\CAlg(\PPP^\otimes))\to \CAlg(\Alg_n(\PPP^\otimes))$
is equivalent to the identity functor $\CAlg(\PPP^\otimes)\to \CAlg(\PPP^\otimes)$.
Thus the induced functor $\CAlg(q_!):\CAlg(\PPP^\otimes)\simeq \CAlg(\Alg_n(\PPP^\otimes))\to \CAlg(\Fun^\otimes((\Mfldf_n)^\otimes,\PPP^\otimes))\simeq \Fun^\otimes((\Mfldf_n)^\otimes,\CAlg(\PPP^\otimes))$
is equivalent to $\CAlg(\xi):\CAlg(\PPP^\otimes)\simeq \CAlg(\CAlg(\PPP^\otimes))\to  \CAlg(\Fun^\otimes((\Mfldf_n)^\otimes,\PPP^\otimes))$
induced by $\xi$.
Note that $\int_M A$ is equivalent to the evaluation of $\CAlg(q_!)(A)$ at $M$.
To observe that $\int_MA\simeq A\otimes M$,
it is enough to show that $\CAlg(\xi):\CAlg(\PPP^\otimes)\to \Fun^\otimes((\Mfldf_n)^\otimes,\CAlg(\PPP^\otimes))$ is equivalent to the functor given by the tensor with the underlying spaces of manifolds.
By $\Gamma\otimes \Gamma \simeq \Gamma$, there exist equivalences $\CAlg(\PPP^\otimes)\stackrel{\sim}{\leftarrow} \CAlg(\CAlg(\PPP^\otimes))\stackrel{\sim}{\to} \CAlg(\PPP^\otimes)$
where the first left arrow $\alpha$ forgets the outer commutative algebra structure,
and the second right arrow $\beta$ forgets the inner commutative algebra structure,
such that the $\beta\circ \alpha^{-1}$ is equivalent to the identity functor.
Using equivalences $\alpha$ and $\beta$, we easily deduce that
$\CAlg(\xi)$ is equivalent to the functor given by the tensor with the underlying spaces of manifolds.
\end{Remark}

\subsection{}

\label{categoricalextension}

The main purpose of Section~\ref{categoricalextension} 
is to prove Theorem~\ref{generalizedeasy}, that provides a categorified
version of Example~\ref{functorcategoryexample} when $M$ is the circle.
We start with the setting of Theorem~\ref{generalizedeasy} 
(in a slightly informal way).
Let $k$ be a commutative ring spectrum.
Let $I$ be the nerve of a small category.
Let $a_I:I\to \CAlg(\Mod_k)=\CAlg_k$, $b^\circ_I:I\to \CAlg_k$ and 
$b_I:I\to \CAlg_k$ be functors which we regard as diagrams
of objects in $\CAlg_k$.
Let 
$a_I\to b_I\leftarrow b_I^\circ$ be morphisms in $\Fun(I,\CAlg_k)$.
For each object $i\in I$ the image of $i$ gives rise to the diagram in $\CAlg_k$
\[
\xymatrix{
k \ar[r] \ar[d] & B^\circ_i \ar[d] \\
A_i \ar[r] &  B_i 
}
\]
where $A_i:=a_I(i)$, $B^\circ_i:=b^\circ_I(i)$ and $B_i:=b_I(i)$.

Let $\textup{ST}_{\CAlg_k}:\CAlg_k\to \wCat$ be the functor which carries
the functor which carries $R\in \CAlg_k$ to
$\ST_{R}$ and carries $R_1\to R_2$ to the base change functor
$\ST_{R_1}\to \ST_{R_2}$ (which is a left adjoint to the restriction 
$\ST_{R_2}\to \ST_{R_1}$ along $\Perf^\otimes_{R_1}\to \Perf_{R_2}^\otimes$).
Let $\mathcal{S}t_{a_I}\to I$, $\mathcal{S}t_{b_I}\to I$ and $\mathcal{S}t_{b^\circ_I}\to I$ be coCartesian fibrations which are classified by functors $\textup{ST}_{\CAlg_k}\circ a_I,\ \textup{ST}_{\CAlg_k}\circ b_I,\ \textup{ST}_{\CAlg_k}\circ b_I^\circ:I\to \wCat$, respectively.
See the discussion before Lemma~\ref{strictsection} for the formulation.

Let $\CCC_I:I\to \mathcal{S}t_{a_I}$ be a section of $\mathcal{S}t_{a_I}\to I$.
For $i\in I$, let $\CCC_i\in \ST_{A_i}$ denote the image of $i$.
Informally, $\CCC_I$ can be regarded as  the assigment
$I\ni i \mapsto \CCC_i\in \ST_{A_i}$
and
$\phi_!:\CCC_i\otimes_{A_i}A_j\to \CCC_j$ in $\ST_{A_j}$ associated to a morphism $\phi:i\to j$ together with homotopy coherence data. Here $\CCC_i\otimes_{A_i}A_j$ means
the base change along $a_I(\phi):A_i\to A_j$.
The section $\CCC_I$ can be thought of as an object of the lax limit of $\textup{ST}_{\CAlg_k}\circ a_I$ in the formulation of $(\infty,2)$-categories.
If $\CCC_I$ carries every morphism in $I$ to a coCartesian
 morphism (equivalently, $\phi_!:\CCC_i\otimes_{A_i}A_j\to \CCC_j$ is an equivalence for any $\phi$), we say that $\CCC_I$ is coCartesian over $a_I$.

We regard $\Fun(I,\Fun(BS^1,\Mod_k))$
as a symmetric monoidal $\infty$-category
whose tensor product is given by termiwise tensor products
induced by that of $\Mod_k$ (cf. \cite[2.1.3.4]{HA}).
Let $E$ be a commutative algebra object in $\Fun(I,\Fun(BS^1,\Mod_k))$,
that is, $E\in \CAlg(\Fun(I,\Fun(BS^1,\Mod_k))$).
Let $E_i\in \CAlg(\Fun(BS^1,\Mod_k))$ denote the evaluation of $E$ at $i\in I$.
Let  $M$ be a module object over $E$.
If $M_i$ denotes the evaluation at $i$, and the
canonical morphism $M_i\otimes_{E_i}E_j\to M_j$
is an equivalence for any morphism $i\to j$ in $I$, 
then we say that $M$ is coCartesian over $E$.

We interpret $\HH_\bullet(-/R)$ as a categorified version of $\int_{S^1}(-)/R$ for stable $\infty$-categories
in Section~\ref{factorizationsection} (cf. Section~\ref{hochschildhomologysection}).
Let $\HH_\bullet(\CCC_I\otimes_{a_I}b_I/b_I^\circ)$ be the diagram of Hochschild homology
defined as a functor $I\to \Fun(BS^1,\Mod_k)=\Mod_k^{S^1}$
informally given by $i\mapsto \HH_\bullet(\CCC_i\otimes_{A_i}B_i/B_i^\circ)$.
Let $\HH_\bullet(\CCC_I/k)$ be the diagram of Hochschild homology
defined as a functor $I\to \Mod_k^{S^1}$
informally given by $i\mapsto \HH_\bullet(\CCC_i/k)$.
Both $\HH_\bullet(\CCC_I\otimes_{a_I}b_I/b_I^\circ)$
and $\HH_\bullet(\CCC_I/k)$ belong to $\Fun(I,\Fun(BS^1,\Mod_k))$.
See Construction~\ref{con5} and Construction~\ref{con6} for the formulation.

\begin{Theorem}
\label{generalizedeasy}
Let $\HH_\bullet(a_I/k)$ and $\HH_\bullet(b_I/b_I^\circ)$
be commutative algebra objects in $\Fun(I,\Fun(BS^1,\Mod_k))$, informally defined
by the formulas $i\mapsto \HH_\bullet(A_i/k)$ and $i\mapsto \HH_\bullet(B_i/B_i^\circ)$, respectively. 
See Construction~\ref{con5} for the formulation. 
Let $\HH_\bullet(a_I/k)\to \HH_\bullet(b_I/b_I^\circ)$
be the morphism induced by $a_I\to b_I\leftarrow b_I^{\circ}$.
Then $\HH_\bullet(\CCC_I/k)$ is promoted to a module over $\HH_\bullet(a_I/k)$,
and there exists
a canonical equivalence
\[
\HH_\bullet(\CCC_I/k)\otimes_{\HH_\bullet(a_I/k)}\HH_\bullet(b_I/b_I^\circ)\simeq \HH_\bullet(\CCC_I\otimes_{a_I}b_I/b_I^\circ)
\]
in $\Mod_{\HH_\bullet(b_I/b_I^\circ)}(\Fun(I,\Mod_k^{S^1}))$.
Moreover, $\CCC_I$ is coCartesian, then $\HH_\bullet(\CCC_I/k)$ is
coCartesian over $\HH_\bullet(a_I/k)$.
\end{Theorem}

\begin{Remark}
The equivalence $\HH_\bullet(\CCC_I/k)\otimes_{\HH_\bullet(a_I/k)}\HH_\bullet(b_I/b_I^\circ)\simeq \HH_\bullet(\CCC_I\otimes_{a_I}b_I/b_I^\circ)$
is functorial in $\CCC_I$.
See Remark~\ref{functorialpsi}.
\end{Remark}



The proof of Theorem~\ref{generalizedeasy}
will require several constructions: Construction~\ref{con1},
Construction~\ref{con2}, Construction~\ref{con3}, Construction~\ref{con4},
Construction~\ref{con5}, Construction~\ref{con6}.

\begin{Construction}
\label{con1}

We will formulate $a_I\to b_I\leftarrow b_I^{\circ}$
in terms of symmetric spectra \cite{HSS}.
We use the theory of symmetric spectra
and model categories of symmetric spectra.
We also use the machinery of symmetric spectra and categories enriched over 
symmetric spectra
in Construction~\ref{con2}, Construction~\ref{con3},  Construction~\ref{con4} and Construction~\ref{con5}.
The theory of symmetric spectra provides a nice model of spectra.
While the reader is not necessary to know the full details of definitions
(and we will use only the model categories of them),
the readers who are not familiar with symmetric spectra
are invited to skip to the proof and grasp the rough strategy on the first reading.

\vspace{2mm}

(i)
Let $\SPS$ denote the symmetric monoidal category of symmetric spectra.
According to \cite[Theorem 2.4, Proposition 2,5]{S},
there exists a combinatorial symmetric monoidal model category
structure on $\SPS$ satisfying the monoid axiom, such that
weak equivalences are stable equivalences (see e.g. \cite{S} for the monoid axiom).
We call it the stable $\mathbb{S}$-model structure.
For a model category $M$, we will denote by $M^c$ the full subcategory spanned by
cofibrant objects.
For a category $C$ and a class of morphisms $W$,
there exist an $\infty$-category $C[W^{-1}]$ (possibly having nonsmall mapping spaces) and a functor $C\to C[W^{-1}]$
such that for any $\infty$-category $\DDD$
the induced funtor $\Fun(C[W^{-1}],\DDD)\to \Fun(C,\DDD)$ is fully faithful
and the essential image consists of those functor which carries 
every morphism in $W$ to an equivalence in $\DDD$.
See \cite[1.3.4, 4.1.7, 4.1.8]{HA}
for localizations with respect to weak equivalences.
If $(\SPS)^c[W^{-1}]$ denotes the $\infty$-category obtained from $(\SPS)^c$
by inverting weak equivalences (stable equivalences),
there exists a symmetric monoidal equivalence $(\SPS)^c[W^{-1}]\simeq \SP$
(see \cite{HA}).

Let $\CAlg(\SPS)$ denote the category of commutative symmetric ring spectra,
endowed with the model structure such that
a morphism is
a fibration or weak equivalence if it is a {\it positive} stable $S$-fibration or weak equivalence in $\SPS$ (see \cite[Section 3]{S} or \cite{PS2}
for stable positive model structure: we here employ the positive variant for a technical reason).
Let $\CAlg((\SPS))^c[W^{-1}]$ be the $\infty$-category obtained from $\CAlg((\SPS))^c$ inverting
by weak equivalences.
There exists a canonical equivalence $\CAlg((\SPS))^c[W^{-1}]\simeq \CAlg(\SP)$  (see Remark~\ref{eqrem}).

\vspace{1mm}

(ii)
Let $\RR$ be an object of $\CAlg(\SPS)$.
Let $\SPS(\RR)$ be the symmetric monoidal category of $\RR$-module objects in $\SPS$, which we refer to as $\RR$-module spectra
In virtue of \cite[Theorem 2.6]{S}
there is a combinatorial symmetric monoidal
projective model structure on $\SPS(\RR)$ satisfying the monoid axiom,
in which a morphism is a weak equivalence (resp. a fibration)
if the underlying morphism in $\SPS$ is a stable equivalence
(resp. a fibration with respect to stable $\mathbb{S}$-model structure).
We refer to this model structure as the stable $\RR$-model structure.
Let $\SPS(\RR)^c[W^{-1}]$ denote the $\infty$-category obtained from $\SPS(\RR)$
by inverting weak equivalences (stable equivalences).
Let $R$ denote the image of $\RR$ in $\CAlg(\SPS[W^{-1}])\simeq \CAlg(\SP)$.
By \cite[4.8.2.19, 4.3.3.17]{HA}
there exists a symmetric monoidal equivalence $(\SPS(\RR))^c[W^{-1}]\stackrel{\sim}{\to} \Mod_{R}((\SPS)^c[W^{-1}])\simeq \Mod_R(\SP)$.
Let $\CAlg(\SPS(\RR))$ be the category of commutative algebra objects
in $\SPS(\RR)$, endowed with the model structure such that
a morphism is
a fibration or weak equivalence if it is a positive stable $R$-fibration or weak equivalence in $\SPS(\RR)$ (see \cite[Theorem 3.2]{S} for the details).
If we write $\CAlg(\SPS(\RR))^c[W^{-1}]$ for the $\infty$-category obtained by inverting weak equivalences,
$\CAlg(\SPS(\RR))^c[W^{-1}]\simeq \CAlg(\Mod_R)=\CAlg_R$
(cf. Remark~\ref{eqrem}).

\vspace{1mm}

(iii)
Let $\KK$ be a cofibrant object in $\CAlg(\SPS)$, which represents $k\in \CAlg(\SP)$.

\vspace{1mm}

(iv)
Let 
$\Fun(I,\CAlg(\SPS(\mathbb{K})))$
denote the ordinary functor category where we abuse notation
by regarding $I$ as an ordinary category.
We put it on the injective model structure
in which weak equivalences (resp. cofibrations)
are objectwise weak equivalences (resp. objectwise cofibrations)
(see e.g. \cite[A 2.8.2]{HTT}).
Let $\Fun(I,\CAlg(\SPS(\mathbb{K})))[W^{-1}]$
denote the $\infty$-category obtained by inverting weak equivalences.
By \cite[Theorem 7.9.8]{Cis}, there exist equivalences of $\infty$-categories
\[
\Fun(I,\CAlg(\SPS(\mathbb{K})))[W^{-1}]\simeq \Fun(I,\CAlg(\SPS(\mathbb{K}))[W^{-1}])\simeq \Fun(I,\CAlg(\Mod_k)).
\]
Let $\mathbf{a}_I:I\to  \CAlg(\SPS(\mathbb{K}))$
be a functor
which  represents $a_I$.
We assume that $\mathbf{a}_I$ is cofibrant with respect to the
injective model structure.
Namely, $\mathbf{a}_I(i)$ is cofibrant in 
$\CAlg(\SPS(\mathbb{K}))$ for any $i\in I$.
Similarly, we take a functor
$\mathbf{b}^\circ_I:I\to  \CAlg(\SPS(\mathbb{K}))$
which represents $b^\circ_I$
and is cofibrant with respect to the
injective model structure.

\vspace{1mm}

(v)
We choose $\mathbf{b}_I:I\to  \CAlg(\SPS(\mathbb{K}))$
which represents $b_I$ and is a fibrant object
in $\Fun(I, \CAlg(\SPS(\mathbb{K}))$ with respect to the injective model structure.
We also choose $\mathbf{f}:\mathbf{a}_I\to \mathbf{b}_I$ and $\mathbf{g}:\mathbf{b}^\circ_I\to \mathbf{b}_I$, which represent $a_I\to b_I$ and $b_I^\circ \to b_I$, respectively.
We assume the property that
$\mathbf{b}_I(i)$ is cofibrant as an object of $\CAlg(\SPS(\mathbf{b}_I^\circ(i)))$
 for any $i\in I$.
Here $\CAlg(\SPS(\mathbf{b}_I^\circ(i)))$ is equipped with
the model structure mentioned in (ii). 
To obtain the final (technical) property
we construct $\mathbf{b}_I$ as follows.
Let $\rho:I\to \wCat$ be the functor
which carries $i$ to $\CAlg(\SPS(\mathbf{b}^\circ_I(i)))$
and carries $i\to j$ to the left Quillen functor $\CAlg(\SPS(\mathbf{b}^\circ_I(i)))\to \CAlg(\SPS(\mathbf{b}^\circ_I(j)))$. The left Quillen functor is
induced by the symmetric monodal base change functor $\SPS(\mathbf{b}^\circ_I(i))\to \SPS(\mathbf{b}^\circ_I(j))$ given by $\otimes_{\mathbf{b}^\circ_I(i)}\mathbf{b}^\circ_I(j)$.
Let $\overline{\rho}:\mathcal{CA}lg(\SPS(\mathbf{b_I^\circ}))\to I$ denote a coCartesian fibration
obtained by applying Grothendieck construction to $\rho$.
Let $\textup{Sect}(\mathcal{CA}lg(\SPS(\mathbf{b_I^\circ})))$
denote the ordinary category of sections of $\overline{\rho}$.
According to \cite[Theorem 2.30]{Bar}, there exists 
the injective model structure
such that $s:I\to \mathcal{CA}lg(\SPS(\mathbf{b}_I^\circ))$
is a weak equivalence (resp. cofibration)
if $s(i)$ is a weak equivalence (resp. cofibration) in $\CAlg(\SPS(\mathbf{b}_I^\circ(i)))$
for any $i\in I$.
We denote by $\mathcal{CA}lg(\SPS(\mathbf{b_I^\circ}))^c$
the full subcategory of $\mathcal{CA}lg(\SPS(\mathbf{b_I^\circ}))$
spanned by fiberwise cofibrant objects.
Given $i\in I$ let $W(i)$ be the class of weak equivalences
in in the fiber $\CAlg(\SPS(\mathbf{b}^\circ_I(i)))^c\subset \mathcal{CA}lg(\SPS(\mathbf{b}_I^\circ))^c$ over $i$.
We regard $\mathcal{CA}lg(\SPS(\mathbf{b_I^\circ}))^c$ with the specified morphisms $\cup_{i\in I}W(i)$
as a marked simplicial set 
in the sence of \cite{HTT}.
Let $\mathcal{CA}lg(\SPS(\mathbf{b_I^\circ}))^c[W^{-1}]\to I$
denote the induced morphism obtained from $\mathcal{CA}lg(\SPS(\mathbf{b_I^\circ}))^c\to I$ by inverting $\cup_{i\in I}W(i)$.
(take a decomposition of $\mathcal{CA}lg(\SPS(\mathbf{b_I^\circ}))^c\to I$
into a trivial cofibration followed by a fibration in the model category of marked simplicial sets). 
According to \cite[Proposition 5.3.2]{Maz}, $\mathcal{CA}lg(\SPS(\mathbf{b_I^\circ}))^c[W^{-1}]\to I$ is a coCaretsian fibration
which classfies the functor $I\to \wCat$
informally given by $i\mapsto \CAlg(\SPS(\mathbf{b}^\circ_I(i)))^c[W^{-1}]\simeq \CAlg(\Mod_{B^\circ_i})$
where $\CAlg(\SPS(\mathbf{b}^\circ_I(i)))^c[W^{-1}]$
is the $\infty$-category obtained by inverting weak (stable) equivalences.
We write
$\textup{Sect}(\CAlg(\SPS(\mathbf{b_I^\circ}))[W^{-1}])$
for the $\infty$-category $\Fun(I,\CAlg(\SPS(\mathbf{b_I^\circ}))[W^{-1}])\times_{\Fun(I,I)}\{\textup{id}_I\}$ of sections.
Let
$\textup{Sect}(\mathcal{CA}lg(\SPS(\mathbf{b_I^\circ})))^c[W^{-1}]$
be the $\infty$-category
obtained from $\textup{Sect}(\mathcal{CA}lg(\SPS(\mathbf{b_I^\circ})))^c$ by inverting weak qeuivalences.
By \cite[Proposition 3.42]{Bal} there exists a canonical equivalence
\[
\textup{Sect}(\mathcal{CA}lg(\SPS(\mathbf{b_I^\circ})))^c[W^{-1}]\simeq \textup{Sect}(\mathcal{CA}lg(\SPS(\mathbf{b_I^\circ}))[W^{-1}]).
\]
Let $\beta$ be an object
of $\textup{Sect}(\mathcal{CA}lg(\SPS(\mathbf{b_I^\circ}))[W^{-1}])$
which is determined by $b_I^\circ\to b_I$.
We choose a cofibrant and fibrant section $\widehat{\mathbf{b}}_I:I\to \mathcal{CA}lg(\SPS(\mathbf{b}_I^\circ))$
which represents $\beta$ through the above equivalence.
The forgetful functors
$\CAlg(\SPS(\mathbf{b}_I^\circ(i)))\to \CAlg(\SPS(\KK))$
determine a functor $\textup{Sect}(\mathcal{CA}lg(\SPS(\mathbf{b_I^\circ})))\to \Fun(I,\CAlg(\SPS(\KK)))$ which preserves 
a weak equivalences and fibrations where the the target category
is endowed with the injective model structure (cf. \cite[Proposition 2.32]{Bar}).
We define $\mathbf{b}_I:I\to  \CAlg(\SPS(\mathbb{K}))$
 to be the image of $\widehat{\mathbf{b}}_I$ under the forgetful functor, that has the desired property by the construction.

\end{Construction}

\begin{Remark}
\label{eqrem}
We explain how the equivalence $\CAlg(\SPS(\RR))^c[W^{-1}]\simeq \CAlg(\SPS(\RR)^c[W^{-1}])$
follows from the work of Pavlov and
Scholbach \cite{PS1}, \cite{PS2}.
Let $Q\textup{Comm}\to \textup{Comm}$ be a cofibrant resolution of the commutative operad
(as the standard symmetric operads in simplicial sets, see \cite{PS2}).
We regard $\SPS(\RR)$ as the symmetric monoidal model category satisfying the monoid axiom, which is
endowed the positive stable $R$-model structure defined in \cite[Proposition 3.1]{S}
(or more generally stable positive model structure defined in \cite{PS2}).
We abuse notation by writing $\Alg_{Q\textup{Comm}}(\SPS(\RR))$
for the model category of $Q\textup{Comm}$-algebras in $\SPS(\RR)$
whose weak equivalences and fibrations are detected by underlying weak equivalences and fibrations in $\SPS(\RR)$.
Then by \cite[Proposition 3.5, Theorem 4.6]{PS2},
it induces a left Quillen equivalence 
$\Alg_{Q\textup{Comm}}(\SPS(\RR))\rightleftarrows \CAlg(\SPS(\RR))$
so that the left adjoint gives rise to an equivalence  $\Alg_{Q\textup{Comm}}(\SPS(\RR))^c[W^{-1}]\stackrel{\sim}{\to} \CAlg(\SPS(\RR))^c[W^{-1}]$ of $\infty$-categories. Here $[W^{-1}]$ indicates the localization with respect to
weak equivalences.
By \cite[Theorem 4.9]{PS2} (see also \cite[Thoerem 7.11]{PS1})
$\Alg_{Q\textup{Comm}}(\SPS(\RR))^c[W^{-1}]\simeq \CAlg(\SPS(\RR)^c[W^{-1}])$. Thus, we have $ \CAlg(\SPS(\RR))^c[W^{-1}]\simeq  \Alg_{Q\textup{Comm}}(\SPS(\RR))^c[W^{-1}]\simeq \CAlg(\SPS(\RR)^c[W^{-1}])$. 
\end{Remark}


\begin{Construction}
\label{con2}
We use the theory of spectral categories.
A spectral category is a category enriched over 
the symmetric monoidal category $\SPS$.
For $\RR\in \CAlg(\SPS)$,
we let $\CAT_{\RR}$ denote the category of $\RR$-spectral categories,
that is, categories enriched over the symmetric monoidal category $\SPS(\RR)$
(see e.g. \cite{BGT1}, \cite{BGT2}, \cite{L}, \cite[Section 6]{I}).
Take $\SPS(\RR)$ as the symmetric monoidal model category
endowed with stable $\RR$-model structure \cite[2.6 (2)]{S}.
Applying the existence theorem of a model structure \cite[Theorem 1.1]{Muro}
to $\SPS(\RR)$ we have a combinatorial model structure
on $\CAT_{\RR}$ such that weak equivalences are Dwyer-Kan equivalences
(see e.g. \cite{Muro} for Dwyer-Kan equivalences).
Any morphism $f:\RR_1\to \RR_2$ in $\CAlg(\SPS)$ gives rise to
a Quillen adjunction $f_!:\CAT_{\RR_1}\rightleftarrows \CAT_{\RR_2}:\textup{rest}(f)$ where the right Quillen functor $\textup{rest}(f)$ is induced by the resctriction
functor $\SPS(\RR_2)\to \SPS(\RR_1)$.
$\RR$-spectral categories up to Morita equivalences are regarded 
as models of $R$-linear stable $\infty$-categories where $R$ is the 
image of $\RR$ in $\CAlg(\SP)$ (see \cite[Section 6]{I}). We will construct
several diagrams of spectral categories parameterized by $I$.

\vspace{2mm}

(i)
We write $\wCatone$ for the full subcategory of $\wCat$,
which consists of ordinary categories.
Consider the (pseudo)functor $\CAT_{\CAlg(\SPS)}:\CAlg(\SPS)\to \wCatone$
which carries $\RR$ to $\CAT_{\RR}$
and carries $f:\RR_1\to \RR_2$
to $f_!:\CAT_{\RR_2}\to \CAT_{\RR_1}$
induced by the symmetric monoidal base change functor $\SPS(\RR_1)\to \SPS(\RR_2)$.
Applying the Grothendieck construction to
$\CAT_{\CAlg(\SPS)}$ we have the corresponding
coCartesian fibration $\overline{\CAT}_{\CAlg(\SPS)}:\mathcal{C}at_{\CAlg(\SPS)}\to \CAlg(\SPS)$.
This coCartesian fibration is also a Cartesian fibration
because each morphism $f$ in $\CAlg(\SPS)$
defines a Quillen adjunction $f_!:\CAT_{\RR_1}\rightleftarrows \CAT_{\RR_2}:\textup{rest}(f)$.
Namely, $\overline{\CAT}_{\CAlg(\SPS)}$ is a biCartesian fibration.
We refer the readers to \cite{HP} for generalities of bifibrations arising from diagrams of model categories. 
We define the coCartesian fibration $\mathcal{C}at_{\mathbf{a}_I}\to I$
to be the base change of $\overline{\CAT}_{\CAlg(\SPS)}$ along
$I\stackrel{\mathbf{a}_I}{\to} \CAlg(\SPS(\KK))\stackrel{\textup{forget}}{\longrightarrow} \CAlg(\SPS)$.
We define
the coCartesian fibrations $\mathcal{C}at_{\mathbf{b}_I}\to I$
and $\mathcal{C}at_{\mathbf{b}^\circ_I}\to I$ in a similar way.

Given $\RR$,
$\CAT_{\RR}^{pc}$ is defined to be the full subcategory of $\CAT_{\RR}$
that consists of pointwise cofibrant $\RR$-spectral categories,
that is, spectral categories whose Hom symmetric spectra are 
cofibrant objects in $\SPS(\RR)$.
Note that
each $f_!:\CAT_{\RR_1}\to \CAT_{\RR_2}$ induces $\CAT_{\RR_1}^{pc}\to \CAT^{pc}_{\RR_2}$ because $f_!:\CAT_{\RR_1}\to \CAT_{\RR_2}$ is induced by the left Quillen base change functor $\SPS(\RR_1)\to \SPS(\RR_2)$.
We define $\mathcal{C}at_{\mathbf{a}_I}^{pc}\to I$
to be the coCartesian fibration obtained by restricting
$\mathcal{C}at_{\mathbf{a}_I}$ to the full subcategory spanned by
pointwise cofibrant categories in all fibers.
Namely, the fiber over $i\in I$ is $\CAT_{\mathbf{a}_I(i)}^{pc}$.
We define $\mathcal{C}at_{\mathbf{b}_I}^{pc}\to I$
and $\mathcal{C}at_{\mathbf{b}^{\circ}_I}^{pc}\to I$
in a similar way.

\vspace{1mm}

(ii)
Let $\RR' \in \CAlg(\SPS)$ be an object, and suppose that
$\RR$ is (promoted to) a cofibrant object of $\CAlg(\SPS(\RR'))$.
Here $\CAlg(\SPS(\RR'))$ is endowed with the {\it positive} stable
$\RR$-model structure defined in \cite[Theorem 3.2]{S}.
Then the image of the forgetful functor 
$\CAT_{\RR}^{pc}\to \CAT_{\RR'}$
is contained in $\CAT_{\RR'}^{pc}$. Namely, we have the induced functor
$\CAT_{\RR}^{pc}\to \CAT_{\RR'}^{pc}$.
To see this, 
it is enough to check that the forgetful functor $\SPS(\RR)\to \SPS(\RR')$ carries cofibrant objects
to cofibrant objects in $\SPS(\RR')$.
For this purpose, we recall that cofibrations in $\SPS$
with respect to the stable $\mathbb{S}$-model structure
is the smallest weakly saturated class \cite[A.1.2.2]{HTT} of morphisms
that contains $\{\mathbb{S}\otimes \phi\}_{\phi\in \operatorname{Mon}}$ where $\operatorname{Mon}$ is
the class of monomorphisms of symmetric sequences, and
$\mathbb{S}\otimes \phi$ denotes the morphism of symmetric spectrum
induced by $\phi$, namely, $\mathbb{S}\otimes(-)$ is the left adjoint of
the forgetful functor from $\SPS$ to the category of symmetric sequences,
see \cite{S}. 
The class of cofibrations in $\SPS(\RR)$
with respect to the stable $\RR$-model structure
is the smallest weakly saturated class of morphisms
containing $\{\RR\otimes\phi:=\RR\wedge (\SSSS\otimes \phi)\}_{\phi\in \operatorname{Mon}}$.
Since we assume that $\RR$ is a cofibrant object of
$\CAlg(\SPS(\RR'))$, it follows from \cite[Corollary 4.3]{S} that $\RR$ is cofibrant in $\SPS(\RR')$
with respect to the stable $\RR'$-model structure
(this consequence is the reason why $\CAlg(\SPS(\RR'))$ is equipped with
the model structure defined in \cite[Theorem 3.2]{S}, see Construction~\ref{con1} (ii)).
Thus, morphisms
$\RR\otimes_{\RR'}(\RR'\otimes \phi)$ are cofibrations in $\SPS(\RR')$
with respect to the stable $\RR'$-model structure.
Since $\SPS(\RR)\to \SPS(\RR')$ preserves colimits,
$\SPS(\RR)\to \SPS(\RR')$ preserves cofibrations.

\vspace{1mm}

(iii)
Set $\AAAA_i:=\mathbf{a}_I(i)$, $\BBBB_i:=\mathbf{b}_I(i)$,
and $\BBBB^\circ_i:=\mathbf{b}^\circ_I(i)$ for $i\in I$.
The composition of $\mathbf{f}$ and $\CAT_{\CAlg(\SPS)}:\CAlg(\SPS)\to \wCatone$
determines a map of coCartesian fibrations $\mathcal{C}at_\mathbf{f}^L:\mathcal{C}at_{\mathbf{a}_I}^{pc}\to \mathcal{C}at_{\mathbf{b}_I}^{pc}$
over $I$ (which preserves coCartesian morphisms).
Similarly, 
$\mathbf{g}:\mathbf{b}^\circ_I\to \mathbf{b}_I$ 
induces a map of coCartesian fibrations $\mathcal{C}at_{\mathbf{g}}^L:\mathcal{C}at_{\mathbf{b}^\circ_I}^{pc}\to \mathcal{C}at_{\mathbf{b}_I}^{pc}$
over $I$.
From the observation in (ii) and the assumption that
$\BBBB_i$ is cofibrant in $\CAlg(\SPS(\BBBB_i^\circ))$
we see that
the restriction $\SPS(\BBBB_i)\to \SPS(\BBBB^\circ_i)$
carries cofibrant objects to cofibrant objects.
It follows that each $\CAT_{\BBBB_i^\circ}^{pc}\to \CAT_{\BBBB_i}^{pc}$, which is
induced by the restriction of $\mathcal{C}at_{\mathbf{g}}^L$ to the fiber over $i$,
admits a right adjoint functor $ \CAT_{\BBBB_i}^{pc}\to  \CAT_{\BBBB^\circ_i}^{pc}$. Consequently, 
$\mathcal{C}at_{\mathbf{g}}^L$ admits a relative right adjoint $\mathcal{C}at_{\mathbf{g}}^R:\mathcal{C}at_{\mathbf{b}_I}^{pc}\to \mathcal{C}at_{\mathbf{b}^\circ_I}^{pc}$ over $I$ (See \cite[Section 7.3.2, 7.3.2.6]{HA} for relative adjunctions).

Consider the forgetful functors 
$\CAT_{\AAAA_i}\to \CAT_{\KK}$.
From the observation in (ii) of Construction~\ref{con2} and the assumption that
$\AAAA_i$ is cofibrant in $\CAlg(\SPS(\KK))$
we see that it induces 
$\CAT_{\AAAA_i}^{pc}\to \CAT_{\KK}^{pc}$.
Thus, the forgetful functors 
$\CAT_{\AAAA_i}^{pc}\to \CAT_{\KK}^{pc}$
determine
a map $\pi_{\mathbf{a}_I}:\mathcal{C}at_{\mathbf{a}_I}^{pc}\to \CAT_{\KK}^{pc}\times I$
over $I$ (such that each fiber defines
the forgetful functor $\CAT_{\AAAA_i}^{pc}\to \CAT_{\KK}^{pc}$).
To be precise, we describe $\pi_{\mathbf{a}_I}$
as a relative right adjoint. Let $\textup{cons}(\KK):I\to \CAlg(\SPS)$ be the constant 
functor taking the value $\KK$. There is a  tautological morphism
$\mathbf{t}:\textup{cons}(\KK)\to \mathbf{a}_I$. As above,
it determines a map of coCartesian fibrations $\mathcal{C}at_{\mathbf{t}}^L:\CAT_{\KK}^{pc}\times I\to \mathcal{C}at_{\mathbf{a}_I}^{pc}$
over $I$.
(Notice that the coCartesian fibration
$\CAT_{\KK}^{pc}\times I\to I$ corresponds to $\CAT_{\CAlg(\SPS)}\circ \textup{cons}(\KK)$.)
By
the general exsitence result
\cite[7.3.2.6]{HA} or the fiberwise construction by hand,
it admits a relative right adjoint functor $\pi_{\mathbf{a}_I}:\mathcal{C}at_{\mathbf{a}_I}^{pc}\to \CAT_{\KK}^{pc}\times I$ over $I$ .
Likewise, the forgetful functors 
$\CAT_{\BBBB^\circ_i}^{pc}\to \CAT_{\KK}^{pc}$
determine
a map $\pi_{\mathbf{b}_I^\circ}:\mathcal{C}at_{\mathbf{b}^\circ_I}^{pc}\to \CAT_{\KK}^{pc}\times I$.

Let us consider the diagram over $I$
\[
\xymatrix{
\mathcal{C}at_{\mathbf{a}_I}^{pc} \ar[r]^{\mathcal{C}at_\mathbf{f}^L} \ar[dr]_{\pi_{\mathbf{a}_I}} &  \mathcal{C}at_{\mathbf{b}_I}^{pc} \ar[r]^{\mathcal{C}at_{\mathbf{g}}^R }  &  \mathcal{C}at_{\mathbf{b}^\circ_I}^{pc} \ar[dl]^{\pi_{\mathbf{b}_I^\circ}} \\
& \CAT_{\KK}^{pc}\times I. & 
}
\]
Note that the composite $\textup{cons}(\KK)\to \mathbf{a}_I\stackrel{\mathbf{f}}{\to} \mathbf{b}_I$
equals to $\textup{cons}(\KK)\to \mathbf{b}^\circ_I\stackrel{\mathbf{g}}{\to} \mathbf{b}_I$.
The unit map of the adjunction $\SPS(\AAAA_i)\rightleftarrows \SPS(\BBBB_i)$ gives rise to
a natural transformation
\[
\sigma:\pi_{\mathbf{a}_I}\to \pi_{\mathbf{b}_I^\circ}\circ \mathcal{C}at_{\mathbf{g}}^R \circ \mathcal{C}at_\mathbf{f}^L
\]
over $I$,
which is given on objects by $\CC_i=\CC_i\otimes_{\AAAA_i}\AAAA_i\to \CC_i\otimes_{\AAAA_i}\BBBB_i$ in $\CAT_{\KK}^{pc}$.
(More precisely, if we replace $\mathcal{C}at_\mathbf{f}^L:\mathcal{C}at_{\mathbf{a}_I}^{pc}\to \mathcal{C}at_{\mathbf{b}_I}^{pc}$ by $\mathcal{C}at_{\mathbf{a}_I}\to \mathcal{C}at_{\mathbf{b}_I}$, then it admits a relative right adjoint
$\mathcal{C}at_\mathbf{f}^R$ so that the unit map of $(\mathcal{C}at_\mathbf{f}^L,\mathcal{C}at_\mathbf{f}^R)$ determines $\sigma:\pi_{\mathbf{a}_I}\to \pi_{\mathbf{a}_I}\circ \mathcal{C}at_\mathbf{f}^R\circ \mathcal{C}at_\mathbf{f}^L\simeq \pi_{\mathbf{b}_I^\circ}\circ \mathcal{C}at_{\mathbf{g}}^R \circ \mathcal{C}at_\mathbf{f}^L$.)
\end{Construction}

We briefly
recall the definition of cyclic objects which determines Hochschild homology
(see \cite[Section 6]{I} or the references therein).
Let $\RR$ be a commutative symmetric ring spectrum and let $\CC$ be a pointwise cofibrant
$\RR$-spectral category. Let $\Lambda$ be the cyclic category (cf. \cite{L}).
The cyclic object $\Lambda^{op}\to \SPS(\RR)$ is defined by the formula
\[
\HH(\CC/\RR)_p:=\bigoplus_{(X_0,\ldots,X_p)} \CC(X_{p-1},X_p)\otimes_\mathbb{R}\ldots \otimes_\mathbb{R} \CC(X_0,X_1)\otimes_\mathbb{R} \CC(X_p,X_0)
\]
(the ``pointwise cofibrant'' condition is necessary for the property that the tensor products in the formula compute the derived ones).
The coproduct is taken over the set of
sequences
$(X_0,\ldots,X_p)$ of objects of $\CC$.
See e.g. \cite{I} for the definition of maps $\HH(\CC/\RR)_p\to \HH(\CC/\RR)_q$
(in a nutshell, degeneracy maps are induced by composition maps,
and face maps are induced by unit maps).
Here we use the symbols $\otimes$ and $\oplus$
instead of the smash product $\wedge$ and the wedge sum $\vee$.
The assignment $\CC\mapsto \HH(\CC/\RR)_{\bullet}$ determines
a symmetric monoidal functor
\[
\HH(-/\RR)_{\bullet}:\CAT_{\RR}^{pc}\to \Fun(\Lambda^{op},\SPS(\RR)^c).
\]
Notice that $\HH(-/\RR)_{\bullet}$ is not $\HH_\bullet(-/\RR)$,
but $\HH_{\bullet}(-/\RR)$ is constructed out of $\HH(-/\RR)_{\bullet}$.

\begin{Construction}
\label{con3}
Next, we generalize $\HH(-/\RR)_{\bullet}$: the domain
will be generalized to a diagram of $\CAT_{\RR}^{pc}$.
For this purpose, we first consider
the functor $\SPS_{\CAlg(\SPS)}:\CAlg(\SPS)\to \wCatone$
which carries $\RR$ to $\SPS(\RR)$
and carries $f_!:\RR_1\to \RR_2$ to
the left Quillen (symmetric monoidal) functor $\otimes_{\RR_1}\RR_2:\SPS(\RR_1)\to \SPS(\RR_2)$.
We obtain a coCartesian and Cartesian fibration which we denote by
$\mathcal{S}p^{\Sigma}_{\CAlg(\SPS)} \to \CAlg(\SPS)$ by Grothendieck
construction (cf. \cite[Section 2.2]{HP}).
Let $J$ be a small category.
Set
\[
\mathcal{S}p^{\Sigma,J}_{\CAlg(\SPS)}:=\Fun(J,\mathcal{S}p^{\Sigma}_{\CAlg(\SPS)})\times_{\Fun(J,\CAlg(\SPS))}\CAlg(\SPS)\to \CAlg(\SPS)
\]
which is a coCartesian fibration classified by the functor $\CAlg(\SPS)\to \wCatone$
given by $\RR\mapsto \Fun(J,\SPS(\RR))$.
If we consider $\CAlg(\SPS)\to \wCatone$
given by $\RR\mapsto \SPS(\RR)^c$
instead of $\SPS_{\CAlg(\SPS)}$,
it gives rise to a coCartesian fibration
$\mathcal{S}p^{\Sigma,J,c}_{\CAlg(\SPS)}\to \CAlg(\SPS)$
which is classified by 
the functor $\CAlg(\SPS)\to \wCatone$
given by $\RR\mapsto \Fun(J,\SPS(\RR)^c)$.

Let $\mathbf{c}:I\to \CAlg(\SPS)$
be a  functor and let $\mathcal{C}at_{\mathbf{c}}^{pc}\to I$
be the coCartesian fibration defined in a similar way to
the case of $\mathcal{C}at_{\mathbf{a}_I}^{pc}\to I$ (cf. Construction~\ref{con3} (i)).
Let $\mathcal{S}p^{\Sigma,\Lambda^{op},c}_{\mathbf{c}}\to I$
denote the base change of $\mathcal{S}p^{\Sigma,\Lambda^{op},c}_{\CAlg(\SPS)}\to \CAlg(\SPS)$ along $\mathbf{c}:I\to \CAlg(\SPS)$ (cf. (i)).
We will extend $\HH(-/\RR)_{\bullet}$ to a map of coCartesian fibrations
$\HH(-/\mathbf{c})_\bullet:\mathcal{C}at_{\mathbf{c}}^{pc}\to \mathcal{S}p^{\Sigma,\Lambda^{op},c}_{\mathbf{c}}$ over $I$ as follows.
If we regard  an object of $\mathcal{C}at_{\mathbf{c}}^{pc}$
as a pair of $(\CC, i\in I)$ such that $\CC$ is a pointwise cofibrant 
$\mathbf{c}(i)$-spectral category,
the image of $(\CC, i\in I)$ is defined to be  the pair
$(\HH(\CC/\mathbf{c}(i))_\bullet:\Lambda^{op}\to \SPS(\mathbf{c}(i))^c,\ i\in I)$ 
regarded as an object of $\mathcal{S}p^{\Sigma,\Lambda^{op},pc}_{\mathbf{c}}$.
Let $(\CC_i,i)\to (\CC_j,j)$ be a morphism in 
$\mathcal{C}at_{\mathbf{c}}^{pc}$ that lies over $i\to j$.
The canonical maps
\[
\CC_i(X_{p-1},X_p)\otimes_{\mathbf{c}(i)}\cdots \otimes_{\mathbf{c}(i)} \CC_i(X_p,X_0)\to \CC_j(X_{p-1},X_p)\otimes_{\mathbf{c}(j)}\cdots \otimes_{\mathbf{c}(j)} \CC_j(X_p,X_0)
\]
induce $\HH(\CC_i/\mathbf{c}(i))_\bullet\to  \HH(\CC_j/\mathbf{c}(j))_\bullet$
as a morphism in $\Fun(\Lambda^{op},\SPS(\mathbf{c}(i)))$.
If $(\CC_i,i)\to (\CC_j,j)$ is coCartesian (i.e., $\CC_i\otimes_{\mathbf{c}(i)}\mathbf{c}(j)\simeq \CC_j$),
 the induced map $\HH(\CC_i/\mathbf{c}(i))_\bullet\otimes_{\mathbf{c}(i)}\mathbf{c}(j)\simeq \HH(\CC\otimes_{\mathbf{c}(i)}\mathbf{c}(j)/\mathbf{c}(j))_\bullet$ is an isomorphism where $\HH(\CC/\mathbf{c}(i))_\bullet\otimes_{\mathbf{c}(i)}\mathbf{c}(j)$ indicates
the image of  $\HH(\CC/\mathbf{c}(i))_\bullet$ under 
$\Fun(\Lambda^{op},\SPS(\mathbf{c}(i)))\to \Fun(\Lambda^{op},\SPS(\mathbf{c}(j)))$. Consequently, the assignment
$(\CC,\ i)\mapsto (\HH(\CC/\mathbf{c}(i))_\bullet,\ i)$ is extended to a map of coCartesian
fibrations (which preserves coCartesian morphisms)
in the natural way.
\end{Construction}

\begin{Construction}
\label{con4}
We will construct a natural transformation
$\theta$ that arise from $\HH(-/-)_\bullet$ and forgetful/restriction functors.

Let $l_{\mathbf{b}_I^\circ}^S:\Fun(\Lambda^{op},\SPS(\KK)^c)  \times I\to \mathcal{S}p^{\Sigma,\Lambda^{op},c}_{\mathbf{b}_I^{\circ}} $ be the map 
of coCartesian fibrations induced by base changes
along $\textup{cons}(\KK)\to \mathbf{b}^\circ_I$.
Since the restriction to each fiber
admits a right adjoint given by the forgetful functor $\Fun(\Lambda^{op},\SPS(\BBBB_i)^c)\to \Fun(\Lambda^{op},\SPS(\KK)^c)$
(see Construction~\ref{con2} (ii) and note $\BBBB_i\in \CAlg(\SPS(\KK))^c$),
$l_{\mathbf{b}_I^\circ}^S$ admits a relative right adjoint
functor $\pi_{\mathbf{b}_I^\circ}^S:\mathcal{S}p^{\Sigma,\Lambda^{op},c}_{\mathbf{b}_I^{\circ}} \to \Fun(\Lambda^{op},\SPS(\KK)^c)  \times I$ over $I$,
that is induced by fiberwise forgetful functors.
Next, we consider
the diagram
\[
\xymatrix{
\mathcal{C}at_{\mathbf{b}_I^\circ}^{pc}  \ar[rr]^(0.4){\mathcal{HH}(-/\mathbf{b}^\circ_I)_\bullet} \ar[d]^{\pi_{\mathbf{b}_I^\circ}}& & \mathcal{S}p^{\Sigma,\Lambda^{op},c}_{\mathbf{b}_I^{\circ}}  \ar[d]^{\pi_{\mathbf{b}_I^\circ}^S} \\
\CAT_{\KK}^{pc}\times I\ar[rr]^(0.4){\mathcal{HH}(-/\KK)_\bullet\times I} & & \Fun(\Lambda^{op},\SPS(\KK)^c)  \times I.
}
\]
Note that $\Fun(\Lambda^{op},\SPS(\KK)^c)  \times I\simeq \mathcal{S}p^{\Sigma,\Lambda^{op},c}_{\textup{cons}(\KK)}$ and $\mathcal{HH}(-/\KK)_\bullet\times I\simeq \HH(-/\textup{cons}(\KK))_\bullet$.
There exists a natural transformation
\[
\theta:(\mathcal{HH}(-/\KK)_\bullet\times I) \circ \pi_{\mathbf{b}_I^\circ}\to \pi_{\mathbf{b}_I^\circ}^S\circ \mathcal{HH}(-/\mathbf{b}_I^\circ)_\bullet
\]
informally given by $\mathcal{HH}(\CC/\KK)_\bullet\to \mathcal{HH}(\CC/\BBBB_i^\circ)_\bullet$ defined by canonical maps
\[
\CC(X_{p-1},X_p)\otimes_\mathbb{K}\ldots \otimes_\mathbb{K} \CC(X_0,X_1)\otimes_\mathbb{K} \CC(X_p,X_0)\to 
\CC(X_{p-1},X_p)\otimes_{\BBBB_i^\circ}\ldots \otimes_{\BBBB_i^\circ} \CC(X_0,X_1)\otimes_{\BBBB_i^\circ} \CC(X_p,X_0)
\]
for $(\CC,\ i\in I)\in \mathcal{C}at_{\mathbf{b}_I^\circ}^{pc}$.
Put another way, we here denote by $l_{\mathbf{b}_I^\circ}$ the relative left adjoint (symmetric monoidal) functor
of $\pi_{\mathbf{b}_I^\circ}$.
The functor $l_{\mathbf{b}_I^\circ}$  is induced by 
the base change of Hom spectra induced by $\textup{cons}(\mathbb{K})\to \mathbf{b}_I^\circ$.
Then $\mathcal{HH}(-/\mathbf{b}_I^\circ)_\bullet \circ l_{\mathbf{b}_I^\circ} \simeq   l_{\mathbf{b}_I^\circ}^S \circ (\mathcal{HH}(-/\KK)_\bullet\times I)$.
Using the unit map $\textup{id}\to \pi_{\mathbf{b}_I^\circ}^S  \circ l_{\mathbf{b}_I^\circ}^S$ and the counit map $l_{\mathbf{b}_I^\circ}  \circ \pi_{\mathbf{b}_I^\circ} \to \textup{id}$, we obtain $\theta$ as the composite 
\[
(\mathcal{HH}(-/\KK)_\bullet\times I) \circ \pi_{\mathbf{b}_I^\circ} \to \pi_{\mathbf{b}_I^\circ}^S \circ l_{\mathbf{b}_I^\circ}^S \circ (\mathcal{HH}(-/\KK)_\bullet \times I) \circ \pi_{\mathbf{b}_I^\circ} \simeq \pi_{\mathbf{b}_I^\circ}^S \circ  \mathcal{HH}(-/\mathbf{b}_I^\circ)_\bullet \circ l_{\mathbf{b}_I^\circ}\circ \pi_{\mathbf{b}_I^\circ} \to  \pi_{\mathbf{b}_I^\circ}^S \circ  \mathcal{HH}(-/\mathbf{b}_I^\circ)_\bullet
\]
that is a symmetric monoidal natural transformation.

\end{Construction}

\begin{Construction}
\label{con5}
We will consider categories of sections of coCartesian fibrations
and will construct natural transformations at the level of categories of sections.

 \vspace{2mm}
 
(i)
Let $p:P\to I$ be a coCartesian fibration whose fibers are ordinary categories. 
Let $\Sect(p)$ or $\Sect(P)$ denote the category of sections of $p$.
Suppose that $p:P\to I$ is classified by a functor $d:I\to \wCatone$
and $d$ is extended to
a functor $d':I\to \CAlg(\wCatone)$ where $\CAlg(\wCatone)$
is the category of commutative monoid objects in $\wCatone$, that is, the $(2,1)$-category
of symmetric monoidal ordinary categories.
In this case, $\Sect(p)$ is promoted to a symmetric monoidal category
in which the tensor product is given by the fiberwise tensor product
(the construction is straightforward and is left to the reader).

We consider symmetric monoidal (ordinary) categories
$\Sect(\mathcal{C}at^{pc}_{\mathbf{a}_I})$, 
$\Sect(\mathcal{C}at^{pc}_{\mathbf{b}_I})$,
$\Sect(\mathcal{C}at^{pc}_{\mathbf{b}^\circ_I})$ and 
$\Sect(\CAT^{pc}_{\KK}\times I)\simeq \Fun(I,\CAT^{pc}_{\KK})$.
Note that relative left adjoint functor $l_{\mathbf{b}_I^\circ}:\CAT^{pc}_{\KK}\times I\to \mathcal{C}at^{pc}_{\mathbf{b}^\circ_I}$ induces a symmetric monoidal
functor $\Sect(\CAT^{pc}_{\KK}\times I)\to \Sect(\mathcal{C}at^{pc}_{\mathbf{b}^\circ_I})$. It follows that $\pi_{\mathbf{b}_I^{\circ}}$
determines a lax symmetric monoidal right adjoint functor 
$\Sect(\pi_{\mathbf{b}_I^\circ}):\Sect(\mathcal{C}at^{pc}_{\mathbf{b}^\circ_I})\to \Sect(\CAT^{pc}_{\KK}\times I)$
in the natural way.
Similarly, $\pi_{\mathbf{a}_I}$ and $ \mathcal{C}at_{\mathbf{g}}^R$ induces lax symmetric monoidal functors
$\Sect(\pi_{\mathbf{a}_I}):\Sect(\mathcal{C}at^{pc}_{\mathbf{a}_I})\to \Sect(\CAT^{pc}_{\KK}\times I)$
and $\Sect( \mathcal{C}at_{\mathbf{g}}^R):\Sect(\mathcal{C}at^{pc}_{\mathbf{b}_I})\to\Sect(\mathcal{C}at^{pc}_{\mathbf{b}^\circ_I}) $.
Let $\Sect(\mathcal{C}at_\mathbf{f}^L)$ denote the symmetric monoidal
functor induced by the symmetric monoidal functor $\mathcal{C}at_\mathbf{f}^L$.
The natural transformation $\sigma$ (cf. Construction~\ref{con2} (iii)) naturally induces a symmetric monoidal
natural transformation
\[
\Sect(\sigma):\Sect(\pi_{\mathbf{a}_I})\to \Sect(\pi_{\mathbf{b}_I^\circ})\circ \Sect(\mathcal{C}at_{\mathbf{g}}^R) \circ \Sect(\mathcal{C}at_\mathbf{f}^L).
\]
We apply the same procedure of taking $\Sect(-)$ to the diagram in Construction~\ref{con4} and the
natural transformation $\theta$. We then obtain the induced symmetric monoidal natural transformation
\[
\Sect(\theta):\Fun(I,\HH(-/\KK)_\bullet) \circ \Sect(\pi_{\mathbf{b}_I^\circ}) \to \Sect(\pi_{\mathbf{b}_I^\circ}^S) \circ \Sect(\mathcal{HH}(-/\mathbf{b}_I^\circ)_\bullet)
\]
(note that $\Sect(\mathcal{HH}(-/\KK)_\bullet\times I)=\Fun(I,\HH(-/\KK)_\bullet)$).
Consider the diagram
\[
\xymatrix{
\Sect(\mathcal{C}at_{\mathbf{a}_I}^{pc}) \ar[r] \ar[rd]_{\Sect(\pi_{\mathbf{a}_I})} & \Sect(\mathcal{C}at_{\mathbf{b}_I^\circ}^{pc})  \ar[rr]^(0.4){\Sect(\mathcal{HH}(-/\mathbf{b}^\circ_I)_\bullet)} \ar[d]^{\Sect(\pi_{\mathbf{b}_I^\circ})}& & \Sect(\mathcal{S}p^{\Sigma,\Lambda^{op},c}_{\mathbf{b}_I^{\circ}})  \ar[d]^{\Sect(\pi_{\mathbf{b}_I^\circ}^S)} \\
  &  \Fun(I,\CAT_{\KK}^{pc})\ar[rr]^(0.4){\Fun(I,\mathcal{HH}(-/\KK)_\bullet)} & & \Fun(I,\Fun(\Lambda^{op},\SPS(\KK)^c)).
}
\]
The composition of $\Sect(\sigma)$ and $\Sect(\theta)$ gives rise to a symmetric monoidal natural transformation 
\[
\tau:\Fun(I,\mathcal{HH}(-/\KK)_\bullet)\circ \Sect(\pi_{\mathbf{a}_I})\to\Sect(\pi_{\mathbf{b}_I^\circ}^S) \circ \Sect(\mathcal{HH}(-/\mathbf{b}_I^\circ)_\bullet)\circ \Sect(\mathcal{C}at_{\mathbf{g}}^R) \circ \Sect(\mathcal{C}at_\mathbf{f}^L).
\]

\vspace{1mm}

(ii)
For $\RR\in \CAlg(\SPS)$ we let $B\RR$ denote the $\RR$-spectral category having one object $\ast$
such that the Hom spectrum is defined by $B\RR(\ast,\ast)=\RR$.
Namely, it is a unit object of $\CAT_{\RR}$.
Let $B\mathbf{a}_I$ denote the section $I\to \mathcal{C}at_{\mathbf{a}_I}^{pc}$
which carries $i\in I$ to $B\AAAA_i$.
In other words,  $B\mathbf{a}_I$ is a unit object of the symmetric monoidal category $\Sect(\mathcal{C}at_{\mathbf{a}_I}^{pc})$.
For ease of notation we abuse notation by writing
$\HH(\mathbf{a}_I/\KK)_\bullet\to \HH(\mathbf{b}_I/\mathbf{b}_I^{\circ})_\bullet$
for the morphism of commutative algebra objects in $\Fun(I,\Fun(\Lambda^{op},\SPS(\KK)^c))$, which is the determined by $B\mathbf{a}_I$ and $\tau$.
(Here we think of $B\mathbf{a}_I$ as a commutative algebra object in 
 $\Sect(\mathcal{C}at_{\mathbf{a}_I}^{pc})$.)
 In particular, $\HH(\mathbf{a}_I/\KK)_\bullet$ is the image of $B\mathbf{a}_I$
 under $\Fun(I,\mathcal{HH}(-/\KK)_\bullet)\circ \Sect(\pi_{\mathbf{a}_I})$,
 and $ \HH(\mathbf{b}_I/\mathbf{b}_I^{\circ})_\bullet$ is the image of $B\mathbf{a}_I$ under $\Sect(\pi_{\mathbf{b}_I^\circ}^S) \circ \Sect(\mathcal{HH}(-/\mathbf{b}_I^\circ)_\bullet)\circ \Sect(\mathcal{C}at_{\mathbf{g}}^R) \circ \Sect(\mathcal{C}at_\mathbf{f}^L)$.
We believe that this notation is not brutal because its evaluation at $i\in I$
is the canonical morphism $\HH(\AAAA_i/\KK)_\bullet\to \HH(\BBBB_i/\BBBB_i^\circ)_\bullet$.

The functor $\Sect(\pi_{\mathbf{b}_I^\circ}^S) \circ \Sect(\mathcal{HH}(-/\mathbf{b}_I^\circ)_\bullet)\circ \Sect(\mathcal{C}at_{\mathbf{g}}^R) \circ \Sect(\mathcal{C}at_\mathbf{f}^L)$ is naturally promoted 
to $\alpha:\Sect(\mathcal{C}at_{\mathbf{a}_I}^{pc})    \to \Mod_{\HH(\mathbf{b}_I/\mathbf{b}_I^{\circ})_\bullet}(\Fun(I,\Fun(\Lambda^{op},\SPS(\KK)^c)))$.
Likewise, $\Fun(I,\mathcal{HH}(-/\KK)_\bullet)\circ \Sect(\pi_{\mathbf{a}_I})$ is naturally promoted 
to $\beta:\Sect(\mathcal{C}at_{\mathbf{a}_I}^{pc})    \to \Mod_{\HH(\mathbf{a}_I/\KK)_\bullet}(\Fun(I,\Fun(\Lambda^{op},\SPS(\KK)^c)))$.
Let $r$ denote the restriction functor along $\HH(\mathbf{a}_I/\KK)_\bullet\to \HH(\mathbf{b}_I/\mathbf{b}_I^{\circ})_\bullet$:
\[
 \Mod_{\HH(\mathbf{b}_I/\mathbf{b}_I^{\circ})_\bullet}(\Fun(I,\Fun(\Lambda^{op},\SPS(\KK)^c)))\to  \Mod_{\HH(\mathbf{a}_I/\KK)_\bullet}(\Fun(I,\Fun(\Lambda^{op},\SPS(\KK)^c))).
\]
Since $\tau$ is a symmetric monoidal natural transformation between lax symmetric monoidal functors, 
there is a natural transformation $\tau':\beta\to r\circ \alpha$ that extends $\tau$.

\vspace{2mm}

(iii)
We consider symmetric monoidal functors
\[
\CAT_{\KK}^{pc}\stackrel{\mathcal{HH}(-/\KK)_\bullet}{\longrightarrow} \Fun(\Lambda^{op},\SPS(\KK)^c)\to \Fun(\Lambda^{op},\SPS(\KK)^c[W^{-1}])\stackrel{L}{\to} \Fun(BS^1,\SPS(\KK)^c[W^{-1}])\simeq \Mod_k^{S^1}.
\]
There is a canonical functor $\SPS(\mathbb{K})^c\to \SPS(\mathbb{K})^c[W^{-1}]$
that induces the second functor. The third functor is the symmetric monoidal functor
determined by left Kan extensions along the groupoid completion $\Lambda^{op}\to BS^1$. 
The composite carries Morita equivalences in $\CAT_{\KK}^{pc}$ to equivalences in $\Mod_k^{S^1}$
(cf. \cite[Lemma 6.11]{I}). 
The composite induces a symmetric monoidal functor $\HH_\bullet(-/k):\CAT_{\KK}^{pc}[M^{-1}]\simeq \ST_k\to \Mod_k^{S^1}$
where $\CAT_{\KK}^{pc}[M^{-1}]$ is the $\infty$-category obtained by inverting Morita equivalences,
and $\CAT_{\KK}^{pc}[M^{-1}]\simeq \ST_k$ is a symmetric monoidal equivalence (see \cite[Proposition 6.7]{I}).

The diagram in Construction~\ref{con5} (i) and the above sequence of symmetric monoidal functors
give rise to
\[
\xymatrix{
\Sect(\mathcal{C}at_{\mathbf{a}_I}^{pc}) \ar[r] \ar[rd] &  \Sect(\mathcal{C}at_{\mathbf{b}^\circ_I}^{pc})  \ar[r] &  \Sect(\mathcal{S}p^{\Sigma,\Lambda^{op},c}_{\mathbf{b}_I^{\circ}})   \ar[d] &  \Fun(I,\Fun(BS^1,\Mod_k)).   \\
   &  \Fun(I,\CAT_{\KK}^{pc}) \ar[r] & \Fun(I,\Fun(\Lambda^{op},\SPS(\KK)^c))  \ar[r] & \Fun(I,\Fun(\Lambda^{op},\SPS(\KK)^c[W^{-1}])) \ar[u] 
}
\]
where $ \Fun(I,\CAT_{\KK}^{pc})\to \cdots\to \Fun(I,\Fun(BS^1,\Mod_k^{S^1}))$
is obtained by applying $\Fun(I,-)$.

We note that the image of $\HH(\mathbf{a}_I/\KK)_\bullet\to \HH(\mathbf{b}_I/\mathbf{b}_I^{\circ})_\bullet$ (see Construction~\ref{con5} (ii)) 
in $\Fun(I,\Fun(BS^1,\Mod_k^{S^1}))$ is the canonical map $\HH_\bullet(\mathbf{a}_I/k)\to \HH_\bullet(b_I/b_I^\circ)$
where $\HH_\bullet(a_I/k)$ and $\HH_\bullet(b_I/b_I^\circ)$ mean the functor informally given by $i\mapsto \HH_\bullet(A_i/k)$
and  $i\mapsto \HH_\bullet(B_i/B_i^\circ)$, respectively
(keep in mind that Hochschild homology is equivalent to a  colimit of a simplicial diagram \cite[Lemma 6.9 (ii)]{I} and the restriction functor $\Mod_{B_i}\to \Mod_k$ preserves colimits).
We easily observe that the morphism $\HH_\bullet(a_I/k)\to \HH_\bullet(b_I/b_I^\circ)$
does not depend on a choice of $\mathbf{a}_I\to \mathbf{b}_I\leftarrow \mathbf{b}_I^\circ$ up to equivalence.
We also remark that by \cite[Lemma 3.5]{IMA}, there exist 
canonical equivalences $\HH_\bullet(A_i/k)\simeq A_i\otimes_kS^1$
and $\HH_\bullet(B_i/B_i^\circ)\simeq B_i\otimes_{B_i^\circ}S^1$ in $\Fun(BS^1,\CAlg_k)$ so that $\HH_\bullet(a_I/k)$ and $\HH_\bullet(b_I/b_I^\circ)$ are naturally equivalent to the functors informally given by $i\mapsto A_i\otimes_kS^1$
and  $i\mapsto B_i\otimes_{B_i^\circ}S^1$, respectively.

Let $\upsilon:\gamma\circ \beta\to \gamma \circ r \circ \alpha$ be the symmetric monoidal natural transformation
between functors
\[
\Sect(\mathcal{C}at_{\mathbf{a}_I}^{pc})\to \Mod_{\HH_\bullet(a_I/k)}(\Fun(I,\Fun(BS^1,\Mod_k)))
\] 
obtained by the composition
of $\tau'$ (see (ii)) and $\gamma:\Mod_{\HH(\mathbf{a}_I/\KK)_\bullet}(\Fun(I,\Fun(\Lambda^{op},\SPS(\KK)^c))) \to \Mod_{\HH_\bullet(a_I/k)}(\Fun(I,\Fun(BS^1,\Mod_k)))$.

\end{Construction}

By a termwise Morita equivalence in $\Sect(\mathcal{C}at_{\mathbf{a}_I}^{pc})$
we mean a morphism $s_1\to s_2$
between two sections $s_1,s_2:I\rightrightarrows \mathcal{C}at_{\mathbf{a}_I}^{pc}$
such that $s_1(i)\to s_2(i)$ is a Morita equivalence in the fiber for any $i\in I$.
Let $\Sect(\mathcal{C}at_{\mathbf{a}_I}^{pc})[M^{-1}]$ denote
the $\infty$-category obtained by inverting termwise Morita equivalences.

Let us consider the $\infty$-categorical description of $\Sect(\mathcal{C}at_{\mathbf{a}_I}^{pc})[M^{-1}]$.
We identify $\CAT_{\AAAA_i}^{pc}$ with the fiber of
$\mathcal{C}at_{\mathbf{a}_I}^{pc}\to I$
over $i\in I$,
we let $M(i)$ denote the class of Morita equivalences in the fiber over $i\in I$.
We write $\mathcal{C}at_{\mathbf{a}_I}^{pc}[M^{-1}]$
for the $\infty$-category obtained by inverting $\cup_{i\in I}M(i)$.
By an argument similar to \cite[Propsition 4.1 (iii)]{BGT2}, any base change functor $\CAT_{\AAAA_i}^{pc}\to \CAT_{\AAAA_j}^{pc}$
preserves Morita equivalences.
Then we apply \cite[5.3.2]{Maz} to $\mathcal{C}at_{\mathbf{a}_I}^{pc}\to I$
and deduce that the induced map $\mathcal{C}at_{\mathbf{a}_I}^{pc}[M^{-1}]\to I$ is a coCartesian fibration which is classified by the functor
$I\to \wCat$ informally given by $i\mapsto \CAT_{\AAAA_i}^{pc}[M(i)^{-1}]\simeq \ST_{A_i}$ (see \cite{Maz} for the formulation), where $\CAT_{\AAAA_i}[M(i)^{-1}]\simeq \ST_{A_i}$ is the $\infty$-category obtained by inverting Morita equivalences \cite[Section 6]{I}. 
We define $\mathcal{S}t_{a_I}\to I$ to be $\mathcal{C}at_{\mathbf{a}_I}^{pc}[M^{-1}]\to I$. Since a stable equivalence $\RR\to \RR'$
induces a Quillen equivalence $\SPS(\RR)\rightleftarrows \SPS(\RR')$,
it follows from \cite[Theorem 1.4]{Muro}
that it gives rise to the Quillen equivalence 
$\CAT_{\RR}\rightleftarrows \CAT_{\RR'}$. Consequently,
there exists an
equivalence
$\CAT^{pc}_{\RR}[M^{-1}]\simeq \CAT^{pc}_{\RR'}[M^{-1}]$
of $\infty$-categories.
Thus,
it is straightforward to check that $\mathcal{S}t_{a_I}\to I$ does not depend on
a choice of  $\mathbf{a}_I$ up to equivalence.
We define $\mathcal{S}t_{b_I}\to I$ and $\mathcal{S}t_{b_I^\circ}\to I$
in a similar way.

\begin{Lemma}
\label{strictsection}
There exists a  canonical equivalence $\Sect(\mathcal{C}at_{\mathbf{a}_I}^{pc})[M^{-1}]\stackrel{\sim}{\to} \Sect(\mathcal{S}t_{a_I})$.
\end{Lemma}

\Proof
Since cofibrant spectral categories are
pointwise cofibrant spectral categories
it follows from \cite[3.42, 3.27]{Bal}
that the canonical functor $\Sect(\mathcal{C}at_{\mathbf{a}_I}^{pc})[M^{-1}]\to \Sect(\mathcal{C}at_{\mathbf{a}_I}^{pc}[M^{-1}])=\Sect(\mathcal{S}t_{a_I})$ is an equivalence.
\QED

\begin{Construction}
\label{con6}
We consider $\CCC_I$ to be an object of $\Sect(\mathcal{C}at_{\mathbf{a}_I}^{pc}[M^{-1}])\simeq \Sect(\mathcal{S}t_{a_I})$.
Let $\CC_I$ be an object of $\Sect(\mathcal{C}at_{\mathbf{a}_I}^{pc})$,
which represents $\CCC_I$ through the equivalence in Lemma~\ref{strictsection}.
Write $\phi:\HH_\bullet(\CCC_I/k) \to \HH_\bullet(\CCC_I\otimes_{a_I}b_I/b_I^\circ)$ for the morphism in $ \Mod_{\HH_\bullet(a_I/k)}(\Fun(I,\Fun(BS^1,\Mod_k)))$, which is 
obtained by applying $\upsilon$ to $\CC_I$ (see Construction~\ref{con5} (iii)).
Notice that $\HH_\bullet(\CCC_I\otimes_{a_I}b_I/b_I^\circ)$ is the resriction
of a module over $\HH_\bullet(b_I/b_I^\circ)$ along $\HH_\bullet(a_I/k)\to \HH_\bullet(b_I/b_I^\circ)$. Thus, passing to the adjunction we have
\[
\psi_{\CCC_I}:\HH_\bullet(\CCC_I/k)\otimes_{\HH_\bullet(a_I/k)}\HH_\bullet(b_I/b_I^\circ) \to \HH_\bullet(\CCC_I\otimes_{a_I}b_I/b_I^\circ)
\]
that is a morphism in $\Mod_{\HH_\bullet(b_I/b_I^\circ)}(\Fun(I,\Fun(BS^1,\Mod_k)))$.
Here $\otimes_{\HH_\bullet(a_I/k)}\HH_\bullet(b_I/b_I^\circ)$ indicates the base change
$\Mod_{\HH_\bullet(a_I/k)}(\Fun(I,\Fun(BS^1,\Mod_k))) \to \Mod_{\HH_\bullet(b_I/b_I^\circ)}(\Fun(I,\Fun(BS^1,\Mod_k))).$
\end{Construction}

\begin{Remark}
\label{functorialpsi}
Consider two functors
$\gamma\circ \beta$ and $\gamma \circ r \circ \alpha$
(see Construction~\ref{con5} (iii)).
According to \cite[Lemma 6.11]{I} both $\gamma\circ \beta$ and $\gamma \circ r \circ \alpha$
send termwise Morita equivalences to equivalences in $\Mod_{\HH_\bullet(a_I/k)}(\Fun(I,\Fun(BS^1,\Mod_k)))$.
Both $\gamma\circ \beta$ and $\gamma \circ r \circ \alpha$
 factor through the canonical projection $\Sect(\mathcal{C}at_{\mathbf{a}_I}^{pc}) \to \Sect(\mathcal{C}at_{\mathbf{a}_I}^{pc})[M^{-1}]$.
Consequently, these two functors
induce functors $\gamma\circ \beta[M^{-1}],\gamma \circ r \circ \alpha[M^{-1}]:\Sect(\mathcal{C}at_{\mathbf{a}_I}^{pc})[M^{-1}]\to \Mod_{\HH_\bullet(a_I/k)}(\Fun(I,\Fun(BS^1,\Mod_k)))$.
Moreover, we have the natural transfromation
$\upsilon_\infty:\gamma\circ \beta[M^{-1}] \to \gamma \circ r \circ \alpha[M^{-1}]$ induced by $\upsilon$.
Since $\gamma \circ r \circ \alpha[M^{-1}]$ factors as
$\Sect(\mathcal{C}at_{\mathbf{a}_I}^{pc})[M^{-1}]\to \Mod_{\HH_\bullet(b_I/b_I^\circ)}(\Fun(I,\Fun(BS^1,\Mod_k))) \to \Mod_{\HH_\bullet(a_I/k)}(\Fun(I,\Fun(BS^1,\Mod_k)))$, $\psi_{\CCC_I}$ can naturally be 
extended to a natural transformation
\[
\Psi:\Delta^1\times \Sect(\mathcal{C}at_{\mathbf{a}_I}^{pc})[M^{-1}]\simeq \Delta^1\times \Sect(\mathcal{S}t_{a_I})\to \Mod_{\HH_\bullet(b_I/b_I^\circ)}(\Fun(I,\Fun(BS^1,\Mod_k))).
\]
\end{Remark}

\vspace{2mm}

{\it Proof of Theorem~\ref{generalizedeasy}.}
As in Construction~\ref{con6}, we take $\CC_I\in \Sect(\mathcal{C}at_{\mathbf{a}_I}^{pc})$ which represents $\CCC_I\in \Sect(\mathcal{S}t_{a_I})$ (cf. Lemma~\ref{strictsection}). 
Taking into account Construction~\ref{con5} (iii) and Construction~\ref{con6},
$\HH_\bullet(\CCC_I/k)$ is a module over $\HH_\bullet(a_I/k)$,
and $\HH_\bullet(\CCC_I\otimes_{a_I}b_I/b_I^\circ)$ is a module over $\HH_\bullet(b_I/b_I^\circ)$.
Now we prove that $\psi_{\CCC_I}$ is an equivalence.
Take a look at the diagram in (iii) of Construction~\ref{con5}.
The point from Construction~\ref{con5}
is that
the natural transformation $\upsilon$ is induced by the natural transformation
 $\tau':\Delta^1\times \Sect(\mathcal{C}at_{\mathbf{a}_I}^{pc})\to  \Mod_{\HH(a_I/k)_\bullet}(\Fun(I,\Fun(\Lambda^{op},\Mod_k)))$.
Let $C\to D$ denote
the image of $\HH(\mathbf{a}_I/\KK)_\bullet\to \HH(\mathbf{b}_I/\mathbf{b}_I^\circ)_\bullet$
in $\CAlg(\Fun(I,\Fun(\Lambda^{op},\Mod_k)))$.
Write $P\to Q$ for
the image of  $\HH(\pi_{\mathbf{a}_I}(\CC_I)/\KK)_\bullet\to \pi^S_{\mathbf{b}_I^\circ}(\HH(\CC_I\otimes_{\mathbf{a}_I}\mathbf{b}_I/\mathbf{b}_I^\circ)_\bullet)$
in $\Fun(I\times \Lambda^{op},\Mod_k)$
where $\CC_I\otimes_{\mathbf{a}_I}\mathbf{b}_I:=\Sect(\mathcal{C}at_\mathbf{f}^L)(\CC_I)$.
Notice that $\Fun(I,\Fun( \Lambda^{op},\Mod_k))\to \Fun(I,\Fun(BS^1,\Mod_k))$ carries $P\otimes_CD$ to $\HH_\bullet(\CCC_I/k)\otimes_{\HH_\bullet(a_I/k)}\HH_\bullet(b_I/b_I^\circ)$ (because it preserves
colimits),
and $Q$ maps to $\HH_\bullet(\CCC_I\otimes_{a_I}b_I/b_I^\circ)$.
It suffices to show that
$P\otimes_CD\to Q$ is an equivalence in $\Fun(I\times \Lambda^{op},\Mod_k)$.
To this end, it is enough to prove that for any $(i,p)\in I\times \Lambda^{op}$, 
the evaluation $(P\otimes_CD)_{i,p}\to Q_{i,p}$ at
$(i,p)$ is an equivalence in $\Mod_k$.
Taking into account the definitions of the cyclic objects $\HH(\AAAA_i/\KK)_\bullet$
and $\HH(-/\BBBB_i^\circ)_\bullet$,
we are reduced to proving that the canonical map
\begin{eqnarray*}
\bigl(\CC_i(X_{p-1},X_p)\otimes_k\ldots \otimes_k \CC_i(X_0,X_1)\otimes_k \CC_i(X_p,X_0)\bigr)\otimes_{(A_i\otimes_k\cdots \otimes_kA_i)} (B_i\otimes_{B_i^\circ}\cdots \otimes_{B_i^\circ}B)\ \ \ \  \ \ \ \ \ \ \  \\
\ \ \ \ \ \ \ \ \ \  \to (\CC_i(X_{p-1},X_p)\otimes_{A_i}B_i)\otimes_{B_i^\circ}\ldots \otimes_{B^\circ_i} (\CC_i(X_0,X_1)\otimes_{A_i}B_i)\otimes_{B_i^\circ} (\CC_i(X_p,X_0)\otimes_{A_i}B_i)
\end{eqnarray*}
is an equivalence where $\CC_i(X_{l},X_m)$ means the image of $\CC_i(X_{l},X_m)$ in $\Mod_k$, and
$(A_i\otimes_k\cdots \otimes_kA_i)=A_i^{\otimes_k (p+1)}\to B_i^{\otimes_{B_i^\circ}(p+1)}$
is the canonical map. This morphism is obviously an equivalence.
The following Lemma shows the final assertion.
\QED

\begin{Lemma}[Gluing property]
\label{gluing}
Suppose that $\CCC_I$ is coCartesian over $a_I$.
Then $\HH_\bullet(\CCC_I/k)$
is coCartesian over $\HH_\bullet(a_I/k)$.
\end{Lemma}

\Proof
It will suffice to prove that $\HH_\bullet(\CCC_i/k)\otimes_{\HH_{\bullet}(A_i/k)}\HH_\bullet(
A_j/k)\to \HH_\bullet(\CCC_i\otimes_{A_i}A_j/k)$ is an equivalence for any morphism $i\to j$ in $I$.
Taking into account the diagram in (iii) of Construction~\ref{con5},
it is enough to show that
$\HH(\CCC_i/k)_\bullet\otimes_{\HH(A_i/k)_\bullet}\HH(
A_j/k)_\bullet\to \HH(\CCC_i\otimes_{A_i}A_j/k)_\bullet$ is an equivalence
in $\Fun(\Lambda^{op},\Mod_k)$. 
Taking into accout the evaluation at each $[p]\in \Lambda^{op}$,
we are reduced to proving that the canonical map is
\begin{eqnarray*}
\bigl(\CC_i(X_{p-1},X_p)\otimes_k\ldots \otimes_k \CC_i(X_0,X_1)\otimes_k \CC_i(X_p,X_0)\bigr)\otimes_{(A_i\otimes_k\cdots \otimes_kA_i)} (A_j\otimes_{k}\cdots \otimes_{k}A_j)\ \ \ \  \ \ \ \ \ \ \  \\
\ \ \ \ \ \ \ \ \ \  \to (\CC_i(X_{p-1},X_p)\otimes_{A_i}A_j)\otimes_{k}\ldots \otimes_{k} (\CC_i(X_0,X_1)\otimes_{A_i}A_j)\otimes_{k} (\CC_i(X_p,X_0)\otimes_{A_i}A_j)
\end{eqnarray*}
is an equivalence. This is clearly an equivalence.
\QED

\begin{Remark}
The proof of Theorem~\ref{generalizedeasy}
is analogous to that of Theorem~\ref{restrictionfactorization}.
To understand this, note that 
$\Lambda^{op}$ is equivalent to
$(\Disk^\dagger_1)_{/\langle S^1\rangle}= \Disk^\dagger_1\times_{\Mfldf_1}(\Mfldf_1)_{/\langle S^1\rangle}$ where $\langle S^1\rangle\to \Mfldf_1$ is the full faithful embedding
of the full subcategory spanned by $S^1$,
where $\Disk_1^\dagger$ is the full subcategory
of $\Diskf_1$ spanned by nonempty spaces, (see e.g., \cite[Lemma 7.3]{I}).
Based on the local-to-global principle of factorization homology, 
the proof of Theorem~\ref{restrictionfactorization} is reduced to the problem about algebras
on $\Diskf_n$. 
Similarly, in the proof of Theorem~\ref{generalizedeasy},
the problem about Hochschild homology is reduced to the problem about algebras on $(\Disk^\dagger_1)_{/\langle S^1\rangle}\simeq \Lambda^{op}$.
In this sense, the proof of Theorem~\ref{generalizedeasy} is
analogous to that of Theorem~\ref{restrictionfactorization}.
 The gluing property also relies on this.
\end{Remark}

\begin{Example}
\label{geometricex2}
We consider the categorified version of Example~\ref{geometricex1}:
an object of $\Alg_n(\QC(X))$ is replaced with a stable $\infty$-category
over $X$. We use the notation in Example~\ref{geometricex1}.
Let $X$ be a derived scheme over $k$.
Suppose that $I$ and $a_I$ are the same with $I$ and $a_I$ in  Example~\ref{geometricex1}.
Set $b_I=b_I^\circ=a_I$, and $a_I\to b_I$ and $b_I^\circ \to b_I$
are identity morphisms. 
Let
$\CCC_X$ be an object of $\Sect(\mathcal{C}at_{\mathbf{a}_I}^{pc}[M^{-1}])\simeq \Sect(\mathcal{S}t_{a_I})$.
Suppose that $\CCC_X$ is coCartesian.
Define $\ST_X=\lim_{i\in I}\ST_{A_i}$ which we refer to as the $\infty$-category of stable $\infty$-categories over $X$.
Namely, $\CCC_X$ is a stable $\infty$-category
over $X$, that is,
an object of 
$\lim_{i\in I}\ST_{A_i}=\ST_X$.
Define $\HH_\bullet(\CCC_X/X)=\HH_\bullet(\CCC_X/a_I)$ (cf. Construction~\ref{con6}). By Lemma~\ref{gluing} and Theorem~\ref{generalizedeasy}, $\HH_\bullet(\CCC_X/a_I)$
is coCartesian over $a_I$ so that it belongs to $\lim_{i\in I}\Mod_{A_i}\simeq \QC(X)$
(cf. Example~\ref{geometricex1}).
By Lemma~\ref{gluing}, $\HH_\bullet(\CCC_X/k)$
is coCartesian over $\HH_\bullet(a_I/k)$.
As in Example~\ref{geometricex1}, $\Mod_{\HH_\bullet(a_I/k)}(\Fun(I,\Mod_k))^{\textup{coCart}}$ is canonically equivalent to $\QC(LX)$.
Let $\HH_\bullet(\CCC_X/X)_L$ denote
the image of $\HH_\bullet(\CCC_X/k)$ in $\QC(LX)$.
It is naturally promoted to $\QC(LX)^{S^1}$.
We write $\iota:X \to LX$ for the morphism given by constant loops. 
The pullback functor  $\iota^*: \QC(LX)\to \QC(X)$
is induced by the base change
$\Mod_{\HH_\bullet(a_I/k)}(\Fun(I,\Mod_k))^{\textup{coCart}}\to \Mod_{a_I}(\Fun(I,\Mod_k))^{\textup{coCart}}$
along $\HH_\bullet(a_I/k)\to \HH_\bullet(a_I/a_I)\simeq a_I$. 
Consequently, we can regard the equivalence in Theorem~\ref{generalizedeasy}
as $\iota^*(\HH_\bullet(\CCC_X/X)_L)\simeq \HH_\bullet(\CCC_X/X)$
in $\QC(X)^{S^1}$.
\end{Example}

We record the consequence of the observation in Example~\ref{geometricex2}.

\begin{Theorem}
\label{preeasy}
Let $X$ be a derived scheme over a commutative ring spectrum $k$.
Let $\CCC_X$ be a stable $\infty$-category over $X$.
We have a canoncially defined object $\HH_\bullet(\CCC_X/X)_L\in \QC(LX)^{S^1}$
such that the pullback $\iota^* \HH_\bullet(\CCC_X/X)_L$ along $\iota:X\to LX$
is naturally equivalent to $\HH_\bullet(\CCC_X/X)\in \QC(X)^{S^1}$.
\end{Theorem}

\begin{Remark} 
Let $\CCC_{I,1}$, $\CCC_{I,2}$ and $\CCC_{I,3}$
be three objects in $\Sect(\mathcal{S}t_{a_I})$.
For $1\le p \le 3$, we let $\CCC_{i,p}\in \ST_{A_i}$ denote the evaluation of $\CCC_{I,p}$ at $i\in I$.
Let $\CCC_{I,1}\to \CCC_{I,2}\to \CCC_{I,3}$ be an exact sequence of $a_I$-linear
small stable $\infty$-categories: it means that passing to Ind-categories,
for each $i\in I$,
(i) $\Ind(\CCC_{i,1})\to \Ind(\CCC_{i,2})$ is fully faithful, (ii) the composite
$\Ind(\CCC_{i,1})\to \Ind(\CCC_{i,2}) \to \Ind(\CCC_{i,3})$ is trivial, and (iii) 
$\Ind(\CCC_{i,2})/\Ind(\CCC_{i,1})\to \Ind(\CCC_{i,3})$ is an equivalence
where $\Ind(\CCC_{i,2})/\Ind(\CCC_{i,1})$ is a quotient $\Ind(\CCC_{i,2})\sqcup_{\Ind(\CCC_{i,1})}0$ as a presentable stable $\infty$-category
(cf. \cite[Section 5]{BGT1}).
In this situation, the localization theorem says that the induced sequence 
\[
\xymatrix{
\HH_\bullet(\CCC_{I,1}/k)\ar[r] \ar[d] & \HH_\bullet(\CCC_{I,2}/k) \ar[d] \\
0\ar[r] & \HH_\bullet(\CCC_{I,3}/k)
}
\]
is a cofiber/fiber sequence in $\Mod_{\HH_\bullet(a_I/k)}(\Fun(I,\Mod_k^{S^1}))$
(see \cite[Theorem 7.1]{BM} for the localization theorem: in {\it loc cit.} spectral categories are  treated, but the proof is applicable to
$\mathbb{K}$-spectral categories).
\end{Remark}

\section{Kodaira-Spencer morphism for a family of stable $\infty$-categories}
\label{sectionKS}

Let $A$ be a connective commutative dg algebra over $k$ and
let $S=\Spec A$ be a derived affine scheme over $k$.
Suppose that $S$ is locally of finite presentation over $k$.
Assume that $A$ is noetherian.

Let $\CCC=\CCC_A$ be an $A$-linear small stable $\infty$-category (cf. Section~\ref{NCsection}).
Informally, we think of $\CCC_A$ as a family of stable $\infty$-categories over $\Spec A$.
In this section, we will construct a Kodaira-Spencer morphism for $\CCC$ over $A$.
It is defined as a morphism of dg Lie algebras over $A$
\[
\TT_{A/k}[-1]\longrightarrow \HH^\bullet(\CCC/A)[1],
\]
where $\HH^\bullet(\DDD_A/A)[1]$ is the shifted  Hochschild cochain complex with a suitable Lie algebra structure,
and $\TT_{A/k}$ is the tangent complex of $A$ over $k$, that is, the dual object of the cotangent complex $\LL_{A/k}$ in $\Mod_A$.

\subsection{}

Consider 
\[
S\stackrel{\Delta}{\longrightarrow} S\times_k S\stackrel{\textup{pr}_1}{\longrightarrow}S
\]
Let $\widehat{S \times_k S}$ denote the formal completion of $S\times_k S\stackrel{\textup{pr}_1}{\longrightarrow}S$ along $\Delta:S\to S\times_kS$ (see Section~\ref{formalstack}).
The formal stack $\widehat{S \times_k S}:\EXT_A\to \SSS$ is informally given by
\[
\EXT_A\ni [A\to R\to A] \mapsto \Map_{(\Aff_k)_{S//S}}(\Spec R, S\times_k S)\in \SSS.
\]
Recall the equivalence $\FST_A\simeq Lie_A$ (\cite[1.5.6]{H}, Section~\ref{formalstack}).
Then
$\widehat{S \times_k S}$ corresponds to an object of $Lie_A$ whose underlying object
is equivalent to the tangent complex
$\Delta^*(\TT_{S\times_kS/S}[-1])$.
Since $\Delta^*(\LL_{S\times_kS/S})\simeq \Delta^*(\textup{pr}_2^*(\LL_{S/k})))\simeq \LL_{S/k}$
in $\Mod_A$,
it follows that $\Delta^*(\TT_{S\times_kS/S}[-1])\simeq \TT_{S/k}[-1]=:\TT_{A/k}[-1]$.
Here $\LL_{(\ /\ )}$ indicates the cotangent complex.

By abuse of notation we write $\TT_{A/k}[-1]$ for the corresponding object in $Lie_A$.

\subsection{}
\label{deformcat}
We review the functor associated to deformations of $\CCC_A=\CCC\in \ST_A$.
For the detailed construction, we refer to the reader to \cite[Section 6.1]{IMA}.
Let $\Alg_2(\Mod_A)$ be the $\infty$-category of $\etwo$-algebras 
and let $\Alg_2^+(\Mod_A)=\Alg_2(\Mod_A)_{/A}$ be the $\infty$-category of augmented $\etwo$-algebas.
We often omit augmentations from the notation. 
Let $\RMod_{\Perf_B}(\ST_A)$ denote the $\infty$-category
of $\Perf_B$-module objects in $\ST_A$ (note that
$B$ is an $\etwo$-algebra so that $\Perf_B$ belongs to $\Alg_1(\ST_A)$).
For a morphism $B\to B'$ in $\Alg_2(\Mod_A)$,
there is the base change functor $\otimes_{\Perf_B}\Perf_{B'}  :\RMod_{\Perf_B}(\ST_A)\to \RMod_{\Perf_{B'}}(\ST_A)$
which is a left adjoint to the functor
$\RMod_{\Perf_{B'}}(\ST_A)\to \RMod_{\Perf_{B}}(\ST_A)$
induced by the restriction along $\otimes_BB':\Perf_B\to \Perf_{B'}$.
For $B\to A\in \Alg_2^+(\Mod_A)$ we consider the $\infty$-category 
\[
\overline{\Def}^{\etwo}_{\CCC}(B):=\RMod_{\Perf_{B}}(\ST_A)\times_{\ST_A}\{\CCC\}
\]
such that $\RMod_{\Perf_{B}}(\ST_A)\to \ST_A\simeq \RMod_{\Perf_{A}}(\ST_A)$
is given by the base change along $B\to A$.
An object of $\overline{\Def}^{\etwo}_{\CCC}(B)$
can be regarded as a $B$-linear small stable $\infty$-category $\CCC_B$
endowed with an equivalence $\CCC_{B}\otimes_{\Perf_B}\Perf_A\simeq \CCC$
in $\ST_A$.
We think of it
as the $\infty$-category of deformations of $\CCC$ to $B$.
We write $\Def^{\etwo}_{\CCC}(B)$ for the largest $\infty$-groupoid (the largest Kan complex)
contained in $\overline{\Def}^{\etwo}_{\CCC}(B)$.
Let
\[
\Def^{\etwo}_{\CCC}:\Alg_2^+(\Mod_A) \longrightarrow \widehat{\SSS}
\]
be the functor
which is informally given by $[B \to A]\mapsto g(\CCC)(B)$.
We call it the $\etwo$-deformation functor of $\CCC$.

Let $\Def_{\CCC}:\EXT_A \to \widehat{\SSS}$
be the functor obtained from
$\Def^{\etwo}_{\CCC}$ by the composition with the forgetful functor
$\EXT_A\to (\CCAlg_A)_{/A}\simeq  (\CCAlg_k)_{A//A} \rightarrow \Alg_2^+(\Mod_A)$.

\subsection{}

We define a morphism $\widehat{S\times_k S}\to \Def_{\CCC}$ in $\Fun(\EXT_A,\SSS)$.
Let $X:\EXT_A\to \SSS$ be a formal prestack, i.e., a functor.
We define the $\infty$-category of small stable $\infty$-categories over $X$ to be 
$\lim_{\Spec R\to X}\ST_R$ where the limit is taken over $(\TSZ_A)_{/X}$.
Set $\ST_{X}:=\lim_{\Spec R\to X}\ST_R$. There is a canonical ``augmentation'' $\ST_X\to \ST_A$.
The formal prestack $X$ is a colimit of $(\TSZ_A)_{/X}\to \TSZ_A\hookrightarrow \Fun(\EXT_A,\SSS)$.
It follows that there exists a natural equivalence
\[
\Map_{\Fun(\EXT_A,\SSS)}(X,\Def_{\CCC})\simeq \ST^{\simeq}_X\times_{\ST^{\simeq}_A}\{\CCC\}.
\]
Now suppose $X=\widehat{S\times_k S}$.
We define a morphism $\overline{KS}_{\CCC}: \widehat{S\times_k S}\simeq t_{S/k} \to \Def_{\CCC}$
to be an functor correponding to the small stable $\infty$-category $\widehat{\textup{pr}_2^*(\CCC)}$ over $\widehat{S\times_k S}$ (with $\Delta^*\widehat{\textup{pr}_2^*(\CCC)}\simeq \CCC$), that is determined by the base change $\textup{pr}_2^*(\CCC)=\CCC\otimes_{\Perf_A}\Perf_{A\otimes_kA}\in \ST_{A\otimes_kA}$ along the second projection $\textup{pr}_2:S\times_kS=\Spec A\otimes_kA\to \Spec A=S$
in the natural way.

\begin{Remark}
\label{KSuniversal}
Informally, $\overline{KS}_{\CCC}$ on objects can be described as follows:
Let $[A \to R \to A] \in \EXT_A$ and let $\alpha\times \beta:\Spec R\to S\times_k S$ be an object of $(\Aff_k)_{S//S}$ which is regarded
as an object of $\widehat{S\times_k S}(R)$.  Then $\overline{KS}_{\CCC}$ sends $\alpha\times \beta$ to the base change
$\CCC\otimes_{\Perf_{A}}\Perf_R$ along $\beta^*:\Perf_A\to \Perf_R$ with the equivalence $(\CCC\otimes_{\Perf_{A}}\Perf_R)\otimes_{\Perf_R}\Perf_A\simeq \CCC$, which is a deformation of
$\CCC$ to $R$.
\end{Remark}

\subsection{}

\label{appro}

Let $\DD_n:\Alg_{\eenu}^+(\Mod_A)^{op} \to \Alg_{\eenu}^+(\Mod_A)$ be the $\eenu$-Koszul duality
functor which sends an augmented $\eenu$-algebra $p:B\to A$ in $\Mod_A$
to the $\eenu$-Koszul dual $\DD_{n}(B)\to A$.
The $\eenu$-Koszul dual $\hat{q}:\DD_{n}(B)\to A$ has the following characterization:
it represents the functor $\Alg_n(\Mod_A)^{op}\to \SSS$ informally given by
\[
[q:C\to A] \mapsto \Map_{\Alg_{n}(\Mod_A)}(B\otimes_AC,A)\times_{\Map_{\Alg_{n}(\Mod_A)}(B,A)\times \Map_{\Alg_{n}(\Mod_A)}(C,A)}\{(p,q)\}.
\]
In particular, there is a tautological morphism $B\otimes_A\DD_n(B)\to A$ in 
$\Alg_{n}(\Mod_A)$ which extends $p$ and $\hat{p}$ up to homotopy.
For Koszul duality of $\eenu$-algebras we refer the reader to \cite{HA} or \cite[X 4.4]{DAG}
(see also \cite[Section 3]{IMA} for the brief review).

Let $\FF_{A\oplus \HH^\bullet(\CCC/A)}^{\etwo}:\Alg_{\etwo}^+(\Mod_A)\to \SSS$ 
be a functor informally given by
\[
[B\to A] \mapsto \Map_{\Alg_{\etwo}(\Mod_A)}(\DD_{2}(B),\HH^\bullet(\CCC/A)).
\]
In other words,  $\FF_{A\oplus \HH^\bullet(\CCC/A)}^{\etwo}$ is the composite of
$\DD_2$ and the functor $\Alg_{\etwo}^+(\Mod_A)^{op}\to \SSS$ represented by
$\textup{pr}_2:A\oplus \HH^\bullet(\CCC/A)\to A$.
Let $\FF_{A\oplus \HH^\bullet(\CCC/A)}:\EXT_A \to \SSS$ 
be the functor obtained from $\FF_{A\oplus \HH^\bullet(\CCC/A)}^{\etwo}$
by the composition with $\EXT_A\hookrightarrow (\CCAlg_{k})_{A//A} \to \Alg_{\etwo}^+(\Mod_A)$.
See \cite[Construction 8.6, Lemma 8.11]{IMA} for $\FF_{A\oplus \HH^\bullet(\CCC/A)}$.

There is a (canonical) morphism
\[
J_{\CCC}^{\etwo}:\Def^{\etwo}_{\CCC}\longrightarrow \FF_{A\oplus \HH^\bullet(\CCC/A)}^{\etwo}
\]
in $\Fun(\Alg^+_{\etwo}(\Mod_A),\SSS)$
(see \cite[Construction 8.6]{IMA} for the precise construction).
The composition with $\EXT_A \to \Alg_{\etwo}^+(\Mod_A)$
induces a morphism
\[
J_{\CCC}:\Def_{\CCC}\longrightarrow \FF_{A\oplus \HH^\bullet(\CCC/A)}
\]
in $\Fun(\EXT_A,\SSS)$.

\begin{Remark}
Roughly speaking, $J_{\CCC}^{\etwo}$ is universal among morphisms
into functors arising from augmented $\etwo$-algebras
(see \cite[X, 5.3.16]{DAG} for a statement: we will not use the universal property, but the argument in {\it loc. cit.}
can be applied to our situation).
Informally, $J_{\CCC}^{\etwo}$ can be described as follows.
Let $B\to A$ be in $\Alg^+_{\etwo}(\Mod_A)$ and let $\overline{\CCC}_B=(\CCC_B,\CCC_B\otimes_{\Perf_B}\Perf_A\simeq \CCC)$ be a deformation of $\CCC$ to $B$. 
The tautological morphism $B\otimes_A \DD_2(B)\to A$
induces a morphism
$\Perf_B\otimes \Perf_{\DD_2(B)}\to \Perf_A$ in $\Alg_1(\ST_A)$,
so that $\Perf_A$ is a left $\Perf_B\otimes \Perf_{\DD_2(B)}$-module.
Thus, we may regard $\Perf_A$ as a $\Perf_B$-$(\Perf_{\DD_2(B)})^{op}$-bimodule
where $(\Perf_{\DD_2(B)})^{op}$ means the ``opposite algebra'' of $\Perf_{\DD_2(B)}\in \Alg_1(\ST_A)$.
The relative tensor product $\CCC_B\otimes_{\Perf_B}\Perf_A$
is a right $(\Perf_{\DD_2(B)})^{op}$-module,
that is, a left $\Perf_{\DD_2(B)}$-module. 
In particular, $\CCC$ is promoted to a left $\Perf_{\DD_{2}(B)}$-module.
By the $\Ind$-construction 
it gives rise to a left $\LMod_{\DD_{2}(B)}$-module structure on
$\DDD:=\Ind(\CCC)\in \PR_A$.
A left $\LMod_{\DD_{2}(B)}$-module structure on $\DDD$ amounts to a morphism $\alpha_{\overline{\CCC}_B}:\LMod^\otimes_{\DD_{2}(B)}\to \mathcal{E}nd_{A}(\DDD)$ in $\Alg_1(\PR_A)$
(see Section~\ref{hochschildhomologysection} for $ \mathcal{E}nd_{A}(\DDD)$).
By the definition of Hochschild cohomology $\HH^\bullet(\CCC/A)$,
there is a canonical equivalence $\Map_{\Alg_2(\Mod_A)}(\DD_2(B),\HH^\bullet(\CCC/A))\simeq \Map_{\Alg_1(\PR_A)}(\LMod_{\DD_{2}(B)},\mathcal{E}nd_{A}(\DDD))$.
Let $\beta_{\overline{\CCC}_B}:\DD_2(B)\to \HH^\bullet(\CCC/A)$ be a morphism corresponding to 
$\alpha_{\overline{\CCC}_B}$
through the equivalence. Then $J_{\CCC}^{\etwo}$ carries $\overline{\CCC}_B$
to $\beta_{\overline{\CCC}_B}$.
\end{Remark}

\subsection{}
Let 
\[
\xymatrix{
U_2:Lie_A  \ar@<0.5ex>[r] &    \Alg_{\etwo}^+(\Mod_A): res_{\etwo/Lie}.   \ar@<0.5ex>[l]  
}
\]
be the adjoint pair where the left adjoint $U_2$ is 
the universal enveloping $\etwo$-algebra functor, that is an ``$\etwo$-generalization''
of th universal enveloping functor (cf. \cite[Section 3.5]{IMA}).
The composite $\Alg_{\etwo}^+(\Mod_A)\stackrel{res_{\etwo/Lie}}{\longrightarrow} Lie_A\stackrel{\textup{forget}}{\longrightarrow} \Mod_A$
is equivalent to the shift functor given by $B\mapsto B[1]$.

Set $\GG_{\CCC}=res_{\etwo/Lie}(A\oplus \HH^\bullet(\CCC/A)\to A)$.
The underlying object of $\GG_{\CCC}$ is equivalent to $\HH^\bullet(\CCC/A)[1]$.
Let $\mathcal{F}_{\GG_{\CCC}}$ be the formal stack associated to $\GG_{\CCC}$.
As proved in \cite[Lemma 8.11]{IMA}, there exists a canonical equivalence
\[
\mathcal{F}_{A\oplus \HH^\bullet(\CCC/A)}\simeq \mathcal{F}_{\GG_{\CCC}}.
\]

\begin{Construction}
We consider the sequence of functors
\[
\mathcal{F}_{\TT_{A/k}[-1]} \simeq  \widehat{S\times_kS} \stackrel{\overline{KS}_{\CCC}}{\longrightarrow}  \Def_{\CCC} \stackrel{J_{\CCC}}{\longrightarrow} \mathcal{F}_{A\oplus \HH^\bullet(\CCC/A)}\simeq \mathcal{F}_{\GG_{\CCC}}
\]
Since $A$ is noetherian, there is the equivalence $\FST_A\simeq Lie_A$.
The composite functor $\mathcal{F}_{\TT_{A/k}[-1]} \to \mathcal{F}_{\GG_{\CCC}}$
is a morphism in $\FST_A$.
Passing to $Lie_A$ 
 we have a morphism in $Lie_A$
\[
KS_{\CCC}:\mathbb{T}_{A/k}[-1]\longrightarrow \GG_{\CCC}.
\]
We refer to $KS_{\CCC}$ as the Kodaira-Spencer morphism for $\CCC$.
\end{Construction}

\begin{Remark}
To define $KS_{\CCC}$, the noetherian condition is not necessary
because we need only the functor $\mathcal{L}:\FST_A\to Lie_A$.
To understand $KS_{\CCC}$, let us consider the simple situation.
Suppose that $R=A\oplus A[n]\in \EXT_A$ is the trivial square extension of $A$ by $A[n]$ ($n\ge0$).
For $L\in Lie_A$, there are equivalences
$\mathcal{F}_{L}(R)\simeq \Map_{Lie_A}(\textup{Free}_{Lie}(A[-n-1]),L)\simeq \Omega^{\infty}(L[n+1])$
in $\SSS$.
We will consider $\mathcal{F}_{\TT_{A/k}[-1]}(R)\simeq \Omega^\infty(\TT_{A/k}[n])\to \Omega^\infty(\HH^\bullet(\CCC/A)[n+2])\simeq \mathcal{F}_{\GG_{\CCC}}(R)$.
The space $\mathcal{F}_{\TT_{A/k}[-1]}(R)$ can be identified with $\Map_{(\CCAlg_k)_{/A}}(id:A\to A, R\to A)$.
Let $p\in \pi_0(\Omega^{\infty}(\TT_{A/k}[n]))$
and suppose that $f_p:A\to R$ in $(\CCAlg_k)_{/A}$
corresponds to $p$. Then we have a deformation $f_p^*\CCC:=\CCC\otimes_{\Perf_A}\Perf_R$
of $\CCC$ to $R$.
The image of $f_p^*\CCC$ under $J_{\CCC}$ determines a class/element $cl(f_p^*\CCC)$ of $\pi_0(\Omega^\infty(\HH^\bullet(\CCC/A)[n+2]))\simeq H^{n+2}(\HH^\bullet(\CCC/A))=:\HH^{n+2}(\CCC/A)$.
By construction, $\pi_0(\Omega^\infty(\TT_{A/k}[n]))= H^n(\TT_{A/k})\to 
\pi_0(\Omega^\infty(\HH^\bullet(\CCC/A)[n+2]))\simeq \HH^{n+2}(\CCC/A)$
sends $p$ to $cl(f^*_p\CCC)$.
\end{Remark}

\section{Lie algebra representations from the Hochschild pair}

\subsection{The Lie algebra action and  the extended Lie algebra action on Hochschild Homology.}

\label{pairreview}

We assume that $A$ is a commutative noetherian dg algebra over $k$.
Suppose that $\CCC$ is an $A$-linear small stable $\infty$-category. 
We first recall an algebraic structure on
the pair $(\HH^\bullet(\CCC/A),\HH_\bullet(\CCC/A))$.
The pair $(\HH^\bullet(\CCC/A),\HH_\bullet(\CCC/A))$
is an algebra over the so-called Kontsevich-Soibelman operad $\KS$.
We use the construction in our previous work \cite{I}.
We will not recall the (colored) operad $\KS$ in this paper,
but we will explain equivalent algebraic data.
According to \cite[Theorem 1.2]{I},
an algebra in $\Mod_A$ over $\KS$ amounts to a triple of the following data
\begin{enumerate}

\item an $S^1$-action on $\HH_\bullet(\CCC/A)$ in $\Mod_A$. That is, $\HH_\bullet(\CCC/A)$ is an object of $\HH_\bullet(\CCC/A)\in \Mod_A^{S^1}=\Fun(BS^1,\Mod_A)$.

\item an $\etwo$-algebra structure on $\HH^\bullet(\CCC/A)$,

\item an $S^1$-equivariant left module action of $\HH_\bullet(\HH^\bullet(\CCC/A)/A))$ on $\HH_\bullet(\CCC/A)$. \\
 (Note that $\HH_\bullet(\HH^\bullet(\CCC/A)/A))\in \Alg_1(\Mod_A^{S^1})$
and $\HH_\bullet(\CCC/A)\in \Mod_A^{S^1}$.) 

\end{enumerate}
More precisely, 
an algebra in $\Mod_A$ over $\KS$ is equivalent to giving an object of
\[
\Alg_{2}(\Mod_A)\times_{\Alg_1(\Mod_A^{S^1})}\LMod(\Mod_A^{S^1})
\]
where $\Alg_{2}(\Mod_A) \to \Alg_1(\Mod_A^{S^1})$ is induced by $\HH_\bullet(-/A)$.

Next, we recall Lie algebra actions on $\HH_\bullet(\CCC/A)$ arising from the pair $(\HH^\bullet(\CCC/A),\HH_\bullet(\CCC/A))$.
For details we refer to \cite[Section 8.1]{IMA}.
Let $\End(\HH_\bullet(\CCC/A))$ denote the endomorphism associative algebra object in $\Mod_A^{S^1}$ (cf. \cite[4.7.1]{HA}).
The above datum (3) is equivalent to giving a morphism in $\Alg_{\eone}(\Mod_A^{S^1})$
\[
\HH_\bullet(\HH^\bullet(\CCC/A)/A) \longrightarrow \End(\HH_\bullet(\CCC/A)).
\]
For $L\in Lie_A$ we let $L^{S1}$ denote the cotensor of $L$ by $S^1$.
More explicitly, $L^{S1}$ is given by the fiber product $L\times_{L\times L} L$ in $Lie_A$.
(However, there is one exception to this: $\End^L(\HH_\bullet(\CCC/A))^{S^1}$ appeared below indicates
the homotopy fixed points of an $S^1$-action.)
Let $U_1:Lie_A\to \Alg_1(\Mod_A)$ denote the universal enveloping algebra functor.
Let $U_2(\GG_{\CCC})\to \HH^\bullet(\CCC/A)$
be the counit map arising from the adjoint pair $U_2:Lie_A\rightleftarrows \Alg_{\etwo}(\Mod_A)$
(we here adopt the unaugmented version).
It gives rise to a sequence of morphisms 
\[
U_1(\GG_{\CCC}^{S^1})\simeq \HH_\bullet(U_2(\GG_{\CCC})/A)\to \HH_\bullet(\HH^\bullet(\CCC)/A)\to \End(\HH_\bullet(\CCC/A))
\]
in $\Alg_1(\Mod_A^{S^1})$,
where the first equivalence is proved in \cite[Propostion 7.21]{IMA}, and 
the second morphism is induced by $U_2(\GG_{\CCC})\to \HH^\bullet(\CCC/A)$.
Let $i:\GG_{\CCC}  \to \GG_{\CCC} ^{S^1}$ be the morphism
in $Lie_A$ induced by $S^1 \to \ast$.
Since $S\to \ast$ is $S^1$-equivariant,
$\GG_{\CCC}  \to \GG_{\CCC} ^{S^1}$ is promoted to
a morphism in $Lie_A^{S^1}=\Fun(BS^1,Lie_A)$ where the $S^1$-action on $\GG_{\CCC}$
is trivial.
Using the composite of the above sequence
we have a sequence  of morphisms
\[
U_1(\GG_{\CCC})\to U_1(\GG_{\CCC}^{S^1})\to \End(\HH_\bullet(\CCC/A)).
\]
in $\Alg_1(\Mod_A^{S^1})$.
Let $\End^L(\HH_\bullet(\CCC/A))\in Lie_A^{S^1}$ denote the dg Lie algebra with $S^1$-action associated to $ \End(\HH_\bullet(\CCC/A))$.
Then these morphisms give rise to morphisms in $Lie_A^{S^1}$:
\[
\GG_{\CCC}\stackrel{i}{\to} \GG_{\CCC}^{S^1} \stackrel{\widehat{A}^L_{\CCC}}{\longrightarrow} \End^L(\HH_\bullet(\CCC/A)).
\]
We write $A^L_{\CCC}$
for the composite $\widehat{A}^L_{\CCC}\circ i:\GG_{\CCC}\to  \End^L(\HH_\bullet(\CCC/A))$.
(This symbol $A^L_{\CCC}$ is different from $A_{\CCC}^{L}$ in \cite{IMA} which means
the induced morphism $\GG_{\CCC}\to  \End^L(\HH_\bullet(\CCC/A))^{S^1}$ in $Lie_A$.)

\begin{Definition}
We call $A^L_{\CCC}$ the canonical Lie algebra action of $\GG_{\CCC}$ on $\HH_\bullet(\CCC/A)$.
We call $\widehat{A}^L_{\CCC}$ the canonical extended Lie algebra action of $\GG_{\CCC}^{S^1}$ on $\HH_\bullet(\CCC/A)$.
\end{Definition}

\subsection{Lie algebra actions of $\mathbb{T}_{A/k}[-1]$ and $\mathbb{T}_{A/k}[-1]^{S^1}$.}

\label{startlie}

\begin{Construction}
\label{const2}
Recall that $\mathbb{T}_{A/k}[-1]$ is the dg Lie algebra that corresponds to $\widehat{S\times_kS}$.
Let $KS_{\CCC}:\mathbb{T}_{A/k}[-1]\rightarrow \GG_{\CCC}$
be the Kodaira-Spencer morphism (see Section~\ref{sectionKS}).
Recall the sequence $\GG_{\CCC}\stackrel{i}{\to} \GG_{\CCC}^{S^1} \stackrel{\widehat{A}_{\CCC}}{\to} \End^L(\HH_\bullet(\CCC/A))$
(see Section~\ref{pairreview}).
We obtain the commutative digaram in $Lie(\Mod_A^{S^1})$:
\[
\xymatrix{
\mathbb{T}_{A/k}[-1] \ar[r] \ar[d]_{KS_{\CCC}} & \mathbb{T}_{A/k}[-1]^{S^1} \ar[d]^{KS^{S^1}_{\CCC}}  &  \\
\GG_{\CCC}  \ar[r]^i & \GG_{\CCC} ^{S^1}\ar[r]^(0.3){\widehat{A}^L_{\CCC}}  & \End^L(\HH_\bullet(\CCC/A))
}
\]
where $\mathbb{T}_{A/k}[-1]^{S^1}$ is the cotensor of $\mathbb{T}_{A/k}[-1]$ by $S^1$,
and $\mathbb{T}_{A/k}[-1] \to \mathbb{T}_{A/k}[-1]^{S^1}$ is induced by $S^1\to \ast$ ($\ast$
is the contractible space). 
The vertical morphism $KS_{\CCC}^{S^1}:\mathbb{T}_{A/k}[-1]^{S^1}\rightarrow \GG^{S^1}_{\CCC}$
is induced by $KS_{\CCC}$.
The dg Lie algebras $\mathbb{T}_{A/k}[-1]$ 
is endowed with the trivial $S^1$-action.
It gives rise to an $S^1$-equivariant Lie algebra action of $\mathbb{T}_{A/k}[-1]$ on $\HH_{\bullet}(\CCC/A)$, that  is, $\mathbb{T}_{A/k}[-1]\to \End^L(\HH_\bullet(\CCC/A))$ in $Lie_A^{S^1}$. Moreover,
this action factors through $\mathbb{T}_{A/k}[-1] \to \mathbb{T}_{A/k}[-1]^{S^1}$.
The $S^1$-action on $\mathbb{T}_{A/k}[-1]^{S^1}$
comes from the action of $S^1$ on $S^1$ defined by the multiplication $S^1\times S^1\to S^1$.
\end{Construction}

We give immediate consequences of Construction~\ref{const2}.
Let $L$ be a dg Lie algebra over $A$.
Let $\Rep(L)(\Mod_A)$ be the stable $\infty$-category
of representations of $L$.
We define $\Rep(L)(\Mod_A)$ to be
$\LMod_{U_1(L)}(\Mod_A)$
where $U_1$ is the universal enveloping algebra functor $Lie_A\to \Alg_1(\Mod_A)$
(see e.g. \cite[Section 3.5]{IMA}).
By the diagram in Construction~\ref{const2}, $\HH_\bullet(\CCC/A)$ 
is promoted to an object of
$\Rep(\mathbb{T}_{A/k}[-1])(\Mod_A)$.
Note that the action of $\mathbb{T}_{A/k}[-1]$ on $\HH_\bullet(\CCC/A)$
is extended to an action of $\mathbb{T}_{A/k}[-1]^{S^1}$.
We note that the diagram in Construction~\ref{const2}
is a diagram in $\Fun(BS^1,\Mod_A)$.
For $L\in \Fun(BS^1,Lie_A)$, we set $\Rep(L)(\Mod^{S^1}_A)=\LMod_{U_1(L)}(\Mod_A^{S^1})$.
Then $\HH_\bullet(\CCC/A)$ is promoted to an object of $\Rep(\mathbb{T}_{A/k}[-1]^{S^1})(\Mod^{S^1}_A)$.


\begin{Definition}
\label{canonicaltangentaction}
We shall  refer to an object of $\Rep(\mathbb{T}_{A/k}[-1])(\Mod_A^{S^1})$ determined by
 $A^L_{\CCC}\circ KS_{\CCC}:\mathbb{T}_{A/k}[-1]\to \End^L(\HH_\bullet(\CCC/A))$
as the canonical $\mathbb{T}_{A/k}[-1]$-module $\HH_\bullet(\CCC/A)$.
We shall  refer to an object of $\Rep(\mathbb{T}_{A/k}[-1]^{S^1})(\Mod_A^{S^1})$ determined by
 $\widehat{A}^L_{\CCC}\circ KS_{\CCC}^{S^1}:\mathbb{T}_{A/k}[-1]^{S^1}\to \End^L(\HH_\bullet(\CCC/A))$ 
as the canonical $\mathbb{T}_{A/k}[-1]^{S^1}$-module $\HH_\bullet(\CCC/A)$.
\end{Definition}

\vspace{2mm}

Let us describe $\mathbb{T}_{A/k}[-1] \to \mathbb{T}_{A/k}[-1]^{S^1}$
in terms of formal stacks.
For this purpose, we first recall the (free) loop space of a derived scheme
from the viewpoint of functors (cf. Section~\ref{derivedscheme}).
Let $X$ be a derived scheme over $k$.
Let $LX$ be the free loop space of $X$ which is defined as $LX\simeq X\times_{X\times_kX}X$.
If we regard it
as a functor,
$LX:\CCAlg_k\to \SSS$ is given by $R\to \Map_{\SSS}(S^1,X(R))\simeq X(R)\times_{X(R)\times X(R)}X(R)$
 where $X(R)$ is the space $\Map_{\Fun(\CCAlg_k, \SSS)}(\Spec R,X)$ of $R$-valued
points.
The derived scheme $LX$ has the $S^1$-action induced by the canonical action on the domain in $\Map_{\SSS}(S^1,X(R))$. The evident map $S^1\to \ast$ to the contractible space $\ast$
induces $X(R)=\Map_{\SSS}(\ast,X(R))\to \Map_{\SSS}(S^1,X(R))$, so that
there is a $S^1$-equivariant canonical morphism $\iota:X\to LX$ where the $S^1$-action on $X$
is trivial. This may be regared as the morphism induced by constant loops.

Let $\widehat{S\times_kLS}$ be the formal stack defined 
as the formal completion of $S\times_kLS$ along
the graph map $S\to S\times_kLS$ induced by $\iota:S\to LS$.
Since $\widehat{S\times_kLS}$ can be obtained from $\widehat{S\times_kS}$
by the cotensor by $S^1$ $\in \FST_A$,
$\widehat{S\times_kLS}$ corresponds to
the dg Lie algebra $\mathbb{T}_{A/k}[-1]^{S^1}$ (obtained from $\mathbb{T}_{A/k}[-1]$ by
the cotensor by $S^1$).
The morphism $\mathbb{T}_{A/k}[-1]\to \mathbb{T}_{A/k}[-1]^{S^1}$ induced by $S^1\to \ast$
corresponds to
\[
\widehat{S\times_k S}\to \widehat{S\times_kLS}
\]
induced by $S\to LS$.

\subsection{Complexes/sheaves on formal stacks.}

\label{sheaves}

Let $X:\EXT_A\to \SSS$ be a functor, that is, a formal prestack.
We define the stable $\infty$-category  $\QC_H(X)$
of quasi-coherent complexes on $X$ (they are also called quasi-coherent sheaves in the literature,
but we prefer to call them complexes).
Let $\EXT_A\to \wCat$ be the functor which carries $R$ to $\Mod_R:=\Mod_R(\Mod_A)$.
A morphism $R\to R'$ maps to the base change functor
$(-)\otimes_RR':\Mod_R\to \Mod_{R'}$.
The functor corresponds to the base change of the coCartesian 
fibration $\Mod(\Mod_A)\to \CAlg(\Mod_A)$ along $\EXT_A\to \CAlg(\Mod_A)$.
Let $\EXT_A\to \Fun(\EXT_A,\SSS)^{op}$ be the fully faithful functor induced by the Yoneda embedding.
Since $\wCat$ has small limits, it follows that
there exists a right Kan extension
\[
\QC_H:\Fun(\EXT_A,\SSS)^{op}\to \wCat
\]
where 
$\QC_H(X)$ is equivalent to $\lim_{\Spec R\to X}\Mod_{R}$ where the limit is taken over
$(\TSZ_A)_{/X}$ ($\TSZ_A$ denotes the opposite category of $\EXT_A$).
If $f:X\to Y$ is a morphism of formal prestacks, we denote
by $f^*:\QC_H(Y)\to \QC_H(X)$ the functor induced by $\QC_H$.
We refer to $f^*$ as the pullback functor.
Since $\EXT_A$ is a full subcategory of $(\CAlg_k)_{A//A}$
and $\QC_H(A)=\Mod_A$, the functor $\QC_H$
can be extended to $\Fun(\EXT_A,\SSS)^{op}\to (\wCat)_{\Mod_A//\Mod_A}$.

Next,
we define the stable $\infty$-category $\Rep_H(X)$.
Let $\Alg_1(\Mod_A)^{op}\to \wCat$ be the functor informally given by $B\mapsto \LMod_B:=\LMod_B(\Mod_A)$, which correponds to the Cartesian fibration $\LMod(\Mod_A)\to \Alg_1(\Mod_A)$.
Consider the composite
\[
u:\EXT_A \stackrel{\textup{forget}}{\longrightarrow} \Alg_1^+(\Mod_A)\stackrel{\DD_1}{\longrightarrow}  \Alg^+_1(\Mod_A)^{op}\to \wCat.
\]
where $\DD_1$ is $\eone$-Koszul duality functor (cf. e.g. \cite[X, 4.4]{DAG}, 
Section~\ref{appro}, \cite{IMA}).
We define
\[
\Rep_{H}:\Fun(\EXT_A,\SSS)^{op}\to \wCat
\]
to be a right Kan extension of $u$ along $\EXT_A\to \Fun(\EXT_A,\SSS)^{op}$.
By definition, $\Rep_H(X)$
is the limit $\lim_{\Spec R\to X}\LMod_{\DD_1(R)}(\Mod_A)$ where 
the limit is taken over $(\TSZ_A)_{/X}$.
The functor
$\Rep_{H}$ has an extension $\Fun(\EXT_A,\SSS)^{op}\to (\wCat)_{\Mod_A//\Mod_A}$
for the same reason as above.
Since $\DD_1(R)\simeq U_1(\DD_{\infty}(R))$ for $R\in \EXT_A$ (see e.g. \cite[Proposition 3.3]{IMA}),
there exist categorical equivalences
$u(R)\simeq \LMod_{\DD_1(R)}(\Mod_A)\simeq \LMod_{U_1(\DD_{\infty}(R))}(\Mod_A)=\Rep(\DD_{\infty}(R))(\Mod_A)$.

\begin{Lemma}
\label{replimit}
Let $L\in Lie_A$.
Let $\FF_L=X$ be the formal stack associated to $L$. 
Then there exists a canonical equivalence
$\Rep_H(X)\simeq \Rep(L)(\Mod_A)$.
\end{Lemma}

\Proof
Let $Lie_A^{f}$ be the full subcategory of $Lie_A$ spanned by
free dg Lie algebras of the form $\textup{Free}_{Lie}(\oplus_{1\le i\le n}A^{\oplus r_i}[p_i])$ ($p_i\le -1$)
(see Section~\ref{formalstack}).
The equivalence $\DD_{\infty}:\TSZ_A\simeq  Lie_A^{f}$
induces $\Fun(\EXT_A,\SSS)\simeq \Fun((Lie_A^f)^{op},\SSS)$.
The fully faithful right adjoint of the composite of localization functors
$\Fun((Lie_A^f)^{op},\SSS)\to \mathcal{P}_{\Sigma}(Lie_A^f)\to \mathcal{P}_{\Sigma}^{st}(Lie_A^f)\simeq Lie_A$
is given by the Yoneda embedding followed by the restriction $Lie_A\to \Fun(Lie_A,\SSS)\to \Fun((Lie_A^f)^{op},\SSS)$ (see \cite[1.2.2]{H}, Section~\ref{stablederived} for $\mathcal{P}_{\Sigma}^{st}(Lie_A^f)\simeq Lie_A$).
It follows that
$\Rep_H(\mathcal{F}_L)$ is a colimit of $(Lie_A^f)_{/L}\to Lie^f_A\stackrel{u\circ Ch^\bullet}{\longrightarrow} \wCat^{op}$.
The functor $u\circ Ch^\bullet$ carries $P\in Lie_A^f$ to $\Rep(P)(\Mod_A)=\LMod_{U_1(P)}$.
By Claim~\ref{sifteddiagram} below, the $\infty$-category $(\TSZ_A)_{/X}\simeq (Lie_A^{f})_{/L}$
is sifted.
By \cite[X, 2.4.32, 2.4.33]{DAG},  a sifted colimit $\colim_{i\in I}L_i=L$ of dg Lie algebras gives rise to a canonical equivalence
$\Rep(L)(\Mod_A)\simeq \lim_{i\in I}\Rep(L_i)(\Mod_A)$.
Thus, we obtain $\Rep(L)(\Mod_A) \stackrel{\sim}{\to}  \lim_{[P\to L]\in (Lie_A^f)_{/L}}\Rep(P)(\Mod_A)=\Rep_H(\mathcal{F}_L)$.
\QED

\begin{Claim}
\label{sifteddiagram}
Let $X$ be a formal stack over $A$.
Then $(\TSZ_A)_{/X}$ is sifted.
\end{Claim}

\Proof
Put $D=(\TSZ_A)_{/X}$.
It is enough to prove that
the diagonal functor $D\to D\times D$ is a right adjoint (in particular, it is cofinal).
We define $D\times D\to D$ by 
$(\Spec A\oplus M\to X, \Spec A\oplus M'\to X)\mapsto \Spec A\oplus (M\oplus M')\to X$
(keep in mind that there is the canonical equivalence
 $X(A\oplus (M\times M'))\simeq X(A\oplus M)\times X(A\oplus M')$
 since $X$ is a formal stack).
This functor is a left adjoint to the diagonal functor.
\QED

\begin{Remark}
\label{liepull}
If $\phi:L\to L'$ is a morphism of dg Lie algebras, it induces a morphism $\mathcal{F}_{\phi}:\mathcal{F}_L\to \mathcal{F}_{L'}$
of formal stacks. 
The contravariant functor $\Rep_H$ gives rise to the pullback functor $\mathcal{F}_{\phi}^!:\Rep_H(\mathcal{F}_{L'})\to \Rep_H(\mathcal{F}_{L})$.
The pullback $\mathcal{F}_{\phi}^!$ can naturally be identified with the restriction functor
$\Rep(L')(\Mod_A)\to \Rep(L)(\Mod_A)$ determined by $\phi$.
Indeed, there is the following diagram of restriction functors
\[
\xymatrix{
 \lim_{[P\to L']\in (Lie_A^f)_{/L'}}\Rep(P)(\Mod_A) \ar[r] &  \lim_{[P\to L]\in (Lie_A^f)_{/L}}\Rep(P)(\Mod_A) \\
 \Rep(L')(\Mod_A) \ar[u]^{\simeq} \ar[r]& \Rep(L)(\Mod_A) \ar[u]^{\simeq}
}
\]
which commutes up to canonical homotopy.
\end{Remark}

Let $X$ be a functor $\EXT_A\to \SSS$.
We will construct a fully faithful functor $\QC_H(X)\to \Rep_H(X)$.
For $R \in \EXT_A$, the functor
$\Mod_{R}(\Mod_A) \to \LMod_{\DD_1(R)}(\Mod_A)$ which sends $P$ to $P\otimes_{R}A$.
More precisely,
$\DD_1(R)\otimes_A R \to A$
exhibits $A$ as a $\DD_1(R)$-$R$-bimodule.
This integral kernel $A$ determines a functor
$I_R:\Mod_{R}(\Mod_A) \to \LMod_{\DD_1(R)}(\Mod_A) \simeq \Rep(\DD_\infty(R))(\Mod_A)$
given the formula $P\mapsto P\otimes_{R}A$.
By \cite[2.3.6]{H},
this functor is fully faithful.
By the construction in \cite[Remark 5.9]{IMA},
$I_R$ is functorial in $R\in \EXT_A$. That is,
there is a natural transformation
$\Mod_{(-)}(\Mod_A) \to \LMod_{\DD_1(-)}(\Mod_A)$ between functors $\EXT_A \to \wCat$, such that the evaluation at each $R\in \EXT_A$ is equivalent to $I_R$.
By the definition of $\QC_H$ and $\Rep_H$ (as right Kan extensions),
it gives rise to a morphism $\gamma:\QC_H\to \Rep_{H}$ in $\Fun(\Fun(\EXT_A,\SSS)^{op}, \wCat)$.
For a functor $X:\EXT_A\to \SSS$
we have a canonical fully faithful functor
\[
\gamma_X:\QC_H(X)=\lim_{\Spec R\to X}\Mod_{R}(\Mod_A)\to \lim_{\Spec R\to X} \LMod_{\DD_1(R)}(\Mod_A)=\Rep_H(X).
\]
Taking into account $\gamma_X$, we often regard
$\QC_H(X)$ as a stable subcategory of $\Rep_H(X)$.

\subsection{}
\label{somediagrams}

Let $N$ be an $A$-module spectrum, that is, an object of $\Mod_A$.
Let $\Def(N):\EXT_A\to \SSS$ be the deformation functor
which assigns $[R\to A]$ to the space of deformations of $N$ to $\Mod_R$.
Namely, $\Def(N)$ is given by $[R\to A]\mapsto \Mod_R^{\simeq}\times_{\Mod_A^\simeq}\{N\}$ 
(see \cite[Section 6.2]{IMA}
or \cite[X, 5.2]{DAG}
for a precise formulation).
As with $\Def_{\CCC}$,
for any formal prestack $X:\EXT_A\to \SSS$, there is a canonical equivalence
\[
\Map_{\Fun(\EXT_A,\SSS)}(X,\Def(N))\simeq \QC_H(X)^{\simeq}\times_{\Mod^{\simeq}_A}\{N\}.
\]
Let $\End(N)\in \Alg_1(\Mod_A)$ be the endomorphism algebra of $N$.
We write $\End^L(N)$ for the dg Lie algebra associated to $\End^L(N)$.
 Let $\FF_{\End^L(N)}$ denote the formal stack associated to
 $\End^L(N)$  (through the equivalence $Lie_A\simeq \FST_A$). 

There exists a canonical morphism
 $J_{N}:\Def(N)\to \FF_{\End^L(N)}$ in $\Fun(\EXT_A,\SSS)$
 such that $\Def(N)(R)\to \FF_{\End^L(N)}(R)$ is a (fully faithful) functor between $\infty$-groupoids/spaces for any $[R\to A]\in \EXT_A$ (cf. \cite{IMA}, \cite{DAG}).
  Here is a quick review of $J_N$.
For any $L\in Lie_A$, we have
\[
\Map_{\Fun(\EXT_A,\SSS)}(\FF_L,\FF_{\End^L(N)})\simeq \Map_{Lie_A}(L,\End^L(N))\simeq \Rep(L)(\Mod_A)^{\simeq}\times_{\Mod_A^{\simeq}}\{N\}.
\]
Through this equivalence, the identity $\FF_{\End^L(N)}\to \FF_{\End^L(N)}$ corresponds to the tautological action of $\End^L(N)$
on $N$.
If $X$ is a functor $\EXT_A\to \SSS$,
it follows from the definition of $\Rep_H(X)$ that 
\[
\Map_{\Fun(\EXT_A,\SSS)}(X,\FF_{\End^L(N)})\simeq \Rep_H(X)^{\simeq}\times_{\Mod^{\simeq}_A}\{N\}.
\]
That is, $\FF_{\End^L(N)}$ represents the functor
$\Rep_H(-)^{\simeq}\times_{\Mod^{\simeq}_A}\{N\}:\Fun(\EXT_A,\SSS)^{op}\to \wSSS$
given by $X\mapsto  \Rep_H(X)^{\simeq}\times_{\Mod^{\simeq}_A}\{N\}$.
On the other hand, $\Def(N)$ represents the functor
$\QC_H(-)^{\simeq}\times_{\Mod^{\simeq}_A}\{N\}:\Fun(\EXT_A,\SSS)^{op}\to \wSSS$
given by $X\mapsto  \QC_H(X)^{\simeq}\times_{\Mod^{\simeq}_A}\{N\}$.
The morphism $\QC_H\to \Rep_H$ in $\Fun(\Fun(\EXT_A,\SSS)^{op}, \wCat)$ (cf. Section~\ref{sheaves}) induces $\QC_H(-)^{\simeq}\times_{\Mod^{\simeq}_A}\{N\}\to \Rep_H(-)^{\simeq}\times_{\Mod^{\simeq}_A}\{N\}$.
By Yoneda lemma, it gives rise $J_N:\Def(N)\to \FF_{\End^L(N)}$.

Consider 
the functor
$M_{\CCC}^{\textup{plain}}:\Def_{\CCC}\rightarrow \Def(\HH_\bullet(\CCC/A))$
which carries a deformation $(\CCC_R,\CCC_R\otimes_{\Perf_R}\Perf_A\simeq \CCC)$ of $\CCC$ to the deformation $(\HH_\bullet(\CCC_R/R),\alpha:\HH_\bullet(\CCC_R/R)\otimes_RA\simeq \HH_\bullet(\CCC/A))$
of $\HH_\bullet(\CCC/A)$, where $\HH_\bullet(\CCC_R/R)$ is the
relative Hochshchild homology.
By \cite[Theorem 8.23]{IMA}, there exists a diagram in $\Fun(\EXT_A,\SSS)$:
\[
\xymatrix{
\FF_{\mathbb{T}_{A/k}[-1]}\simeq \widehat{S\times_kS} \ar[r]^(0.65){\overline{KS}_{\CCC}} \ar[rd] & \Def_{\CCC} \ar[r]^(0.4){M_{\CCC}} \ar[d]^{J_{\CCC}} & \Def^{S^1}(\HH_\bullet(\CCC/A)) \ar[d]^{J_{\HH_\bullet(\CCC/A)}^{S^1}} \ar[r] & \Def(\HH_\bullet(\CCC/A))   \ar[d]^{J_{\HH_\bullet(\CCC/A)}}  \\
   &     \FF_{\GG_{\CCC}} \ar[r] & \FF_{\End^L(\HH_\bullet(\CCC/A))^{S^1}}  \ar[r] & \FF_{\End^L(\HH_\bullet(\CCC/A))} \\
}
\]
whose squares and triangle commute up to canonical homotopy.
Let us briefly explain this diagram.
The morphism
 $\FF_{\mathbb{T}_{A/k}[-1]}\to \FF_{\GG_{\CCC}}$
is induced by the Kodaira-Spencer morphism $KS_{\CCC}$.
The functor $\Def^{S^1}(\HH_\bullet(\CCC/A)):\EXT_A\to \SSS$
is given by $[R\to A] \mapsto (\Mod_R^{S^1})^\simeq \times_{(\Mod_A^{S^1})^{\simeq}}\{\HH_\bullet(\CCC/A)\}$.
The map  $\ast=\{\HH_\bullet(\CCC/A)\}\to \Mod_A^{S^1}=\Fun(BS^1,\Mod_A)$
is determined by the Hochschild homology with the $S^1$-action.
Namely, it describe the $S^1$-equivariant deformations of $\HH_\bullet(\CCC/A)$.
Informally, $M_{\CCC}:\Def_{\CCC}  \to \Def^{S^1}(\HH_\bullet(\CCC/A))$
can be descirbed as follows: for any $R\in \EXT_A$,
$\Def_{\CCC}(R)  \to \Def^{S^1}(\HH_\bullet(\CCC/A))(R)\to  \Def(\HH_\bullet(\CCC/A))(R)$
is naturally equivalent to
\[
\ST_R^{\simeq}\times_{\ST_A^{\simeq}}\{\CCC\} \to (\Mod_R^{S^1})^\simeq \times_{(\Mod_A^{S^1})^{\simeq}}\{\HH_\bullet(\CCC/A)\} \stackrel{\textup{forget}}{\longrightarrow}  \Mod_R^{\simeq}\times_{\Mod_A^{\simeq}}\{\HH_\bullet(\CCC/A)\} 
\]
where the first arrow is induced by $\HH_\bullet(-/R):\ST_R\to \Mod_R^{S^1}$ and $\HH_\bullet(-/A):\ST_A\to\Mod_A^{S^1}$.
Note that $A^L_{\CCC}:\GG_{\CCC}\to \End^L(\HH_\bullet(\CCC/A))$
is a morphism in $Lie_A^{S^1}$, where the $S^1$-action on $\GG_{\CCC}$ is trivial. 
It gives rise to
\[
\GG_{\CCC}\to \End^L(\HH_\bullet(\CCC/A))^{S^1}\to \End^L(\HH_\bullet(\CCC/A))
\]
in $Lie_A$ where $\End^L(\HH_\bullet(\CCC/A))^{S^1}$ denotes the homotopy fixed points
of the $S^1$-action on $\End^L(\HH_\bullet(\CCC/A))$
(it is not the cotensor by $S^1$).
Passing to formal stacks, we have 
$\FF_{\GG_{\CCC}} \to \FF_{\End^L(\HH_\bullet(\CCC/A))^{S^1}}  \to F_{\End^L(\HH_\bullet(\CCC/A))}$. 

\begin{Construction}
\label{firstobject}
Applying $\QC_H$ to the above diagram we obtain the following diagram 
\[
\xymatrix{
\QC_H(\widehat{S\times_k S})   &   \QC_H(\Def_{\CCC})  \ar[l] & \QC_H(\Def^{S^1}(\HH_\bullet(\CCC/A))) \ar[l]  & \QC_H(\Def(\HH_\bullet(\CCC/A))) \ar[l] \\
   &   \QC_H(\FF_{\GG_{\CCC}}) \ar[ul]  \ar[u] &  \QC_H(\FF_{\End^L(\HH_\bullet(\CCC/A))^{S^1}}) \ar[u] \ar[l] &   \QC_H(\FF_{\End(\HH_\bullet(\CCC/A))}) \ar[u] \ar[l] \\
}
\]
where each morphism is a pullback functor.
If we replace $\QC_H$ with $\Rep_H$, then we have a similar diagram.

By Lemma~\ref{replimit} and Remark~\ref{liepull},
the functors
\[
 \QC_H(\FF_{\End^L(\HH_\bullet(\CCC/A))}) \to \QC_H(\FF_{\GG_{\CCC}}) \to  \QC_H(\FF_{\mathbb{T}_{A/k}[-1]})\simeq\QC_H(\widehat{S\times S})
\]
is extended to
\[
 \Rep(\End^L(\HH_\bullet(\CCC/A)))(\Mod_A)  \to \Rep(\GG_{\CCC})(\Mod_A)\to \Rep(\mathbb{T}_{A/k}[-1])(\Mod_A)
\]
where functors are induced by the restriction of $\mathbb{T}_{A/k}[-1]\to \GG_{\CCC}\to \End^L(\HH_\bullet(\CCC/A))$
(cf. Construction~\ref{const2}).
We note that the composition $\mathbb{T}_{A/k}[-1]\to \GG_{\CCC}\to \End^L(\HH_\bullet(\CCC/A))$
factors as the sequence of $S^1$-equivariant morphisms $\mathbb{T}_{A/k}[-1] \to \mathbb{T}_{A/k}[-1]^{S^1}\to \End^L(\HH_\bullet(\CCC/A))$.
Moreover, $\widehat{S\times_k LS}$ corresponds to $\mathbb{T}_{A/k}[-1]^{S^1}$ (cf. Section~\ref{startlie}).
Thus, we have pullback functors
\[
\QC_H(\FF_{\End^L(\HH_\bullet(\CCC/A))})\to \QC_H(\widehat{S\times_kLS})\to \QC_H(\widehat{S\times_kS}), \]
and the restriction functors
\[
\Rep(\End^L(\HH_\bullet(\CCC/A)))(\Mod_A)\to \Rep(\mathbb{T}_{A/k}[-1]^{S^1})(\Mod_A)\to \Rep(\mathbb{T}_{A/k}[-1])(\Mod_A)
\]
where the first sequence can be contained in the second sequence.

Let $\mathcal{U}$ be the universal deformation of $\HH_\bullet(\CCC/A)$.
That is, $\mathcal{U}$ is the tautological object in
\[
\QC_H(\Def(\HH_\bullet(\CCC/A)))\times_{\Mod_A}\{\HH_\bullet(\CCC/A)\}
\]
which corresponds to
the identity map $\Def(\HH_\bullet(\CCC/A))\to \Def(\HH_\bullet(\CCC/A))$
through the equivalence
$\Map_{\Fun(\EXT_A,\SSS)}(\Def(\HH_\bullet(\CCC/A)),\Def(\HH_\bullet(\CCC/A)))\simeq \QC_H(\Def(\HH_\bullet(\CCC/A)))^{\simeq}\times_{\Mod^{\simeq}_A}\{\HH_\bullet(\CCC/A)\}$.
(Informally, the image of $\mathcal{U}$ in $\QC_H(\Def(\HH_\bullet(\CCC/A)))$ is an object
which associates to each $f:\Spec R\to \Def(\HH_\bullet(\CCC/A))$
the $R$-module corresponding to $f$.)

Let $\mathcal{V}$ denote  the image of $\mathcal{U}$ under $\QC_H(\Def(\HH_{\bullet}(\CCC/A)))\to \QC_H(\widehat{S\times_kS})$.

\end{Construction}

By the definition of $\overline{KS}_{\CCC}$ classified by $\widehat{\textup{pr}_2^*(\CCC)}$ over $\widehat{S\times_k S}$
and the construction of $M^{\textup{plain}}_{\CCC}$ induced by the relative Hochschild homology functor, we see:

\begin{Lemma}
\label{projectionidentify}
The object $\mathcal{V}\in \QC_H(\widehat{S\times_kS})$
is obtained from the base change $\HH_\bullet(\CCC/A)\otimes_A(A\otimes_kA)=\textup{pr}_2^*(\HH_\bullet(\CCC/A))
\in \Mod_{A\otimes_kA}=:\QC(S\times_kS)$
along the second projection $\textup{pr}_2:S\times_kS\to S$.
That is, $\mathcal{V}$ is the image of $\textup{pr}_2^*(\HH_\bullet(\CCC/A))$ through the canonical map $\QC(S\times_kS)\to \lim_{\Spec R\to S\times_kS}\Mod_{R}=\QC_H(\widehat{S\times_kS})$.
\end{Lemma}

Next, we observe:

\begin{Proposition}
\label{looppromotion}
The object $\mathcal{V} \in \QC_H(\widehat{S\times_kS})$
is promoted to an object $\mathcal{V}_L$ of $\Rep_H(\widehat{S\times_kLS})$.
Moreover, $\mathcal{V}_L$ is promoted to $\Rep_H(\widehat{S\times_kLS})^{S^1}$.
Here the $S^1$-action $\Rep_H(\widehat{S\times_kLS})$ comes from the action on $LS$.
\end{Proposition}

\begin{Lemma}
\label{universalmodule}
The universal deformation $\mathcal{U}\in \QC_H(\Def(\HH_\bullet(\CCC/A)))\times_{\Mod_A}\{\HH_\bullet(\CCC/A)\}$ 
is naturally equivalent to the image of the
$\End(\HH_\bullet(\CCC/A))$-module $\HH_\bullet(\CCC/A)$
(defined by the tautological action) under the pullback along $J_{\HH_\bullet(\CCC/A)}$
\[
\Rep(\End(\HH_\bullet(\CCC/A)))\times_{\Mod_A}\{\HH_\bullet(\CCC/A)\}\to \Rep_H(\Def(\HH_\bullet(\CCC/A)))\times_{\Mod_A}\{\HH_\bullet(\CCC/A)\}.
\]
\end{Lemma}

{\it Proof of Lemma~\ref{universalmodule}}.
If $X$ is a functor $\EXT_A\to \SSS$,
\[
\Map_{\Fun(\EXT_A,\SSS)}(X,\FF_{\End^L(\HH_\bullet(\CCC/A))})\simeq \Rep_H(X)^{\simeq}\times_{\Mod^{\simeq}_A}\{\HH_\bullet(\CCC/A)\}.
\]
By the definition of $J_{\HH_\bullet(\CCC/A)}$
(see the second paragraph of Section~\ref{somediagrams}), the canonical morphism $J_{\HH_\bullet(\CCC/A)}:\Def(\HH_\bullet(\CCC/A))\to \FF_{\End^L(\HH_\bullet(\CCC/A))}$
is classified by the universal deformation
\[
\mathcal{U}\in \QC_H(\Def(\HH_\bullet(\CCC/A)))\times_{\Mod_A}\{\HH_\bullet(\CCC/A)\}\subset \Rep_H(\Def(\HH_\bullet(\CCC/A))))\times_{\Mod_A}\{\HH_\bullet(\CCC/A)\}.
\]
Thus, our claim follows.
\QED

{\it Proof of Proposition~\ref{looppromotion}.}
Note that $\mathbb{T}_{A/k}[-1]\to \End^L(\HH_\bullet(\CCC/A))$
factors as $\mathbb{T}_{A/k}[-1]\to\mathbb{T}_{A/k}[-1]^{S^1}\to  \End^L(\HH_\bullet(\CCC/A))$. Thus,
$\HH_\bullet(\CCC/A)$ endowed with the tautological action of $\End(\HH_\bullet(\CCC/A))$
gives rise to a $\mathbb{T}_{A/k}[-1]^{S^1}$-module $\HH_\bullet(\CCC/A)$, that is,
an object $\mathcal{V}_L$ of $\Rep(\mathbb{T}_{A/k}[-1]^{S^1})(\Mod_A)$ whose underlying module is $\HH_\bullet(\CCC/A)$.
We have the diagram
\[
\xymatrix{
\mathcal{F}_{\TT_{A/k}[-1]}\simeq \widehat{S\times_kS} \ar[r] \ar[d] & \mathcal{F}_{\TT_{A/k}[-1]^{S^1}}\simeq \widehat{S\times_kLS} \ar[d] \ar[rd] &      \\
 \mathcal{F}_{\GG_{\CCC}} \ar[r] & \mathcal{F}_{\GG_{\CCC}^{S^1}} \ar[r]  &
\mathcal{F}_{\End^L(\HH_\bullet(\CCC/A))}.
}
\]
Lemma~\ref{universalmodule} shows that the image of $\mathcal{V}_L$
in $\Rep_H(\widehat{S\times_kS})$
is naturally equivalent to $\mathcal{V}$
because $\widehat{S\times_kS}\to \mathcal{F}_{\End^L(\HH_\bullet(\CCC/A))}$
factors through $\Def(\HH_\bullet(\CCC/A))\to \mathcal{F}_{\End^L(\HH_\bullet(\CCC/A))}$.
Remember that $\mathbb{T}_{A/k}[-1]^{S^1}\to \End^L(\HH_\bullet(\CCC/A))$
is promoted to an $S^1$-equivariant morphism of dg Lie algebras.
Thus, $\mathcal{V}_L$ is promoted to $(\Rep(\mathbb{T}_{A/k}[-1]^{S^1})(\Mod_A))^{S^1}
\simeq \LMod_{U_1(\mathbb{T}_{A/k}[-1]^{S^1})}(\Mod_A^{S^1})$
where 
$\mathbb{T}_{A/k}[-1]^{S^1}$ is regarded as an object of $Lie_A^{S^1}=\Fun(BS^1,Lie_A)$.
Note that the $S^1$-action on $\widehat{S\times_kLS}$ amounts to the canonical $S^1$-action
on $\mathbb{T}_{A/k}[-1]^{S^1}$ since both $\widehat{S\times_kLS}$ and $\mathbb{T}_{A/k}[-1]^{S^1}$
are obtained by the cotensor by $S^1$.
Hence we see that $\mathcal{V}_L$ belongs to $\Rep_H(\widehat{S\times_kLS})^{S^1}$.
\QED

In other words, Proposition~\ref{looppromotion}
goes as follows (see also Construction~\ref{firstobject}).
The object $\mathcal{V}_L\in \Rep_H(\widehat{S\times_kLS})\simeq \Rep(\mathbb{T}_{A/k}[-1]^{S^1})(\Mod_A)$ and $\mathcal{V}\in \QC_H(\widehat{S\times_kS})\subset \Rep(\mathbb{T}_{A/k}[-1])(\Mod_A)$
are obtained from the $S^1$-equivariant morphisms of dg Lie algebras 
\[
\mathbb{T}_{A/k}[-1]\to  \mathbb{T}_{A/k}[-1]^{S^1}\to \End^L(\HH_\bullet(\CCC/A))
\]
and the tautological action $\End^L(\HH_\bullet(\CCC/A))$ on $\HH_\bullet(\CCC/A)$
 (note that this tautological action defines an object of $\Rep(\End^L(\HH_\bullet(\CCC/A)))(\Mod_A^{S^1})$).
Consequently, $\mathcal{V}_L$ and $\mathcal{V}$ are promoted to objects in
\[
\mathcal{V}_L'\in \Rep(\mathbb{T}_{A/k}[-1]^{S^1})(\Mod_A^{S^1})\ \ \ \ \textup{and} \ \ \ \ \ \mathcal{V}'\in \QC_H(\widehat{S\times_kS})^{S^1}\subset \Rep(\mathbb{T}_{A/k}[-1])(\Mod_A^{S^1})
\]
such that $\mathcal{V}_L'$ maps to $\mathcal{V}'$.
Namely, $\mathcal{V}_L'$ and $\mathcal{V}'$ correspond to
the canoncial $\mathbb{T}_{A/k}[-1]^{S^1}$-module and the canonical $\mathbb{T}_{A/k}[-1]$-module,
respectively
(cf. Defintion~\ref{canonicaltangentaction}).

\begin{Proposition}
\label{naiveequivariant}
There exists an equivalence $\mathcal{V}'\simeq \widehat{\textup{pr}_2^*(\HH_\bullet(\CCC/A))}$ in 
$\QC_H(\widehat{S\times_kS})^{S^1}\subset \Rep(\mathbb{T}_{A/k}[-1])(\Mod_A^{S^1})$.
Combined with the data $\mathcal{V}_L'\mapsto \mathcal{V}'$,
it determines an object $\mathcal{V}_\dagger$ of
\[
\QC(S)^{S^1}\times_{\Rep(\mathbb{T}_{A/k}[-1])(\Mod_A^{S^1})} \Rep(\mathbb{T}_{A/k}[-1]^{S^1})(\Mod_A^{S^1}).
\]
where $\QC(S)^{S^1}\to \Rep(\mathbb{T}_{A/k}[-1])(\Mod_A^{S^1})$
is determined by the pullback along the second projection
$\QC(S)\to \QC_H(\widehat{S\times S})\subset \Rep(\mathbb{T}_{A/k}[-1])(\Mod_A)$.
\end{Proposition}

\Proof
Note that $\mathcal{V}'\in \Rep(\mathbb{T}_{A/k}[-1])(\Mod_A^{S^1})\simeq\Fun(BS^1, \Rep(\mathbb{T}_{A/k}[-1])(\Mod_A))$
is obtained from  the tautological action of $\End(\HH_\bullet(\CCC/A))$ on $\HH_\bullet(\CCC/A)$
and the $S^1$-equivariant morphism of dg Lie algebras
\[
\mathbb{T}_{A/k}[-1]\to \End^L(\HH_\bullet(\CCC/A))^{S^1}\to  \End^L(\HH_\bullet(\CCC/A))
\]
where the actions on the left and middle ones are trivial.
Set
\[
\QC_H(\Def^{S^1}(\HH_\bullet(\CCC/A))^{S^1}:=\Fun(BS^1,\QC_H(\Def^{S^1}(\HH_\bullet(\CCC/A)))).
\]
Let $\mathcal{U}'\in \QC_H(\Def^{S^1}(\HH_\bullet(\CCC/A)))^{S^1}\times_{\Mod_A^{S^1}}\{\HH_\bullet(\CCC/A)\}$
be the universal $S^1$-equivariant deformation of $\HH_\bullet(\CCC/A)$
(using the canonical equivalence
\[
\Map_{\Fun(\EXT_A,\SSS)}(X,\Def^{S^1}(\HH_{\bullet}(\CCC/A)))\simeq (\QC_H(X)^{S^1})^{\simeq}\times_{(\Mod_A^{S^1})^{\simeq}}\{\HH_\bullet(\CCC/A)\}
\]
we define $\mathcal{U}'$ in a similar way to $\mathcal{U}$).
By abuse of notation, we write $\mathcal{U}'$ for the image of $\mathcal{U}'$ in
$\QC_H(\Def^{S^1}(\HH_\bullet(\CCC/A)))^{S^1}$.
The morphism of dg Lie algebras
$\End^L(\HH_\bullet(\CCC/A))^{S^1}\to  \End^L(\HH_\bullet(\CCC/A))$
gives us an $\End^L(\HH_\bullet(\CCC/A))^{S^1}$-module $\HH_\bullet(\CCC/A)$ endowed with
the $S^1$-action, which we think of as an object $\mathcal{R}$ of $\Rep(\End^L(\HH_\bullet(\CCC/A))^{S^1})(\Mod_A^{S^1})\simeq (\Rep_H(\mathcal{F}_{\End^L(\HH_\bullet(\CCC/A))^{S^1}})(\Mod_A))^{S^1}$.
The pullback along $J_{\HH_\bullet(\CCC/A)}^{S^1}$
\[
(\Rep_H(\mathcal{F}_{\End^L(\HH_\bullet(\CCC/A))^{S^1}})(\Mod_A))^{S^1}
\to \Rep_H(\Def^{S^1}(\HH_\bullet(\CCC/A))))^{S^1}\supset \QC_H(\Def^{S^1}(\HH_\bullet(\CCC/A))))^{S^1} 
\]
sends $\mathcal{R}$ to $\mathcal{U}'$ (this is a simple generarization of Lemma~\ref{universalmodule} to the $S^1$-equivariant situation, see Remark~\ref{equivariantuniversalmodule}).
Consequently, $\mathcal{V}'$ is naturally equivalent to
the image of $\mathcal{U}'$ under $\QC_H(\Def^{S^1}(\HH_\bullet(\CCC/A)))^{S^1}\to \QC_H(\widehat{S\times_kS})^{S^1}$. 
Now the same argument as in Lemma~\ref{projectionidentify}
shows that $\mathcal{V}'$ is naturally equivalent to the object 
$\widehat{\textup{pr}_2^*(\HH_\bullet(\CCC/A))}$ obtained from the pullback
$\textup{pr}_2^*(\HH_\bullet(\CCC/A))\in \QC(S\times_kS)^{S^1}$ of $\HH_\bullet(\CCC/A)\in \QC(S)^{S^1}=\Mod_A^{S^1}$.
\QED

\begin{Remark}
\label{equivariantuniversalmodule}
Consider the adjoint pair $\textup{triv}:\Mod_A=\Fun(\ast,\Mod_A)\rightleftarrows \Mod_A^{S^1}=\Fun(BS^1,\Mod_A)$ where
the left adjoint induced by the composition with $BS^1\to \ast$ (so that the right adjoint is given by $(-)^{S^1}$). 
By this adjoint pair and the universal property of 
$\End(\HH_\bullet(\CCC/A))\in \Alg_{\eone}(\Mod_A^{S^1})$, we have
\begin{eqnarray*}
\Map_{Lie_A}(L, \End^L(\HH_\bullet(\CCC/A))^{S^1}) &\simeq& \Map_{\Alg_{\eone}(\Mod_A)}(U_1(L),\End(\HH_\bullet(\CCC/A))^{S^1}) \\ 
 &\simeq& \Map_{\Alg_{\eone}(\Mod_A^{S^1})}(\textup{triv}(U_1(L)),\End(\HH_\bullet(\CCC/A))) \\
&\simeq& \LMod_{\textup{triv}(U_1(L))}(\Mod_A^{S^1})^{\simeq}\times_{(\Mod_A^{S^1})^{\simeq}} \{\HH_\bullet(\CCC/A)\} \\
&\simeq& (\Rep(L)(\Mod_A)^{S^1})^{\simeq}\times_{(\Mod_A^{S^1})^{\simeq}} \{\HH_\bullet(\CCC/A)\}
\end{eqnarray*}
 for any $L\in Lie_A$.
The identity map of $\End^L(\HH_\bullet(\CCC/A))^{S^1}$
maps to the tautological action of $\End^L(\HH_\bullet(\CCC/A))^{S^1}$ on $\HH_\bullet(\CCC/A)$.
As with $\mathcal{F}_{\End(\HH_\bullet(\CCC/A)}$, there is a canonical equivalence
\[
\Map_{\Fun(\EXT_A,\SSS)}(X,\FF_{\End^L(\HH_\bullet(\CCC/A))^{S^1}}) \simeq (\Rep_H(X)^{S^1})^{\simeq}\times_{(\Mod_A^{S^1})^{\simeq}}\{\HH_\bullet(\CCC/A)\}
\]
for any $X\in \Fun(\EXT_A,\SSS)$.
The morphism $\Def^{S^1}(\HH_\bullet(\CCC/A))\to \FF_{\End(\HH_\bullet(\CCC/A))^{S^1}}$ is
determined by
\begin{eqnarray*}
\QC_H(-)^{S^1}\times_{\Mod^{S^1}_A}\{\HH_\bullet(\CCC/A)\} \to \Rep_H(-)^{S^1}\times_{\Mod^{S^1}_A}\{\HH_\bullet(\CCC/A)\}.
\end{eqnarray*} 
Thus, we see that the universal deformation $\mathcal{U}'$ can naturally be identified with
the pullback of the tautological action of $\End^L(\HH_\bullet(\CCC/A))^{S^1}$ on $\HH_\bullet(\CCC/A)$
along $\Def^{S^1}(\HH_\bullet(\CCC/A))\to \FF_{\End^L(\HH_\bullet(\CCC/A))^{S^1}}$.
\end{Remark}

\section{Extensions via Lie algebra actions}

\subsection{}
\label{Indreviewsec}

We use a theory of
 formal stacks in the formulation of Ind-coherent complexes/sheaves, which is extensively developed in Gaitsgory and Rozenblyum \cite[Vol. II]{Gai2}.
In Section~\ref{Indreviewsec}, we will give a minimal review of definitions and some results which we will use.
Suppose that $A$ is almost of finite type over $k$: $H^0(A)$ is a usual commutative algebra
of finite type over $k$, and each $H^i(A)$ is a finitely generated $H^0(A)$-module.
(By derived Hilbert basis theorem 
\cite[7.2.4.31]{HA}, the condition of almost of finite type over $k$ is
equivalent to the condition of almost of finite presentation over $k$ in \cite[Definition 7.2.4.26]{HA}.)
Furthermore, assume that $A$ is eventually coconnective, that is, there is a nonpositive integer
$n$ such that $H^{i}(A)=0$ for $i<n$.

{\it $!$-Formal stacks.}
Let $\CCAlgft_k$ be the full subcategory of $\CCAlg_k$ spanned by those objects
which are almost of finite type over $k$.
Let $\CCAlgftec_k$ be the full subcategory of $\CCAlg_k$ spanned by those objects
which are eventually coconnective and almost of finite type over $k$.
For $R\in \CCAlg_k$ we let $R_{red}$ denote the reduced ring $H^0(R)_{red}$.
Let $\BEXT_A$ be the full subcategory of $(\CCAlgftec_k)_{A//A}$
that consists of objects of the form $A\stackrel{u}{\to} R\stackrel{v}{\to} A$ such that 
the associated ring homomorphisms of usual reduced rings
$u_{red}:A_{red}\to R_{red}$ and $v_{red}:R_{red}\to A_{red}$ are isomorphisms.
We let $\NIL_A$ or $\NIL_{S}$ denote the opposite category of $\BEXT_A$.
Since $A$ is eventually coconnective, 
there are natural fully faithful inclusions $\EXT_A\subset \BEXT_A$ and $\TSZ_A\subset \NIL_A$.
We define the $\infty$-category $\GFST_A$ to be the full subcategory of $\Fun(\BEXT_A,\SSS)$ which consists of
those functors $F:\BEXT_A\to \SSS$
satisfying the following conditions (see \cite[Vol. II, Chap.5, 1.5]{Gai2} for the detailed account):
\begin{itemize}
\item $F(A)$ is a contractible space,

\item If $T_1\sqcup_{T}T'$ is a pushout in $\NIL_S$ such that $T\to T'$ is a square zero extension,
then the canonical morphism $F(T_1\sqcup_{T}T')\to F(T_1)\times_{F(T)}F(T')$ is an equivalence.

\end{itemize}
We will dub an object of $\GFST_A$ as a pointed $!$-formal stack over $A$ (or $S=\Spec A$)
(in \cite{Gai2}, objects of $\GFST$  are referred to as pointed formal moduli problems).
The Yoneda embedding $\NIL_A\to \GFST_A$ exhibits $\NIL_A$ as a full subcategory.

\begin{Remark}
The $\infty$-category $\GFST_A$ is one of the formulations of the $\infty$-category of 
pointed formal moduli
problems over $\Spec A$ in \cite{Gai2}.
See \cite[Ch.5, 1.2.2, 1.4.2, 1.5.2]{Gai2} for several formulations.
\end{Remark}

{\it Ind-coherent complexes.}
We review some definitions concerning Ind-coherent complexes/sheaves.
The following is based on \cite{Ind}, \cite{Gai2}.
Suppose that $B\in \CCAlg_k$ is almost of finite type $k$. 
Let $\Coh(B)$ be the full subcategory of $\Mod_B=\QC(B)$ spanned by
bounded complexes with coherent cohomologies. 
We define the stable $\infty$-category of Ind-coherent complexes over $\Spec B$
to be $\QC_!(B):=\Ind(\Coh(B))$.
The stable presentable $\infty$-categories $\QC_!(B)$ have a functoriality
given by $!$-pullbacks:
a morphism $f:\Spec B'\to \Spec B$ induces a colimit-preserving functor
$f^!:\QC_!(\Spec B)=\QC_!(B)\to \QC_!(B')=\QC_!(\Spec B')$.
Let $f_*:\Coh(B')\to \QC^+(B)$ be the pushforward determined by the restriction along $B\to B'$.
Here
$\QC^+(B)\subset \QC(B)$ is the full subcategory spanned by left-bounded objects
with respect to the standard $t$-structure.
Since $\QC^+(B)\simeq \Ind(\Coh(B))^+\subset \QC_!(B)$, we have $\Coh(B')\to \QC^+(B)\simeq \Ind(\Coh(B))^+$.
Passing to the Ind-category, it gives rise to a colimit-preserving functor
$f_*^{\textup{IndCoh}}:\QC_!(B')=\Ind(\Coh(B'))\to \QC_!(B)=\Ind(\Coh(B))$
which extends $\Coh(B')\to \Ind(\Coh(B))$ in an essentially unique way.
When $f$ is proper (the induced morphism of classical schemes is proper),
$f^!$ is defined to be a right adjoint of  $f_*^{\textup{IndCoh}}$.
When $f$ is an open immersion, $f^!$ is defined to be a left adjoint of $f_*^{\textup{IndCoh}}$.
If $f$ is arbitrary, we decompose $f$ as $j\circ p$ such that $j$ is an open immersion and $p$
is proper,
and set $f^!=p^!\circ j^!$.
It is necessary to prove that the definition is independent of the choice of 
decomposition.
It is a difficult task to give a functorial and canonical definition of $\QC_!(B)$.
This was achieved in \cite{Gai2}. In particular, from \cite{Gai2} we have
a functor
\[
\QC_!:\CCAlgft_{k}\longrightarrow \PR_k
\]
which carrris $B$ to $\QC_!(B)$, and carries $f$ to $f^!$ (in {\it loc. cit.}, the symbol $\textup{IndCoh}_{\textup{Sch}_{\textup{aft}}^{\textup{aff}}}$ is used).
There is a functor
\[
\varUpsilon_B: \Mod_B=\QC(B)\longrightarrow \QC_!(B)
\]
which sends $M$ to $M\otimes_B \omega_B$
where $\omega_B$ a dualizing object given by $\omega_B=p^!(k)$ where
$p:\Spec B\to \Spec k$ is the structure morphism, see \cite[Vol. I]{Gai2}, \cite[Section 5]{Ind}.
When $B$ is eventually cocconnective, $\varUpsilon_B$ is fully faithful.
When $B$ is smooth over $k$, $\varUpsilon_B$ is an equivalence.
The functor $\varUpsilon_B$ is functorial with respect to $B$:
if $\QC: \CCAlgft_{k}\longrightarrow \PR_k$ denotes the functor
given by $B\mapsto \Mod_B=\QC(B)$,
there is $\Upsilon:\QC\to \QC_!$ between functors $\CCAlgft_{k}\longrightarrow \PR_k$.
In particular, $f:\Spec B'\to \Spec B$ induces the commutative diagram in $\PR_k$:
\[
\xymatrix{
\QC(B)  \ar[r]^{\Upsilon_B}\ar[d]^{f^*}  & \QC_!(B) \ar[d]^{f^!}  \\
\QC(B') \ar[r]^{\Upsilon_{B'}}   &   \QC_!(B').
}
\]
The stable $\infty$-category $\QC_!(B)$ admits a symmetric monoidal structure,
and $\Upsilon_B$ is promoted to a symmetric monoidal functor.
Moreover, $\QC_!:\CCAlgft_{k}\longrightarrow \PR_k$ is promoted to
$\CCAlgft_{k}\rightarrow \CAlg(\PR_k)$.

Let $\Fun(\CCAlg_k,\SSS)_{\textup{laft}}\subset\Fun(\CCAlg_k,\SSS)$
be the full subcategory of functors/prestacks locally almost of finite type over $k$, see  \cite[Vol.I, Chap.2, 1.7]{Gai2}.
Here we do not recall the definition of the condition of locally almost of finite type,
but the restriction functor along $\CCAlgftec_{k} \to \CCAlg_k$
induces an equivalence $\Fun(\CCAlg_k,\SSS)_{\textup{laft}}\stackrel{\sim}{\to} \Fun(\CCAlgftec_k,\SSS)$.
For example, a derived affine scheme almost of finite type is an example of a prestack
locally almost of finite type over
$k$. Thus, $\Spec B$ such that $B$ is
almost of finite type, is completely determined by the functor $\CCAlgftec_k\to \SSS$
defined by the formula $R\mapsto \Map_{\CCAlg_k}(B,R)$.
We define $\Fun(\CCAlgftec_k,\SSS)^{op}\to \PR_k$ as a right Kan extention of
\[
\CCAlgftec_k\hookrightarrow \CCAlgft_k\stackrel{\QC_!}{\longrightarrow} \PR_k
\]
along the Yoneda embedding $\CCAlgftec_k\to \Fun(\CCAlgftec_k,\SSS)^{op}$.
We abuse notation by writing $\QC_!$ for the resulting functor
$\Fun(\CCAlgftec_k,\SSS)^{op}\to \PR_k$.
We refer to $\QC_!(F)$ as the stable $\infty$-category of Ind-coherent complexes on $F$.

According to \cite[Vol.II, Chap.7, 3.1.4, 3.3.2]{Gai2}, there is a categorical equivalence
\[
Lie_A^!\simeq \GFST_A
\]
where $Lie_A^!$ is the $\infty$-category of Lie algebra objects $\Lie(\QC_!(A))$.
The symmetric monoidal functor $\varUpsilon_A$ induces $Lie_A \longrightarrow Lie_A^!$.
Since we assume that $A$ is eventually coconnective, $\varUpsilon_A$ is a fully faithful left adjoint functor, see \cite[Section 9.6]{Ind}.
Hence
$Lie_A\to Lie_A^!$ is also a fully faithful left adjoint functor.
From Section~\ref{zigzagconstructionsec}, we will assume that $A$ is smooth over $k$
so that in the forthcoming sections we consider the case
when $Lie_A\to Lie_A^!$ induced by $\varUpsilon_A$ is an equivalence.

There exists a fully faithful functor
$\GFST_A\hookrightarrow (\Fun(\CCAlg_k,\SSS)_{\textup{laft}})_{S//S}$
(actually, in \cite[Vol.II, Chapter 5]{Gai2}, $\GFST_A$ is defined as a full subcategory of
 $(\Fun(\CCAlg_k,\SSS)_{\textup{laft}})_{S//S}$).
 We define
\[
\QC_!|_{\GFST_A}:(\GFST_A)^{op}\hookrightarrow (\Fun(\CCAlg_k,\SSS)_{\textup{laft}})^{op}_{S//S}\stackrel{\textup{forget}}{\longrightarrow} \Fun(\CCAlg_k,\SSS)^{op}_{\textup{laft}}\stackrel{\QC_!}{\to} \PR_k.
\]
For $F\in \GFST_A$, $\QC_!(F)$ is defined as the image of $F$ under the composite. 
Suppose that $\colim_{i\in I}F_i=F$ is a colimit of a sifted diagram $I\to \GFST_A$.
By the compatibility of $\QC_!$ with sifted colimits in $\GFST_A$ (\cite[Vol.I, Chap.7, 5.3]{Gai2}),  
the diagram of $!$-pullback functors yields an equivalence of symmetric monoidal $\infty$-categories
$\QC_!(F)\stackrel{\sim}{\to}  \lim_{i\in I}\QC_!(F_i)$.
In particular, if each $F_i$ is representable by $\Spec B_i$,
we have
$\QC_!(F)\simeq \lim_{i\in I} \QC_!(B_i)$.

\subsection{}

\label{compareformal}

In Section~\ref{compareformal},
we prove some results
concerning comparisons of
pointed formal stacks and pointed $!$-formal stacks.

\begin{Construction}
\label{transformalstack}
We will construct
$\varTheta_A:\FST_A\rightarrow \GFST_A$
that extends the fully faithful embedding $\TSZ_A\hookrightarrow \NIL_A \subset \GFST_A$
induced by the Yoneda embedding.
Recall the definition $\FST_A=\PPP_{\Sigma}^{st}(\TSZ_A)$ and the universal property
of $\PPP_{\Sigma}^{st}(\TSZ_A)$ (cf. Section~\ref{stablederived} and Section~\ref{formalstack}).
 To construct $\FST_A\rightarrow \GFST_A$, it will suffice to construct
 $\TSZ_A\to \GFST_A$
 which preserve
finite coproducts and carry each morphism in $T=\{\Spec A\sqcup_{\Spec(A\oplus M)}\Spec A\to \Spec (A\oplus M[-1])\}_{A\oplus M[-1]\in \EXT_A}$ to an equivalence.
Consider the Yoneda embedding $\TSZ_A\subset \NIL_A\hookrightarrow \GFST_A$.
This functor preserves finite coproducts (by the definition \cite[Ch. 5, 1.5.2 (b)]{Gai2}).
Moreover, it carries each morphism in $T$ to an equivalence since (pointed) objects lying in $\GFST_A$
admit deformation theory in the sense of \cite[Ch.1, 7.1.2]{Gai2} (in particular,
objects in $\GFST_A$ have cotangent complexes). 
Consequently, $\TSZ_A\subset \NIL_A\hookrightarrow \GFST_A$
is (uniquely) extended to a colimit-preserving functor $\FST_A\rightarrow \GFST_A$,
that we will denote by $\varTheta_A$.
\end{Construction}

\begin{Lemma}
\label{fullyfaithfullemma}
The functor $\varTheta_A:\FST_A\rightarrow \GFST_A$ 
is fully faithful.
The composite $Lie_A\simeq \FST_A \rightarrow \GFST_A\simeq Lie_A^!$
carries $\Free_{Lie}(N)\in Lie_A^f$ to $\varUpsilon_A(\Free_{Lie}(N))$.
See Section~\ref{formalstack} for $Lie_A^f$.
\end{Lemma}

\Proof
We first prove the second assertion.
We note that the restriction of $\FST_A\simeq Lie_A$
determines an equivalence $\TSZ_A\simeq Lie_A^f$.
If $M$ is a connective $A$-module of the form $M=\oplus_{1\le i \le  n} A^{\oplus r_i}[p_i]$ 
($p_i\ge0$),
$\TSZ_A\simeq Lie_A^f$ sends $\Spec (A\oplus M)$ to $\Free_{Lie}(N)$ where $N=M^\vee[-1]$ ($M^\vee$ is
the dual object of $M$ in the symmetric monoidal $\infty$-category $\Mod_A$).
According to \cite[Vol. II, Chap.7, 3.7.1 (3.13), 3.7.10]{Gai2}, the equivalence
$\GFST_A \stackrel{\sim}{\to} Lie^!_A$ carries $\Spec (A\oplus M)$ to
$\Free_{Lie^!}(N\otimes \omega_A)\simeq \Free_{Lie^!}(\Upsilon_A(N))\simeq \Upsilon_A(\Free_{Lie}(N))$
where $\Free_{Lie^!}:\QC_!(A)\to Lie_A^!$ is the free functor, that is, a left adjoint to the forgetful functor.
Thus the second claim follows.

Next, we prove the first assertion.
To this end, it is enough to show that the composite $F:Lie_A\simeq \FST_A \to \GFST_A\simeq Lie_A^!$ is fully faithful.
Observe that $\TSZ_A \hookrightarrow \FST_A\to \GFST_A$ determines an equivalence to
its essential image which we denote by $\varTheta_A(\TSZ_A)$.
Through the equivalence $\GFST_A\simeq Lie_A^!$, 
$\varTheta_A(\TSZ_A)$ corresponds to the full subcategory of $Lie_A^!$
which consists of objects of the form $\varUpsilon_A(\Free_{Lie}(N))$ such that $\Free_{Lie}(N)\in Lie_A^f$.
We write $Lie_A^{!,f}$ for this full subcategory.
The composite $F:Lie_A \to Lie_A^!$ induces an equivalence $\theta:Lie^f_A\stackrel{\sim}{\to} Lie_A^{!,f}$.
By \cite[1.2.2]{H},
$Lie_A^f \hookrightarrow Lie_A$ determines an equivalence
$\mathcal{P}_{\Sigma}^{st}(Lie_A^f)\simeq Lie_A$
where $\mathcal{P}_{\Sigma}^{st}(Lie_A^f)\to Lie_A$ is induced by the universal property
of $\mathcal{P}_{\Sigma}^{st}(Lie_A^f)$ (see Section~\ref{stablederived}).
Recall that $\varUpsilon_A:Lie_A\to Lie_A^!$ is fully faithful.
It can be identified with
$\mathcal{P}^{st}_\Sigma(Lie_A^f){\simeq} \mathcal{P}^{st}_\Sigma(Lie_A^{!,f})\hookrightarrow Lie_A^!$
obtained from
 $Lie_A^f\stackrel{\varUpsilon_A}{\simeq}Lie_A^{!,f}\hookrightarrow Lie_A^!$.
The composite $F:Lie_A\to Lie_A^!$ is a colimit-preserving functor that extends the restrcition
$f:Lie_A^f\subset Lie_A\to Lie^!_A$ in an essentially unique way.
In particular, $F$ is given by
$\mathcal{P}^{st}_\Sigma(Lie^f_A)\stackrel{\mathcal{P}^{st}_\Sigma(\theta)}{\simeq} \mathcal{P}^{st}_\Sigma(Lie^{!,f}_A)\hookrightarrow Lie_A^!$
so that 
$F$ is fully faithful.
This proves the first assertion.
\QED

\begin{Remark}
\label{smoothremark}
When $A$ is smooth usual algebra  over $k$, the proof of Lemma~\ref{fullyfaithfullemma} shows that
$\varTheta_A:\FST_A\rightarrow \GFST_A$ is an equivalence because the inclusion $Lie_A^{!,f}\subset Lie_A^!$
induces an equivalence $\mathcal{P}_\Sigma^{st}(Lie_A^{!,f})\stackrel{\sim}{\to} Lie_A^!$.
\end{Remark}

Let $X$ and $S$ be functors $\CCAlg_k\to \SSS$
and let $i:S\to X$ be a morphism in $\Fun(\CCAlg_k,\SSS)$.
We define $X_i^\wedge$ (or simply $X_S^\wedge$) which belongs to $\Fun(\CCAlg_k,\SSS)$
(see \cite[Vol.II, Chap. 4]{Gai2}).
The functor $X_i^\wedge$ may be regarded as the formal completion of $X$
along $X$ along $i$.
As we shall observe, it is closely related to $\widehat{X}$ in Section~\ref{formalstack} in a suitable situation.
Let $red:\CCAlg_k\to \CCAlg_k$ be the functor 
given by $R\mapsto R_{red}$, that is a left adjoint of 
the inclusion $\CAlg_k^{red}\hookrightarrow \CCAlg_k$ where
$\CAlg_k^{red}$ is the ordinary category of reduced commutative $k$-algebras.
Let $X_{\DR}:\CCAlg_k\to \SSS$ be the composite functor $X\circ red$:
it is given by $R\mapsto X(R_{red})$.
We refer to $X_{\DR}$ as the de Rham prestack of $X$.
There is a natural transformation (a unit map of the adjoint pair) from the identity functor to $red$.
It gives rise to a canonical morphism $X\to X_{\DR}$.
We define $X^\wedge_i:\CCAlg_k\to \SSS$ 
to be $X\times_{X_{\DR}}S_{\DR}$.
Suppose that $S=\Spec A$ and $S\to X\to S$ is a pointed derived scheme locally almost of finite type
over $S$. (cf. \cite[Vol. I]{Gai2} for the condition of locally almost of finite type).
We regard $S\to X\times_{X_{\DR}}S_{\DR} =X_S^\wedge \to S$ 
as an object of $\Fun(\CCAlgftec_k,\SSS)_{S//S}$ (by restricting to $\CCAlgftec_k$).
The fully faithful functor $\NIL_A\hookrightarrow \Fun(\CCAlgftec_k,\SSS)_{S//S}$
induced by Yoneda embedding and
the functor $H_{S\to X_S^\wedge\to S}$ 
represented by
$S\to X_S^\wedge \to S$
define a pointed $!$-formal stack $\BEXT_A\hookrightarrow (\Fun(\CCAlgftec_k,\SSS)_{S//S})^{op}\stackrel{H_{S\to X_S^\wedge\to S}}{\longrightarrow} \SSS$.
By abuse of notation, we will write $X^\wedge_S$ also for this pointed $!$-formal stack over $A$.


The formal completion of $S\to X\to S$ in the sense of Section~\ref{formalstack} also gives rise to the pointed formal stack
$\widehat{X}\in \FST_A$. We will compare $\widehat{X}$ and $X^\wedge_S$.  
We continue to assume $S=\Spec A\in \PST$.
Let $S\to F\to S$ be an object of $\PST_{S//S}$ (i.e., a pointed functor over $S$).
Consider $\TSZ_A$ as a full subcategory of $\PST_{S//S}$,
which consists of objects of the form $\Spec A\to \Spec (A\oplus M)\to \Spec A$.
We write $h_{S\to F\to S}$ (or simply $h_F$)
for the functor $\TSZ_A^{op}=\EXT_A\to \SSS$ which is represented by $S\to F\to S$.
As above, $S\to X\to S$ is a pointed derived scheme locally almost of finite type over $S$.
We regard it as an object of  $\PST_{S//S}$.
There is a canonical morphism $X^\wedge_S\to X$ in $\PST_{S//S}$.

\begin{Lemma}
The canonical morphism $X^\wedge_S\to X$ induces
an equivalence $h_{S\to X^\wedge_S\to S}\to h_{S\to X\to S}$ in $\Fun(\EXT_A,\SSS)$.
\end{Lemma}

\Proof
The mapping space from $S_M=\Spec (A\oplus M)\in \TSZ_A$ to $X_S^\wedge$ in $\PST_{S//S}$
is given by
\begin{eqnarray*}
&\ &\{\ast\}\times_{\Map(S,X)}\Map(S_M,X)\times_{\Map(S_M,X_{\DR})}\Map(S_M,S_{\DR})\times_{\Map(S,S_{\DR})}\{\star\} \\
&\simeq& \{\ast\}\times_{\Map(S,X)}\Map(S_M,X)\times_{\Map(S_{red},X)}\Map(S_{red},S)\times_{\Map(S_{red},S)}\{\star\} \\
&\simeq& \{\ast\}\times_{\Map(S,X)}\Map(S_M,X)\times_{\Map(S_{red},X)}\{\bullet\} \\
&\simeq& \{\ast\}\times_{\Map(S,X)}\Map(S_M,X)
\end{eqnarray*}
where $\Map(-,-)$ indicates the mapping space in $\PST_{/S}$.
Here $\ast$ is the morphism $S\to X$, $\star$ is the canonical morphism $S\to S_{\DR}$,
and $\bullet$ is the composite $S_{red}\to S\to X$.
We remark that the third equivalence follows from the fact
that $\Map(S_{red},X)$ is a discrete space because $X$ is a derived scheme (or algebraic space).
We have a canonical equivalence
\[
\Map_{\PST_{S//S}}(S_M,X_S^\wedge)\simeq \Map_{\PST_{S//S}}(S_M,X)
\]
which is natural in $S_M$. Thus Lemma follows.
\QED

\begin{Corollary}
\label{formalvsgeom}
There exists an equivalence between 
\[
(\TSZ_A)_{/X_S^\wedge}\subset (\GFST_A)_{/X_S^\wedge}\ \ \ \  \textup{and} 
\ \ \  \ (\TSZ_A)_{/X}\subset \PST_{S//S}.
\]
\end{Corollary}


\begin{Proposition}
\label{formalcompatibility}
Suppose that $A$ is smooth over $k$.
Let $X$ be a derived scheme over $S$, which is endowed with
a section morphism $i:S\to X$. 
Then the functor $\varTheta_A:\FST_A\stackrel{\sim}{\to} \GFST_A$ carries $\widehat{X}$ to $X_S^\wedge$.
More precisely, there is an equivalence $\varTheta_A(\widehat{X})\to X_S^\wedge$
which is constructed in the proof.
\end{Proposition}


\Proof
For any $S_M=\Spec (A\oplus M) \in \TSZ_A$, by the definition of $\widehat{X}$
there is an equivalence
$\Map_{\FST_A}(S_M,\widehat{X})\stackrel{\sim}{\to} \Map_{\Fun(\CCAlgftec_k,\SSS)_{S//S}}(S_M,X)$.
Thus, by Corollary~\ref{formalvsgeom}  we have canoncial equivalences
\[
\Map_{\PST_{S//S}}(S_M,X_S^\wedge)\simeq \Map_{\PST_{S//S}}(S_M,X)\simeq \Map_{\FST_A}(S_M,\widehat{X}).
\]
Consequently,
$(\TSZ_A)_{/\widehat{X}}\to \FST_{A}\to \GFST_A$
is equivalent to $(\TSZ_A)_{/X^\wedge_S} \to \GFST_A$
 (by abuse of notation, we here regard $\TSZ_A$ as full subcategories of $\FST_A$ and $\GFST_A$:
 keep in mind that $\TSZ_A\subset \FST_A\to \GFST_A$ is naturally equivalent
 to the Yoneda embedding.)
Note that by definition $\varTheta_A(\widehat{X})$ is a colimit of $(\TSZ_A)_{/\widehat{X}}\to \FST_{A}\to \GFST_A$.
We may and will assume that $\varTheta_A(\widehat{X})$ is a colimit of
$(\TSZ_A)_{/X^\wedge_S} \to \GFST_A$.
The universal property of the colimit
determines a morphism $\varTheta_A(\widehat{X})\to X_S^\wedge$.
We prove that the morphism is an equivalence.
Note that the restriction of the equivalence $\GFST_A\simeq Lie_A^!$ induces 
$\TSZ_A\simeq Lie_A^{!,f}$. Let $L$ be an object of $Lie_A^!$ which corresponds
to $X_S^\wedge$.
It is enough to prove that $L$ is a colimit of $(Lie_A^{!,f})_{/L}\to Lie_A^!$.
Note the equivalence $Lie_A\simeq \FST\simeq \GFST_A\simeq Lie^{!}_A$
(cf. Remark~\ref{smoothremark}). Moreover, it induces $Lie_A^f\simeq Lie_A^{!,f}$.
Suppose that $P\in Lie_A$ maps to $L$.
Taking account of the induced equivalence
$(Lie^{f}_A)_{/P}\simeq (Lie_A^{!,f})_{/L}$, we are reduced to showing that
$P$ is a colimit of $(Lie^{f}_A)_{/P}\to Lie_A$. This is clear
since $Lie_A\simeq \mathcal{P}_{\Sigma}^{st}(Lie^{f}_A)$.
\QED

\subsection{}
\label{zigzagconstructionsec}

In Section 7, in what follows,
{\it $A$ is a commutative (ordinary) algebra which is smooth over the base field $k$}.
If $\overline{R}$ is the image of $R\in \EXT_A$ under the forgetful functor $\EXT_A\to \CAlg_A$,
then $\Map_{\CAlg_A}(\overline{R},A)$ is a contractible space. Thus, the forgetful functor
$\EXT_A\to \CAlg_A$ is a fully faithful functor into the essential image. 
We often identify $\EXT_A$ with the essential image in $\CAlg_A$.

In Section~\ref{zigzagconstructionsec}, we establish
Corollary~\ref{quickconstruction}, which will be used in Construction~\ref{equivariantobject}.
By abuse of notation, we write 
$\QC:\CAlg_A^+\to \wCat$
for the functor $\CAlg_A^+\stackrel{\textup{forget}}\to \CAlg_A\to \wCat$
where the second functor corresponds to the coCartesian fibration
$\Mod(\Mod_A)\to \CAlg_A$.
Namely, $\QC$ carries $C\in \CAlg_A^+$ to $\Mod_C$,
and we write $\QC(C)$ for $\Mod_C$.
Let $\QC_!|_{\EXT_A}:\EXT_A\to \wCat$
be the functor given on objects by $C\mapsto \QC_!(\Spec C)=\Ind(\Coh(C))$.
For $f:\Spec C'\to \Spec C$ in $\TSZ_A$, it carries $f$ to
the $!$-pullback functor $f^!:\QC_!(\Spec C)\to \QC_!(\Spec C')$
which is the right adjoint to the proper pushforward functor 
 $f^{\textup{IndCoh}}_*:\QC_!(\Spec C')\to \QC_!(\Spec C)$.
This functor is the restriction of the functor $\QC_!:\FST_!\to \wCat$
constructed in \cite{Ind}, \cite{Gai2} (see also Section~\ref{Indreviewsec}).
Let 
$\LMod\circ \DD_1:\Alg_1^+(\Mod_A)\stackrel{\DD_1}{\to} \Alg_1(\Mod_A)^{op} \stackrel{\LMod}{\to} \wCat$
denote the composite
where the first fuctor is given by the $\eone$-Koszul duality functor $\DD_1:\Alg_1^+(\Mod_A)\to \Alg_1^+(\Mod_A)^{op}$ followed by the forgetful functor $\Alg_1^+(\Mod_A)\to \Alg_1(\Mod_A)$ (here we slightly abuse notation),
and the second functor $\LMod$ indicates the functor corresponding to
the Cartesian fibration
$\LMod(\Mod_A)\to \Alg_1(\Mod_A)$.
We consider three functors
 $\EXT_A\to \wCat$
given by the restrictions
$\QC|_{\EXT_A}$, $\QC_!|_{\EXT_A}$, and $\LMod\circ \DD_1 |_{\EXT_A}$
induced by the forgetful functor $\EXT_A\to \Alg_1^+(\Mod_A)$.
Since
$\Rep:(Lie_A)^{op}\to \wCat$ is the composite of $\LMod:\Alg_1(\Mod_A)^{op}\to \wCat$ and the universal enveloping algebra functor $U_1:Lie_A\to \Alg_1(\Mod_A)$,
and  there exists $U_1\circ \DD_\infty\simeq \DD_1$ between functors
$\EXT_A\to \Alg_1^+(\Mod_A)$ (see \cite[Proposition 3.3]{IMA}),
it follows that
$\LMod\circ \DD_1 |_{\EXT_A}$ is equivalent to $\Rep\circ \DD_\infty|_{\EXT_A}:\EXT_A\to \wCat$.

We define $\QC_H':(\GFST_A)^{op}\to \wCat$
to be the right Kan extension of 
$\QC|_{\EXT_A}$
along
$\EXT_A=(\TSZ_A)^{op}\hookrightarrow (\GFST_A)^{op}$.
Similarly, $\QC_!':(\GFST_A)^{op}\to \wCat$ and $\Rep_H':(\GFST_A)^{op}\to \wCat$
are defiend to be right Kan extensions of $\QC_!|_{\EXT_A}$ and $\LMod\circ \DD_1 |_{\EXT_A}$, respectively.

Next, we recall that
for $R \in \EXT_A$, the duality functor
$\QC(R)=\Mod_{R}(\Mod_A) \to \LMod_{\DD_1(R)}(\Mod_A)$ sends $P$ to $P\otimes_{R}A$.
More precisely,
$\DD_1(R)\otimes_A R \to A$
exhibits $A$ as a $\DD_1(R)$-$R$-bimodule.
This integral kernel $A$ determines a functor
$I_R:\Mod_{R}(\Mod_A) \to \LMod_{\DD_1(R)}(\Mod_A) \simeq \Rep(\DD_\infty(R))(\Mod_A)$
given by $P\mapsto P\otimes_{R}A$.
By \cite[2.3.6]{H},
this functor is fully faithful.
By the construction in \cite[Remark 5.9]{IMA},
$I_R$ is functorial in $R\in \EXT_A$. That is,
there is a natural transformation
\[
\mathcal{I}:\QC|_{\EXT_A} \to \LMod \circ \DD_1|_{\EXT_A}
\]
between functors $\EXT_A \to \wCat$, such that the evaluation at each $R\in \EXT_A$ is equivalent to $I_R$.
We let $\Upsilon|_{\EXT_A}:\QC|_{\EXT_A}\to \QC_!|_{\EXT_A}$ denote the natural transformation
between functors $\EXT_A\to \wCat$, which is
induced by $\Upsilon$.

\begin{Definition}
\label{initialtrans}
We define 
\[
\QC_!' \stackrel{\Upsilon'}{\longleftarrow} \QC_H' \stackrel{\mathcal{I}'}{\longrightarrow} \Rep_H'
\]
to be the diagram obtained from 
$\QC_!|_{\EXT_A} \stackrel{\Upsilon|_{\EXT_A}}{\longleftarrow}\QC|_{\EXT_A} \stackrel{\mathcal{I}}{\longrightarrow} \LMod \circ \DD_1|_{\EXT_A}$
by taking the right Kan extension
along $(\TSZ_A)^{op}\hookrightarrow (\GFST_A)^{op}$.
\end{Definition}

\begin{Remark}
\label{fullyremark}
Given $W\in \GFST_A$, 
$\QC_!' \leftarrow \QC_H' \to \Rep_H'$
determines $\QC_!'(W)\leftarrow \QC_H'(W) \to \Rep_H'(W)$
which is equivalent to the diagram of limits in $\wCat$
\[
\lim_{\Spec C\in (\TSZ_A)_{/W}}\QC_!(\Spec C)\leftarrow \lim_{\Spec C\in (\TSZ_A)_{/W}}\QC(\Spec C)\to\lim_{\Spec C\in (\TSZ_A)_{/W}}\LMod_{\DD_1(C)}.
\]
Since $\Upsilon_C:\QC(\Spec C)\to \QC_!(\Spec C)$ and $\QC(\Spec C)\to \LMod_{\DD_1(C)}$
are fully faithful for any $C\in \EXT_A$ (cf. \cite[Section 9]{Ind}, \cite[2.3.6]{H}),
both $\QC_H'(W)\rightarrow \QC_!'(W)$
and $\QC_H'(W) \to \Rep_H'(W)$ are fully faithful.

The restriction $\QC_H|_{\FST_A}:(\FST_A)^{op}\to \wCat$ is equivalent to
the right Kan extension of $\QC_H|_{\EXT_A}$ so that  $\QC_H|_{\FST_A}:(\FST_A)^{op}\to \wCat$ is naturally equivalent to
the composite $\QC_H' \circ (\Theta_A)^{op}:(\FST_A)^{op} \simeq (\GFST_A)^{op} \to \wCat$.
Similarly, the restriction $\Rep_H|_{\FST_A}:(\FST_A)^{op}\to \wCat$ is naturally equivalent to
the composite $\Rep_H' \circ (\Theta_A)^{op}:(\FST_A)^{op}\simeq (\GFST_A)^{op}\to \wCat$. The fully faithful functor $\QC_H'(W)\to \Rep_H'(W)$ can naturally be identified with
$\gamma_W:\QC_H(W)\to \Rep_H(W)$ in Section~\ref{sheaves}.
\end{Remark}

\begin{Definition}
\label{circlecompletion}
Let
$\otimes_AS^1:\CCAlg_A\to \CCAlg_A$ 
be the functor determined by the tensor by $S^1\in \SSS$ in $\CCAlg_A$.
It naturally gives rise to $(\CCAlg_A)^+\to (\CCAlg_A)^+$ which we also denote by $\otimes_AS^1$.
If we regard the constant functor  $(\CCAlg_A)^+\to (\CCAlg_A)^+$
as the functor $\otimes_A\ast$ induced by the cotensor by the contractible space $\ast$,
the obvious map $S^1\to \ast$ determines the natural transformation
$\otimes_AS^1\to \otimes_A\ast$ defined as $\tau:\Delta^1\times (\CCAlg_A)^+\to (\CCAlg_A)^+$.
Let $\textup{comp}:((\CCAlg_A)^+)^{op}\to \FST_A$
be the functor which carries $\Spec C$ to its formal completion $\widehat{\Spec C}$
(see \cite[2.2.5]{H}, Section~\ref{formalstack}),
where we think  of $((\CCAlg_A)^+)^{op}$ as the $\infty$-category
of pointed affine derived schemes over $A$.
Let $\Theta_A:\FST_A\to \GFST_A$ be the functor defined in Construction~\ref{transformalstack},
which is an equivalence since $A$ is smooth (cf. Remark~\ref{smoothremark}).
We define $\xi$ to be the composite 
\[
\Delta^1\times (\CCAlg_A)^+\stackrel{\tau}{\longrightarrow} (\CCAlg_A)^+\stackrel{\textup{comp}^{op}}{\longrightarrow} (\FST_A)^{op}\stackrel{\Theta_A^{op}}{\longrightarrow}(\GFST_A)^{op}.
\]

\end{Definition}

\begin{Remark}
\label{circlecompletionrem}
We write $(\Spec C)^\wedge$ for the image of $\widehat{\Spec C}$
under the equivalence $\Theta_A:\FST_A\stackrel{\sim}{\to} \GFST_A$.
The functor $\xi|_{\{0\}\times(\CCAlg_A)^+}$ carries $C$ to $(\Spec C\otimes_AS^1)^\wedge$.
The functor $\xi|_{\{1\}\times(\CCAlg_A)^+}$ carries $C$ to $(\Spec C)^\wedge$.
Note that $\otimes_AS^1:(\CCAlg_A)^+\to (\CCAlg_A)^+$
is extended to $(\CCAlg_A)^+\to \Fun(BS^1,(\CCAlg_A)^+)$
in the obvious way.
The functor $\otimes_A\ast:(\CCAlg_A)^+\to (\CCAlg_A)^+$
is extended to $(\CCAlg_A)^+\to \Fun(BS^1,(\CCAlg_A)^+)$
given by trivial $S^1$-actions.
Consequently, $\xi$ is extended to $\Delta^1\times (\CCAlg_A)^+\to \Fun(BS^1,(\GFST_A)^{op})$
which we also denote by $\xi$.
\end{Remark}

\begin{Construction}
\label{secondtrans}
Using Definition~\ref{initialtrans} and Definition~\ref{circlecompletion}
we will construct a zig-zag diagram in $\Fun(BS^1,\wCat)$.
Applying $\Fun(BS^1,-)$ to the diagram in Defintion~\ref{initialtrans}
we have functors
\[
\Fun(BS^1,\QC_H'), \Fun(BS^1, \QC_!'), \Fun(BS^1, \Rep_H'):\Fun(BS^1,(\GFST_A)^{op})\to \Fun(BS^1,\wCat)
\]
and the diagram
\[
\Fun(BS^1,\QC_!') \stackrel{\Upsilon'}{\longleftarrow} \Fun(BS^1,\QC_H') \stackrel{\mathcal{I}'}{\longrightarrow} \Fun(BS^1,\Rep_H')
\]
in $\Fun(\Fun(BS^1,(\GFST_A)^{op}),\Fun(BS^1,\wCat))$.
Let $X$ be a pointed formal stack over $A$, that is, an object of $\FST_A$.
Set $(\TSZ_A)_{/X}=\TSZ_A\times_{\FST_A}(\FST_A)_{/X} $.
Consider the composite $\gamma$
\[
\Delta^1\times \bigl((\TSZ_A)_{/X}\bigr)^{op}\stackrel{\textup{forget}}{\longrightarrow} \Delta^1\times (\CCAlg_A)^+ \stackrel{\xi}{\longrightarrow} \Fun(BS^1,(\GFST_A)^{op})\to \Fun(BS^1,\wCat)
\]
where the second functor $\xi$ is defined in Definition~\ref{circlecompletion} and
Remark~\ref{circlecompletionrem}, and the third functor is either 
$\Fun(BS^1,\QC_H'), \Fun(BS^1, \QC_!')$ or  $\Fun(BS^1, \Rep_H')$.
Consider $\bigl((\TSZ_A)_{/X}\bigr)^{op}\to \Fun(\Delta^1,\Fun(BS^1,\wCat))$
which are adjoint to $\gamma$.
We define
\[
\xymatrix{
\QC_!^{\circlearrowleft\wedge}(X) \ar[d] & 
\QC_H^{\circlearrowleft\wedge}(X)  \ar[r] \ar[l] \ar[d] & \Rep_H^{\circlearrowleft\wedge}(X) \ar[d] \\
\QC_!^{\wedge}(X) & 
\QC_H^{\wedge}(X) \ar[r] \ar[l] & \Rep_H^{\wedge}(X)
}
\] 
in $\Fun(BS^1,\wCat)$ to be the diagram obtained by taking the limits of three
functors
$\bigl((\TSZ_A)_{/X}\bigr)^{op}\to \Fun(\Delta^1,\Fun(BS^1,\wCat))$.
The diagram $\QC_!^{\circlearrowleft\wedge}(X)\leftarrow \QC_H^{\circlearrowleft\wedge}(X)\to  \Rep_H^{\circlearrowleft\wedge}(X)$ is informally given by
\[
\xymatrix{
\lim_{\Spec C\in (\TSZ_A)_{/X}}\QC_H'((\Spec C\otimes_AS^1)^\wedge) \ar[r] \ar[d] & \lim_{\Spec C\in (\TSZ_A)_{/X}}\Rep_H'((\Spec C\otimes_AS^1)^\wedge) \\
\lim_{\Spec C\in (\TSZ_A)_{/X}}\QC'_!((\Spec C\otimes_AS^1)^\wedge).
}
\]
The diagram
$\QC_!^{\wedge}(X)\leftarrow \QC_H^{\wedge}(X)\to  \Rep_H^{\wedge}(X)$
also has a similar form (replace $(\Spec C\otimes_AS^1)^\wedge$ by $\Spec C$).
Note that $\QC_H'((\Spec C\otimes_AS^1)^\wedge)  \to \QC'_!((\Spec C\otimes_AS^1)^\wedge)$ and $\QC_H'((\Spec C\otimes_AS^1)^\wedge) \to \Rep_H'((\Spec C\otimes_AS^1)^\wedge)$
are fully faithful for any $C\in \EXT_A$ (cf. Remark~\ref{fullyremark}).
Consequently, both $ \QC_H^{\circlearrowleft\wedge}(X)\rightarrow \QC_!^{\circlearrowleft\wedge}(X)$
and $\QC_H^{\circlearrowleft\wedge}(X) \to \Rep_H^{\circlearrowleft\wedge}(X)$ are fully faithful.
Similarly, $ \QC_H^{\wedge}(X)\rightarrow \QC_!^{\wedge}(X)$
and $\QC_H^{\wedge}(X) \to \Rep_H^{\wedge}(X)$ are fully faithful.
\end{Construction}

Let $(S\times_kS)_{S}^\wedge$ denote $(S\times_kS)\times_{(S\times_kS)_{\DR}}S_{\DR}$,
that is determined by the diagonal $\Delta:S\to S\times_k S$
(see \cite[Vol.II, Chap.4]{Gai2} or the review after Remark~\ref{smoothremark}).
We think of it as
the pointed $!$-formal stack 
obtained from $S\times_kS$ by taking the formal completion along the diagonal $S\to S\times_kS$.
Let $LS=\Spec A\otimes_kS^1=S\times_{S\times_kS}S$ be the free loop space of derived scheme $S=\Spec A$ over $k$.
Let $\iota:S\to LS$ be the morphism induced by constant loops, so that $\iota$ is
an $S^1$-equivariant morphism.
Let $(S\times_kLS)_{S}^\wedge$ denote $(S\times_kLS)\times_{(S\times_kLS)_{\DR}}S_{\DR}$,
determined by $\textup{id}\times \iota:S\to S\times_k LS$.
The equivalence $\Theta_A:\FST_A \stackrel{\sim}{\to} \GFST_A$ carries $\widehat{S\times_kS}$ and $\widehat{S\times_kLS}$
to $(S\times_kS)^\wedge_{S}$ and $(S\times_kLS)^\wedge_{S}$, respectively (see Proposition~\ref{formalcompatibility}). Thus, we will regard $(S\times_kS)^{\wedge}_S$ and $(S\times_k LS)^{\wedge}_S$
as the images of $\widehat{S\times_kS}$ and $\widehat{S\times_kLS}$, respectively.

\begin{Proposition}
\label{descent1}
There exists a canonical equivalence
$\QC_!((S\times_kLS)_S^\wedge)\stackrel{\sim}{\to} \QC_!^{\circlearrowleft\wedge}(\widehat{S\times_kS})$
and $\QC_!((S\times_kS)_S^\wedge)\stackrel{\sim}{\to} \QC_!^{\wedge}(\widehat{S\times_kS})$
in $\Fun(BS^1,\wCat)$.
Here the $S^1$-action on $\QC_!((S\times_kLS)_S^\wedge)$ is determined by that on $LS$,
and  the $S^1$-action on $\QC_!((S\times_kS)_S^\wedge)$ is trivial.
\end{Proposition}

\begin{Proposition}
\label{descent2}
There exists canonical equivalences
$\Rep(\mathbb{T}_{A/k}[-1]^{S^1})(\Mod_A)\simeq \Rep_H^{\circlearrowleft\wedge}(\widehat{S\times_kS})$
and $\Rep(\mathbb{T}_{A/k}[-1])(\Mod_A)\simeq \Rep_H^{\wedge}(\widehat{S\times_kS})$
in $\Fun(BS^1,\wCat)$. Here $S^1$-actions on $\Rep(\mathbb{T}_{A/k}[-1])(\Mod_A)$ and $\Rep_H^{\wedge}(\widehat{S\times_kS})$ are trivial.
\end{Proposition}

Before the proofs of Propositions, we apply
Construction~\ref{secondtrans}, Proposition~\ref{descent1}, and Proposition~\ref{descent2}
to $X=\widehat{S\times_kS}$ and obtain:

\begin{Corollary}
\label{quickconstruction}
There exists a canonical diagram
\[
\xymatrix{
\QC_!((S\times_kLS)_S^\wedge)
 \ar[d] & \QC_H^{\circlearrowleft\wedge}(\widehat{S\times_kS}) \ar[r] \ar[l] \ar[d] & \Rep(\mathbb{T}_{A/k}[-1]^{S^1})(\Mod_A) \ar[d] \\ 
 \QC_!((S\times_kS)_S^\wedge)
 & \QC_H^{\wedge}(\widehat{S\times_kS}) \ar[r] \ar[l] & \Rep(\mathbb{T}_{A/k}[-1])(\Mod_A) 
}
\]
in $\Fun(BS^1, \wCat)$, where the vertical functor on the right side is
determined by the restriction along the diagonal morphism
$\mathbb{T}_{A/k}[-1]\to\mathbb{T}_{A/k}[-1]^{S^1}$.
The vertical functor on the left side is the $!$-pullback functor
along $(\textup{id}_S\times\iota)_S^\wedge:(S\times_kS)_S^\wedge\to (S\times_kLS)_S^\wedge$.
Moreover, horizontal functors are fully faithful (see Construction~\ref{secondtrans}).
If we take $S^1$-invariants, we obtain
\[
\xymatrix{
\QC_!((S\times_kLS)_S^\wedge)^{S^1}
 \ar[d] & \QC_H^{\circlearrowleft\wedge}(\widehat{S\times_kS})^{S^1} \ar[r] \ar[l] \ar[d] & \Rep(\mathbb{T}_{A/k}[-1]^{S^1})(\Mod_A^{S^1}) \ar[d] \\ 
 \QC_!((S\times_kS)_S^\wedge)^{S^1}
 & \QC_H^{\wedge}(\widehat{S\times_kS})^{S^1} \ar[r] \ar[l] & \Rep(\mathbb{T}_{A/k}[-1])(\Mod_A^{S^1}) 
}
\]
where horizontal functors are fully faithful functors.
\end{Corollary}

For the proofs of above Propositions, we need some Lemmata.

\begin{Lemma}
\label{prepare1}
\begin{enumerate}
\item Consider
\[
\beta:\QC_!((S\times_kS)^\wedge_S)\to  \lim_{\Spec C\in (\TSZ_A)_{/\widehat{S\times_kS}}} \QC_!(\Spec C)
\]
determined by $\QC_!$ with the $!$-pullback functoriality on $(\TSZ_A)_{/\widehat{S\times_kS}}\to(\FST_A)_{/\widehat{S\times_kS}}\stackrel{\sim}{\to} (\GFST_A)_{/(S\times_kS)^\wedge_S}$.
Then $\beta$ is an equivalence.

\item Consider
\[
\beta_L:\QC_!((S\times_kLS)^\wedge_S)\to  \lim_{(\TSZ_A)_{/\widehat{S\times_kS}}} \QC_!((\Spec C\otimes_AS^1)^\wedge)
\]
determined by $\QC_!$ with the $!$-pullback functoriality on $(\TSZ_A)_{/\widehat{S\times_kS}}\to(\FST_A)_{/\widehat{S\times_kLS}}\simeq  (\GFST_A)_{/(S\times_kLS)^\wedge_{S}}$,
where $(\TSZ_A)_{/\widehat{S\times_kS}}\to(\FST_A)_{/\widehat{S\times_kLS}}$
is induced by the cotensor by $S^1$ in $\FST_A$, that is, the functor informally
given by $[\Spec C\to \widehat{S\times_kS}]\mapsto [\widehat{\Spec C\otimes_AS^1}\to \widehat{S\times_k LS}]$.
Then $\beta_L$ is an equivalence.
Moreover, $\beta_L$ is an $S^1$-equivariant equivalence, that is,
an equivalence in $\Fun(BS^1,\wCat)$.
\end{enumerate}
\end{Lemma}

\Proof
We prove (1). Write $U_A$ for $(\TSZ_A)_{/\widehat{S\times_kS}}$.
Taking into account the compatibility of $\QC_!$ with sifted colimits in $\GFST_A$ (see \cite[Vol. II, Chap. 7, 5.3.2]{Gai2}) and the fact of $U_A$ is sifted (see Claim~\ref{sifteddiagram}),
we are reduced to showing that $(S\times S)^\wedge_S$ is a colimit of the functor $s:U_A\to\TSZ_A\to \GFST_A$
given by $U_A\ni \Spec C \mapsto \Spec C$ where the latter $\Spec C$
indicates the object of $\GFST_A$, corepresented by $C$.
Observe that $\widehat{S\times_k S}$ is a colimit of
the functor $u:U_A\to \TSZ_A\hookrightarrow  \FST_A$ given by the
forgetful functor.
Using the equivalences $Lie_A\simeq \FST_A$ and $(Lie_A^f)_{/\mathbb{T}_{A/k}[-1]}\simeq (\TSZ_A)_{/\widehat{S\times S}}$,
$u$ corresponds to the forgetful functor
$(Lie_A^f)_{/\mathbb{T}_{A/k}[-1]}\to Lie_A$ whose colimit is $\mathbb{T}_{A/k}[-1]$.
It follows that
$\widehat{S\times_k S}$ is a colimit of
$u$. Since $\FST_A\simeq \GFST_A$ carries $\widehat{S\times_kS}$ to
$(S\times S)^\wedge_S$, $(S\times S)^\wedge_S$ is a colimit of $s$.

We prove (2).
If $\beta_L$ is a categorical equivalence, the final assertion is clear
since the diagram of pointed formal stacks is promoted to an $S^1$-equivariant diagram
in the obvious way. We will prove that $\beta_L$ is an equivalence in $\wCat$.
Taking into account the compatibility of $\QC_!$ with sifted colimits in $\GFST_A$ again,
it suffices to show that $(S\times LS)^\wedge_S$ is a colimit of the functor $t:U_A\to\TSZ_A\to \GFST_A$
given by $U_A\ni \Spec C \mapsto (\Spec C\otimes_AS^1)^\wedge$.
Let $v:U_A\to  \FST_A$ be the functor 
given by $\Spec C\mapsto \widehat{\Spec C\otimes_AS^1}$.
Notice that $\widehat{\Spec C\otimes_AS^1}$ is naturally equivalent to
the fiber product $\Spec C\times_{\Spec C\times \Spec C}\Spec C$
in $\FST_A$.
Therefore, through the equivalences $Lie_A\simeq \FST_A$ and $(Lie_A^f)_{/\mathbb{T}_{A/k}[-1]}\simeq (\TSZ_A)_{/\widehat{S\times S}}$,
 $v:U_A\to  \FST_A$
corresponds to 
$(Lie_A^f)_{/\mathbb{T}_{A/k}[-1]}\to Lie_A$ given by
$[L\to \mathbb{T}_{A/k}[-1]]\mapsto L^{S^1}=L\times_{L\times L}L$
whose colimit is $\mathbb{T}_{A/k}[-1]^{S^1}$
since the functor $(-)^{S^1}:Lie_A\to Lie_A$ given by cotensor by $S^1$
preserves sifted colimits (since finite limits commute with small
colimits in $\Mod_A$, and $Lie_A\to \Mod_A$ preserves sifted colimits).
It follows from the correspondence between
$\mathbb{T}_{A/k}[-1]^{S^1}$ and $\widehat{S\times_k LS}$ via $Lie_A\simeq \FST_A$
that $\widehat{S\times_k LS}$ is a colimit of $v$.
The equivalence $\FST_A\stackrel{\sim}{\to} \GFST_A$.
sends $\widehat{S\times_k LS}$
to $(S\times_kLS)^\wedge_S$.
It follows that $(S\times_kLS)^\wedge_S$ is a colimits of $t$, as desired.
\QED

{\it Proof of Propostion~\ref{descent1}.}
We first prove $\QC_!((S\times_kLS)_S^\wedge)\stackrel{\sim}{\to} \QC_!^{\circlearrowleft\wedge}(\widehat{S\times_kS})$. 
By the equivalence $\beta_L$ in Lemma~\ref{prepare1} (2)
it will suffice to construct an equivalence 
$\lim_{(\TSZ_A)_{/\widehat{S\times_kS}}} \QC_!((\Spec C\otimes_AS^1)^\wedge)\simeq 
\lim_{(\TSZ_A)_{/\widehat{S\times_kS}}} \QC'_!((\Spec C\otimes_AS^1)^\wedge)$.
Recall that $\QC_!':(\GFST_A)^{op} \to \wCat$ is 
a right Kan extension of the restriction of $\QC_!|_{\GFST_A}:(\GFST_A)^{op}\to \wCat$ to $\EXT_A$. 
Thus, there exists a natural transformation $\QC_!\to \QC_!'$.
It is enough to prove that for any $C\in \EXT_A$ the induced functor
$\QC_!((\Spec C\otimes_AS^1)^\wedge)\to  \QC_!'((\Spec C\otimes_AS^1)^\wedge)\simeq \lim_{\Spec R\in (\TSZ_A)_{/(\Spec C\otimes_AS^1)^\wedge}}\QC_!(R)$
is an equivalence.
If $L$ is a dg Lie algebra corresponding to $(\Spec C\otimes_AS^1)^\wedge$,
then there exists $(\TSZ_A)_{/(\Spec C\otimes_AS^1)^\wedge}\simeq (Lie^f_A)_{/L}$.
Since the sifted colimit of $(Lie^f_A)_{/L}\to Lie_A\simeq \GFST_A$
is $(\Spec C\otimes_AS^1)^\wedge$, it follows from the compatibility of $\QC_!$
with sifted colimits in $\GFST_A$ \cite[Vol. II, Chap. 7, 5.3.2]{Gai2}
that  $\QC_!((\Spec C\otimes_AS^1)^\wedge)\to \lim_{\Spec R\in (\TSZ_A)_{/(\Spec C\otimes_AS^1)^\wedge}}\QC_!(R)$ is an equivalence.

The proof of $\QC_!((S\times_kS)_S^\wedge)\stackrel{\sim}{\to} \QC_!^{\wedge}(\widehat{S\times_kS})$ 
is similar and easier.
\QED

\begin{Lemma}
\label{vsrepresentation}
Let $W$ be a pointed formal stack.
We write $L_W$ for the dg Lie algebra $\mathcal{L}(W)$ corresponding to $W$.
Then there exist canonical equivalences
$\lim_{\Spec C\in (\TSZ_A)_{/W}} \Rep(\DD_\infty(C\otimes_AS^1))(\Mod_A)\simeq \Rep(L_W^{S^1})(\Mod_A)$ and $\lim_{\Spec C\in (\TSZ_A)_{/W}} \Rep(\DD_\infty(C))(\Mod_A)\simeq \Rep(L_W)(\Mod_A)$
in $\Fun(BS^1,\wCat)$. Here $S^1$-actions on $\Rep(\DD_\infty(C))(\Mod_A)$ and $\Rep(L_W)(\Mod_A)$ are trivial.
\end{Lemma}

\Proof
We will prove the first equivalence. The proof of the second equivalence is similar. 
Notice first that the equivalence $\FST_A\simeq Lie_A$ induces $ (\TSZ_A)_{/W}\simeq (Lie_A^f)_{/L_W}$.
The functor $ (\TSZ_A)_{/W}\to \TSZ_A\to \Fun(BS^1,Lie_A)$ given by $\Spec C\mapsto \DD_\infty(C\otimes_AS^1)$ can be identified with
$(Lie_A^f)_{/L_W}\to \Fun(BS^1,Lie_A)$ induced by cotensor by $S^1$ (i.e., $L\mapsto L^{S^1}$).
Since the functor $(-)^{S^1}:Lie_A\to Lie_A$ given by cotensor by $S^1$
preserves sifted colimits and the forgetful functor
$\Fun(BS^1,Lie_A)\to Lie_A$ preserves colimits, we see that
$L_W^{S^1}$ is a colimit of $ (\TSZ_A)_{/W}\to \TSZ_A\to \Fun(BS^1,Lie_A)$.
Using the fact that  the formulation of $\LMod$ (with restriction functors) carries sifted colimits in $\Alg_1(\Mod_A)$ to limits \cite[X, 2.4.32]{DAG},
we have $\lim_{\Spec C\in  (\TSZ_A)_{/W}} \Rep(\DD_\infty(C\otimes_AS^1))(\Mod_A)
\simeq\Rep(L_W^{S^1})(\Mod_A)$ in $\Fun(BS^1,\wCat)$, as desired.
\QED

\begin{Lemma}
\label{vsrepresentation2}
There exists a canonical equivalence
$\Rep(\DD_\infty(B))(\Mod_A)\stackrel{\sim}{\to}\Rep_H'((\Spec B)^\wedge)$
which is functorial in $B\in (\CCAlg_A)^+$.
\end{Lemma}

\Proof
Combine Lemma~\ref{replimit} and Remark~\ref{fullyremark}
to deduce $\Rep(\DD_\infty(B))(\Mod_A)\simeq \Rep_H(\widehat{\Spec B})\simeq  \Rep_H'((\Spec B)^\wedge)$.
\QED

{\it Proof of Proposition~\ref{descent2}.}
According to Lemma~\ref{vsrepresentation} and Lemma~\ref{vsrepresentation2},
we have canonical equivalences
\begin{eqnarray*}
\Rep_H^{\circlearrowleft\wedge}(\widehat{S\times_kS}) &\simeq& \lim_{\Spec C\in (\TSZ_A)_{/\widehat{S\times_kS}}} \Rep_H'((\Spec C\otimes_AS^1)^\wedge)(\Mod_A)  \\ 
&\simeq& \lim_{\Spec C\in (\TSZ_A)_{/\widehat{S\times_kS}}} \Rep(\DD_\infty(C\otimes_AS^1))(\Mod_A) 
\simeq \Rep(\mathbb{T}_{A/k}[-1]^{S^1})(\Mod_A),
\end{eqnarray*}
as desired.
Replacing $C\otimes_AS^1$ by $C$ we also obtain $\Rep_H^{\wedge}(\widehat{S\times_kS})\simeq \Rep(\mathbb{T}_{A/k}[-1])(\Mod_A)$.
\QED

\begin{Remark}
It is possible to construct categorical equivalences
$\Rep(\mathbb{T}_{A/k}[-1]^{S^1})(\Mod_A)\simeq \QC_!((S\times_kLS)^\wedge_S)$
and $\Rep(\mathbb{T}_{A/k}[-1])(\Mod_A)\simeq \QC_!((S\times_kS)^\wedge_S)$.
Indeed, in the earlier version we constructed such equivalences.
We prefer to the zig zag diagram in Corollary~\ref{quickconstruction}
because this diagram is clearly related to $\varUpsilon$ and
the duality functor $\Mod_C\to \LMod_{\DD_1(C)}$.
\end{Remark}

The following Lemma will be proved in \cite[Corollary 6.8]{DC2} (see also Remark~\ref{mild}).

\begin{Lemma}
\label{postponed}
The object $\mathcal{V}_L' \in \Rep(\mathbb{T}_{A/k}[-1]^{S^1})(\Mod_A^{S^1})$
in Proposition~\ref{naiveequivariant} lies in the essential image of the fully faithful functor
$\QC^{\circlearrowleft\wedge}_H(\widehat{S\times_kS})^{S^1}\to \Rep(\mathbb{T}_{A/k}[-1]^{S^1})(\Mod_A^{S^1})$.
\end{Lemma}

\begin{Remark}
\label{mild}
There is a quick proof of Lemma~\ref{postponed}, under a mild condition, which is different from the argument in \cite{DC2}.
Here is the sketch of the proof. We will assume that $J_{\HH_\bullet(\CCC/A)}:\Def(\HH_\bullet(\CCC/A))\to \FF_{\End^{L}(\HH_\bullet(\CCC/A))}$ is an equivalence. 
For example, if $\HH_\bullet(\CCC/A)$ is right-bounded, i.e., $H^n(\HH_\bullet(\CCC/A))=0$ for $n>>0$,
then the proof of \cite[X, 5.2.14]{DAG} shows that $J_{\HH_\bullet(\CCC/A)}$ is an equivalence.
By the equivalence $J_{\HH_\bullet(\CCC/A)}:\Def(\HH_\bullet(\CCC/A)) \stackrel{\sim}{\to} \FF_{\End^{L}(\HH_\bullet(\CCC/A))}$, $\widehat{S\times_kLS} \to \FF_{\End^{L}(\HH_\bullet(\CCC/A))}$
can be promoted to $\widehat{S\times_kLS} \to \Def(\HH_\bullet(\CCC/A))$ (see the diagram
in the proof of Proposition~\ref{looppromotion}). Moreover, $\widehat{S\times_kLS} \to \Def(\HH_\bullet(\CCC/A))$
is promoted to an $S^1$-equivariant morphism.
Consider the universal deformation $\mathcal{U}\in \QC_H(\Def(\HH_\bullet(\CCC/A)))\times_{\Mod_A}\{\HH_\bullet(\CCC/A)\}$.
By Proposition~\ref{looppromotion} and its proof, $\mathcal{V}_L\in \Rep_H(\widehat{S\times_kLS})=\Rep(\mathbb{T}_{A/k}[-1]^{S^1})(\Mod_A)$ is the image of $\mathcal{U}$ under 
$\QC_H(\Def(\HH_\bullet(\CCC/A)))\to \QC_H(\widehat{S\times_kLS})\subset \Rep_H(\widehat{S\times_kLS})$.
Taking into account Proposition~\ref{naiveequivariant}, we see that $\mathcal{V}_L'$ belongs to $\QC_H(\widehat{S\times_kLS})^{S^1}\subset \Rep_H(\widehat{S\times_kLS})^{S^1}$.
Lemma~\ref{postponed} follows by observing that the fully faithful functor
$\QC_H(\widehat{S\times_kLS})^{S^1}\hookrightarrow \Rep_H(\widehat{S\times_kLS})^{S^1}\simeq \Rep(\mathbb{T}_{A/k}[-1]^{S^1})(\Mod_A^{S^1})$
factors through $\QC^{\circlearrowleft\wedge}_H(\widehat{S\times_kS})^{S^1}\hookrightarrow \Rep(\mathbb{T}_{A/k}[-1]^{S^1})(\Mod_A^{S^1})$.
\end{Remark}


\subsection{}
\label{indfcompletion}
The goal of Section~\ref{indfcompletion} is Theorem~\ref{cyclicextension}.
Let $LS^\wedge_S$ be the formal completion $LS\times_{(LS)_{\DR}}S_{\DR}$
along the morphism $\iota:S\to LS$.
Here we consider $(S\times_kS)_{S}^\wedge$, $(S\times_kLS)_{S}^\wedge$ and $LS^\wedge_S$
as functors $\CCAlgftec_k\to \SSS$ (i.e., prestacks locally almost of finite type).
Observe that $S_{\DR}\to (LS)_{\DR}$ is an equivalence so that the
canonical morphism $LS^\wedge_S\to LS $ is an equivalence.
Indeed, by definition, for $R\in \CCAlg_k$, $\Map_{\PST}(\Spec R,(LS)_{\DR})$ is canonically equivalent to
$\Map_{\PST}(\Spec R_{red}, \Map(S^1,S))\simeq \Map_{\SSS}(S^1,S(R_{red}))\simeq S(R_{red})$.
The final equivalence follows from the facts that $S^1$ is connected and $S(R_{red})$
is a discrete space because $S$ is a derived scheme (or algebraic space).
Thus, $S_{\DR}\simeq (LS)_{\DR}$.
The derived (affine) schemes $S\times_kS$, $S\times_{k} LS$ and $LS$
are almost of finite type over $k$.
(By 
\cite[7.2.4.31]{HA}, the condition of almost of finite type over $k$ is
equivalent to the condition of almost of finite presentation over $k$ in \cite[Definition 7.2.4.26]{HA}.)
Indeed, since we assume that $A$ is almost of finite type over $k$, 
it follows from \cite[Vol. I, Chap.2, 1.6.6, 1.7.10]{Gai2}
that the finite limits $LS=S\times_{S\times_kS}S$, $S\times_kS$ and $S\times_k(S\times_{S\times_kS}S)$
are also almost of finite type over $k$.

Consider  the square diagram in $\Fun(\CCAlgftec_k,\SSS)$:
\[
\xymatrix{
(S\times_kS)_{S}^\wedge \ar[r]^{(\textup{id}\times \iota)^\wedge_S} \ar[d] & (S\times_kLS)^\wedge_{S} \ar[d] \\
S \ar[r]_{\iota} & LS_S^{\wedge}
}
\]
where the vertical morphisms are second projections.
The functor $\QC_!$  gives rise to
\[
G:\QC_!(LS)\simeq \QC_!(LS_S^\wedge)\to \QC_!(S)\times_{\QC_!((S\times_kS)_{S}^\wedge)}\QC_!((S\times_kLS)^\wedge_{S})
\]

\begin{Proposition}
\label{gluingloop}
The functor $G$ is and equivalence of $\infty$-categories.
Moreover, this functor is promoted to an $S^1$-equivariant functor.
\end{Proposition}

\Proof
Consider the trivial action of $S^1$ on $(S\times_kS)_{S}^\wedge$
and $S$. The actions of $S^1$ on $LS \simeq LS_S^\wedge$ and $ (S\times_kLS)^\wedge_{S}$
are induced by the canonical $S^1$-action.
The above square is promoted to an $S^1$-equivariant diagram.
Thus, $G$ is also promoted to an $S^1$-equivariant diagram.
It will suffice to show that the underlying functor is an equivalence.
For this purpose, we first note that the square diagram can be regarded
as a diagram in $\GFST_A$. The morphisms from $S$ is given
by canonical maps $\textup{id}:S\to S$, $\Delta_S:S\to (S\times_kS)_{S}^\wedge$, $\textup{id}\times \iota:S\to (S\times_k LS)_{S}^\wedge$ and $\iota:S\to LS\simeq LS_S^{\wedge}$, respectively.
The morphisms to $S$ is given by
$\textup{id}:S\to S$, the second projection
$(S\times_kS)_{S}^\wedge\to S$, $(S\times_k LS)_{S}^\wedge\to LS\to S$ and $LS\simeq LS_S^{\wedge}\to S$, respectively.
Then the square diagram in $\GFST_A$ induces the square diagram in $Lie_A^!$
\[
\xymatrix{
L_{(S\times_kS)_{S}^\wedge} \ar[r] \ar[d] & L_{(S\times_kLS)^\wedge_{S}} \ar[d] \\
0 \ar[r] & L_{LS_S^{\wedge}}
}
\]
For ease of notation, we put $L_1=L_{(S\times_kS)_{S}^\wedge}$, $L_2= L_{(S\times_kLS)^\wedge_{S}}$,
$L_3=L_{LS_S^{\wedge}}$.
Thanks to \cite[Vol.II, Chap.7, 5.1.2, 5.2]{Gai2}, $G$ can be identified with the morphism induced by the restriction functors
\[
G_{\Rep}:\Rep(L_3)(\QC_!(A))\to \QC_!(A) \times_{\Rep(L_1)(\QC_!(A))}\Rep(L_2)(\QC_!(A)).
\]
According to \cite[Lemma 2.4.32, 2.4.33]{DAG} (one can apply the argument to our situation
by replacing $\Mod_k$ with $\QC_!(A)$ in {\it loc. cit.}), if the morphism from the pushout $0\sqcup_{L_1}L_2\to L_3$
in $Lie_A^!$ is an equivalence (or $A\sqcup_{U_1(L_1)}U_1(L_2)\simeq U_1(L_3)$ in $\Alg_1(\Mod_A)$), then $G_{\Rep}$ is an equivalence.
Since $((S\times_kS)_{S}^\wedge)\times_S(LS_S^{\wedge})\simeq (S\times_kLS)^\wedge_{S}$ in $\GFST_A$, the sequence $L_1\to L_2 \to L_3$ is $L_1\times\{0\}\hookrightarrow L_1\times L_3\stackrel{\textup{pr}_2}{\to} L_3$. Thus, $0\sqcup_{L_1}L_2=L_2/L_1\simeq L_3$.
In other words, passing to universal enveloping algebras
$U_1(0)\sqcup_{U_1(L_1)}U_1(L_1\times L_3)\simeq A\sqcup_{U_1(L_1)}(U_1(L_1)\otimes_A U_1(L_3)) \simeq A\sqcup_{A\oplus \overline{U_1(L_1)}}((A\oplus \overline{U_1(L_1)})\otimes_A U_1(L_3))\simeq U_1(L_3)$.
 Therefore, our assertion follows.
\QED

\begin{Construction}
\label{equivariantobject}

We will define an object of $\QC_!(LS)^{S^1}$, which will be denoted by $\mathcal{H}_\circlearrowleft(\CCC)$. 
We consider the diagram
\[
\xymatrix{
\QC(S)^{S^1}\ar[r]^(0.4){\textup{pr}_2^*} \ar[d]_{\varUpsilon_A} & \QC_H^\wedge(\widehat{S\times_kS})^{S^1}  \ar[d] & \QC_H^{\circlearrowleft\wedge}(\widehat{S\times_kS})^{S^1} \ar[l] \ar[d]  \\
\QC_!(S)^{S^1}\ar[r]_(0.4){\textup{pr}_2^!} &  \QC_!((S\times_kS)_{S}^\wedge)^{S^1} & \QC_!((S\times_kLS)^\wedge_{S})^{S^1} \ar[l]
}
\]
which commutes up to canonical homotopy, where 
the right square comes from Corollary~\ref{quickconstruction}.
The functor $\textup{pr}_2^*$ is induced by
$\QC(S)\to \lim_{\Spec C\in (\TSZ_A)_{/\widehat{S\times S}}}\QC(C)\simeq \QC_H^\wedge(\widehat{S\times_kS})$
determined by the second projection $\textup{pr}_2:S\times_kS\to S$
and $*$-pullbacks.
Likewise, $\textup{pr}_2^!$ is
$\QC_!(S)\to \lim_{\Spec C\in (\TSZ_A)_{/\widehat{S\times S}}}\QC_!(\Spec C)\simeq  \QC_!((S\times_kS)_{S}^\wedge)$ induced by th second projection and $!$-pullbacks.
This diagram induces the fully faithful functor
\[
\QC(S)^{S^1}\times_{\QC_H^\wedge(\widehat{S\times_kS})^{S^1} }  \QC_H^{\circlearrowleft\wedge}(\widehat{S\times_kS})^{S^1}\hookrightarrow \QC_!(S)^{S^1}\times_{\QC_!((S\times_kS)_{S}^\wedge)^{S^1}} \QC_!((S\times_kLS)^\wedge_{S})^{S^1}\simeq \QC_!(LS)^{S^1}.
\]

By Proposition~\ref{naiveequivariant} and Lemma~\ref{postponed},
we have constructed
\[
\mathcal{V}_\dagger \in \QC(S)^{S^1}\times_{\QC_H^\wedge(\widehat{S\times_kS})^{S^1} }  \QC_H^{\circlearrowleft\wedge}(\widehat{S\times_kS})^{S^1}
\]
(see also Remark~\ref{fullyremark}).
We let $\mathcal{H}_\circlearrowleft(\CCC)$ denote the image of $\mathcal{V}_\dagger$ in $\QC_!(LS)^{S^1}$. 
\end{Construction}

We deduce the following result from the above Construction.

\begin{Theorem}
\label{cyclicextension}
Let $\CCC$ be an $A$-linear small stable $\infty$-category.
We have constructed an object $\mathcal{H}_\circlearrowleft(\CCC)$ of $\QC_!(LS)^{S^1}$
having the following properties:
\begin{enumerate}
\renewcommand{\labelenumi}{(\roman{enumi})}
\item
the image in $\QC_!(A)^{S^1}$ obtained by the $!$-pullback along $\iota:S\to LS$
is naturally equivalent to $\HH_{\bullet}(\CCC/A)\otimes_A\omega_A$
(so that it maps to $\HH_\bullet(\CCC/A)\in \Mod_A=\QC(S)^{S^1}$ through the equivalence
$\QC(S)^{S^1}\simeq \QC_!(S)^{S^1}$),

\item the image in $\QC_!((S\times_kLS)^\wedge_{S})^{S^1}$
obtained by the $!$-pullback along $(S\times_kLS)^\wedge_{S}\to LS$
is naturally equivalent to the object
corresponding to the canonical $\mathbb{T}_{A/k}[-1]^{S^1}$-module $\HH_\bullet(\CCC/A)$
(see Definition~\ref{canonicaltangentaction}) 
through the zig-zag of fully faithful functors
\[
\QC_!((S\times_kLS)^\wedge_{S})^{S^1}\hookleftarrow  \QC_H^{\circlearrowleft\wedge}(\widehat{S\times_kS})^{S^1} \hookrightarrow \Rep(\mathbb{T}_{A/k}[-1]^{S^1})(\Mod_A^{S^1}).
\]

\item the image in $\QC_!((S\times_k S)^\wedge_{S})^{S^1}$
is naturally equivalent to the object
corresponding to the canonical $\mathbb{T}_{A/k}[-1]$-module $\HH_\bullet(\CCC/A)$
(see Definition~\ref{canonicaltangentaction}) 
through the zig-zag of fully faithful functors
\[
\QC_!((S\times_kS)^\wedge_{S})^{S^1} \hookleftarrow  \QC_H^{\circlearrowleft\wedge}(\widehat{S\times_kS})^{S^1} \hookrightarrow \Rep(\mathbb{T}_{A/k}[-1])(\Mod_A^{S^1}).
\]

\end{enumerate}

\end{Theorem}


\section{D-modules}
\label{dmodule}

Let $S$ be a derived scheme locally almost of finite type over a field $k$
of characteristic zero. Suppose that $S$ is quasi-compact.
Let $\CCC_S$ be a stable $\infty$-category over $S$
(see Example~\ref{geometricex2}).
In this section, we apply Theorem~\ref{preeasy} and Theorem~\ref{cyclicextension}
to prove that the periodic cyclic homology/complex (that computes complex computing periodic cyclic homology) admits a $D$-module structure.

\vspace{2mm}

\subsection{}

We first consider 
\[
 \xymatrix{
S \ar[r]^{\iota} \ar[d]_i &  LS \ar[d] \ar[dl]_{\pi} \\
S_{\DR} \ar[r]^{\sim}_{\iota_{\DR}} & (LS)_{\DR}
 }
\]
where the vertical functors are the canonical functors to de Rham prestacks.
As observed in Section~\ref{indfcompletion}, $\iota_{\DR}$ is an equivalence.
We denote by $\pi$ the composite $LS\to (LS)_{\DR}\stackrel{\iota_{\DR}^{-1}}{\to} S_{\DR}$.
Let $\Coh(LS)\subset \QC_!(LS)$ be the full subcategory of coherent complexes.
From \cite[4.4.3]{P}, $\QC_!(LS)=\Ind(\Coh(LS))$ so that $\Coh(LS)$ coincides with the full subcategory of compact objects. 
We write $\Coh(LS)^{S^1}\otimes_{k[t]}k[t,t^{-1}]$ for the Tate construction $\Coh(LS)^{S^1}\otimes_{\Perf_{k[t]}}\Perf_{k[t,t^{-1}]}$ where $\Perf_{k[t]}$ is the full subcategory of $\Mod_{k[t]}$
spanned by compact objects, and $k[t]$ is the free commutative dg algebra generated by one element $t$ of cohomological degree two
(see e.g. \cite{BN}, \cite{P} for the Tate construction).
Let $\Ind(\Coh(LS)^{S^1})\otimes_{k[t]}k[t,t^{-1}]$ denote $\Ind(\Coh(LS)^{S^1})\otimes_{\Ind(\Perf_{k[t]})}\Ind(\Perf_{k[t,t^{-1}]})=\Ind(\Coh(LS)^{S^1})\otimes_{\Mod_{k[t]}}\Mod_{k[t,t^{-1}]}$.
Let $\Coh(S_{\DR})\otimes_k k[t]$ be $\Coh(S_{\DR})\otimes_k \Perf_{k[t]}$ which can be identified with
the limit of  the trivial $S^1$-action (homotopy fixed points)
on $\Coh(S_{\DR})$ (the full subcategory $\Coh(S_{\DR})$ may and will be identified with the full subcategory of the compactly generated $\infty$-category
$\QC_!(S_{\DR})$, which consists of compact objects, see \cite{Gai2} or \cite[4.4.4]{P}).
By \cite[4.4.4, 4.5.4]{P}, there exists
the standard equivalence $\Coh(S_{\DR})^{S^1}\simeq \Coh(S_{\DR})\otimes_k\Perf_{k[t]}$
that comes from the trivial $S^1$-action.
Similarly, we let $\Coh(S_{\DR})\otimes_k k[t,t^{-1}]$ denote $\Coh(S_{\DR})\otimes_k \Perf_{k[t,t^{-1}]}$.
Now we use the relation between loop spaces and $D$-modules (cf. \cite{BN} and \cite[Theorem 1.3.5]{P}).
We use the result presented in \cite{P}:
the pushforward functor $(\pi_*)^{S^1}:\Coh(LS)^{S^1}\to \Coh(S_{\DR})^{S^1}\simeq \Coh(S_{\DR})\otimes_k\Perf_{k[t]}$
with the base change to $k[t,t^{-1}]$ gives an equivalence
\[
(\pi_*)^{S^1}\otimes_{k[t]}k[t,t^{-1}]:\Coh(LS)^{S^1}\otimes_{k[t]} k[t,t^{-1}] \stackrel{\sim}{\longrightarrow} \Coh(S_{\DR})\otimes_k k[t,t^{-1}]
\]
of $k[t,t^{-1}]$-linear small stable $\infty$-categories. 
Passing to $\Ind$-categories, it gives rise to
\[
\textup{LD}_S:\Ind(\Coh(LS)^{S^1})\otimes_{k[t]}k[t,t^{-1}]\stackrel{\sim}{\longrightarrow} \Ind(\Coh(S_{\DR}))\otimes_k k[t,t^{-1}]=\Ind(\Coh(S_{\DR}))\otimes_k \Mod_{k[t,t^{-1}]}.
\]
We define the $\infty$-category of crystals (right D-modules) on $S$ to be
$\Ind(\Coh(S_{\DR}))$.
Set $\textup{Crys}(S):=\Ind(\Coh(S_{\DR}))$.
We then regard an object of $\textup{Crys}(S)\otimes_k k[t,t^{-1}]$
as a $\ZZ/2\ZZ$-periodic crystal/right D-module.
Therefore, this equivalence brings us a relation
between $\Ind(\Coh(LS)^{S^1})$ and D-modules up to Tate construction.

\begin{Construction}
\label{indequivariant}
We continue to suppose that $\CCC_S$ is a stable $\infty$-category over $S$.
For ease of notation, we write $\CCC=\CCC_S$.
We apply the construction of $\Ind$-categories to the exact functor $\Coh(LS)^{S^1} \hookrightarrow \QC_!(LS)^{S^1}$
to obtain a colimit-preserving functor $\Psi_{LS}:\Ind(\Coh(LS)^{S^1}) \to \QC_!(LS)^{S^1}$.
By adjoint functor theorem \cite{HTT}, there is a right adjoint functor
\[
\Phi_{LS}:\QC_!(LS)^{S^1}\to \Ind(\Coh(LS)^{S^1}).
\]
Consider the diagram
\[
\xymatrix{
\Ind(\Coh(LS)^{S^1})\ar[r] \ar[d] &  \Ind(\Coh(S_{\DR}))\otimes_k \Mod_{k[t]}  \ar[d] \\
\Ind(\Coh(LS)^{S^1})\otimes_{k[t]}k[t,t^{-1}]\ar[r]^{\textup{LD}_S} & \Ind(\Coh(S_{\DR}))\otimes_k \Mod_{k[t,t^{-1}]}
}
\]
Let $\mathcal{H}'_{\circlearrowleft}(\CCC)$ be the image of $\mathcal{H}_\circlearrowleft(\CCC)$ (see
Construction~\ref{equivariantobject}) under $\Phi_{LS}:\QC_!(LS)^{S^1}\to \Ind(\Coh(LS)^{S^1}).$
Let $\mathsf{H}_{\circlearrowleft}(\CCC)$ be the image of $\HH_\bullet(\CCC_X/X)_L\in \QC(LS)^{S^1}$
under $\QC(LS)^{S^1}\stackrel{\Upsilon_{LS}}{\to} \QC_!(LS)^{S^1}\stackrel{\Phi_{LS}}{\to}\Ind(\Coh(LS)^{S^1})$.
\end{Construction}

\begin{Definition}
By Theroem~\ref{preeasy} (see also Example~\ref{geometricex2}),
there is the canonically defined object $\HH_\bullet(\CCC/S)_L\in \QC(LS)^{S^1}$.
Let $\Omega^{\circ, \textup{int}}(\CCC)$
be $\Phi_{LS}(\HH_\bullet(\CCC/S)_L)$ in $\Ind(\Coh(LS)^{S^1})$.
Let $\Omega^\circ(\CCC)$ 
be the images of $\HH_\bullet(\CCC/S)_L$ in $\Ind(\Coh(S_{\DR}))\otimes_k \Mod_{k[t,t^{-1}]}$.

Assume that $S$ is (smooth) affine over $k$.
Let $\Omega^{\bullet,\textup{int}}(\CCC)$
be $\Phi_{LS}(\mathcal{H}_{\circlearrowleft}(\CCC))$ in $\Ind(\Coh(LS)^{S^1})$.
Let $\Omega^\bullet(\CCC)$
be the images of $\mathcal{H}_{\circlearrowleft}'(\CCC)$ in $\Ind(\Coh(S_{\DR}))\otimes_k \Mod_{k[t,t^{-1}]}$.
\end{Definition}

\begin{Remark}
In a subsequent paper \cite{DC2}, we will prove that 
$\mathcal{H}_{\circlearrowleft}(\CCC)\simeq \varUpsilon_{LS}(\HH_\bullet(\CCC/S)_L)$ when $S$ is affine and smooth.
Consequently, in Theorem~\ref{maindmodule} below we have $\Omega^\circ(\CCC)\simeq \Omega^\bullet(\CCC)$.
\end{Remark}

\subsection{}
We briefly overview the definition of the periodic cyclic homology/complex.
In what follows, we assume further that $S$ is a smooth scheme over $k$.
Let $\Perf(S)$ be the full subcategory of $\QC(S)$ spanned by perfect complexes.
Let $\Perf(S)^{S^1}\to \QC(S)^{S^1}$ be the canonical inclusion
and let $\Psi^{\textup{QC}}_S:\Ind(\Perf(S)^{S^1})\to \QC(S)^{S^1}$ be the (essentially unique)
colimit-preserving functor which extends the inclusion.
We define
\[
\Phi^{\textup{QC}}_S:\QC(S)^{S^1} \longrightarrow \Ind(\Perf(S)^{S^1})
\]
to be the right adjoint of $\Psi_S$.
We note that there exists an equivalence 
$\Coh(S)^{S^1}\stackrel{\sim}{\to} \Coh(S)\otimes_k\Perf_{k[t]}$
which carries $M$ to the homotopy fixed points $M^{S^1}$ (cf. \cite[Lemma 4.5.4]{P}). Since $S$ is smooth, $\Perf(S)=\Coh(S)$.
It follows that
\[
\Ind(\Perf(S)^{S^1})\simeq \Ind(\Perf(S))\otimes_k\Ind(\Perf_{k[t]})\simeq \QC(S)\otimes_k\Mod_{k[t]}.
\]

Let $\HH_\bullet(\CCC/S)$ be the Hochschild homology of $\CCC$
over $S$, that is an object of $\QC(S)^{S^1}$ (see Example~\ref{geometricex2}).
We define the negative cyclic homology/complex $\mathcal{HN}_\bullet(\CCC/S)$
in $\QC(S)\otimes_k\Mod_{k[t]}$
to be the image
$\Phi_S^{\textup{QC}}(\HH_\bullet(\CCC/S))$ in  $\QC(S)\otimes_k\Mod_{k[t]}$.
The periodic cyclic homology/complex $\mathcal{HP}_\bullet(\CCC/S)$ in $\QC(S)\otimes_k\Mod_{k[t,t^{-1}]}$ is defined as the image of $\mathcal{HN}_\bullet(\CCC/S)$
under the canonical functor 
$\QC(S)\otimes_k\Mod_{k[t]}  \to \QC(S)\otimes_k\Mod_{k[t,t^{-1}]}$.
As with Hochschild homology, we refer to the {\it complex}
$\mathcal{HP}_\bullet(\CCC/S)$ as the periodic cyclic {\it homology}.

\begin{Remark}
When $S=\Spec A$, the above definition of negative cyclic homology and
periodic cyclic homology recover the standard definitions in term of
$S^1$-fixed points and Tate construction. 
In this situation, $\Psi^{\textup{QC}}_S$ is the composite
\[
\Ind(\Perf_A^{S^1})\simeq \Ind(\Perf_A\otimes_k\Perf_{k[t]}) \to \Mod_A^{S^1} \simeq \Mod_A\otimes_k \Mod_{k[\epsilon]}\simeq \Mod_{A[\epsilon]}
\]
where $k[\epsilon]$ is the free commutative dg algebra generated by one element $\epsilon$ of homological degree one.
Here, if we write $A[t]$ for $A\otimes_kk[t]$
we use an equivalence $\Perf_A^{S^1}\simeq \Coh(A)^{S^1}\stackrel{\sim}{\to} \Coh(A)\otimes_k\Perf_{k[t]}\simeq \Perf_{A[t]}$
which carries $M$ to the homotopy fixed points $M^{S^1}$
with the canonical module structure over $A^{S^1}=A[t]$,
where we use the smoothness of $A$ which implies $\Perf_A\simeq \Coh(A)$.
Then $\Perf_A^{S^1}\simeq \Perf_{A[t]}$
gives rise to
\[
\Ind(\Perf_A^{S^1})\simeq \Ind(\Perf_{A[t]})\simeq \Mod_{A[t]}\simeq \Mod_A\otimes_k\Mod_{k[t]}.
\]
We note that the inverse equivalence $\Perf_{A[t]}\to \Perf_A^{S^1}$
is given by $A\otimes_{A[t]}(-):N\mapsto A\otimes_{A[t]}N$ with the module structure over the endomorphism algebra $\End_{\Mod_{A[t]}}(A)\simeq A[\epsilon]$.
Thus, we can regard $\Ind(\Perf_A^{S^1})\to \Mod_A^{S^1}\simeq \Mod_{A[\epsilon]}$ as
$\Mod_{A[t]} \to \Mod_{A[\epsilon]}$ given by $N\mapsto A\otimes_{A[t]}N$.
The right adjoint functor of $\Mod_{A[t]} \to \Mod_{A[\epsilon]}$ is given by $(-)^{S^1}:\Mod_A^{S^1}\simeq \Mod_{A[\epsilon]}\to \Mod_{A[t]};\ N\mapsto N^{S^1}$ (where $N^{S^1}$ can be represented by the (derived) Hom complex $\mathcal{H}om_{\Mod_{A[\epsilon]}}(A,N)$).
Namely,  $(\Phi^{\textup{QC}}_S,\Psi^{\textup{QC}}_S)$ can be identified with
\[
\xymatrix{
A\otimes_{A[t]}(-):\Mod_{A[t]}\simeq \Ind(\Perf_A^{S^1}) \ar@<0.5ex>[r] &   \Mod_A^{S^1} :(-)^{S^1}.\ar@<0.5ex>[l]  
}
\]
Consequently, the right adjoint $\Phi^{\textup{QC}}_S:\Mod_A\otimes_k \Mod_{k[\epsilon]}\to \Mod_A\otimes \Mod_{k[t]}$ of $\Psi^{\textup{QC}}_S$ is equivalent to the functor $(-)^{S^1}$.
This functor sends  $N\in \Mod_A^{S^1}$ to the $S^1$-invariants $N^{S^1}\in \Mod_{A[t]}$.
Thus, it sends $\HH_\bullet(\CCC/A)$ to $\HH_\bullet(\CCC/A)^{S^1}$.
\end{Remark}

\subsection{}

\begin{Lemma}
Let $H\in \QC_!(LS)^{S^1}$ and suppose that 
$(\iota^{!})^{S^1}(H)\in \QC_!(S)^{S^1}$ is equivalent to
$\varUpsilon_A(\HH_\bullet(\CCC/S))$ with the (standard) $S^1$-action.
Let $X$ denote the image of $H$ in $\QC_!(S_{\DR})\otimes_k \Mod_{k[t,t^{-1}]}$.
Then the forgetful functor $u_S:\QC_!(S_{\DR})\otimes_k \Mod_{k[t,t^{-1}]}\to \QC_!(S)\otimes_k \Mod_{k[t,t^{-1}]}$
sends $X$ to $\Upsilon_A(\mathcal{HN}_\bullet(\CCC/S))\otimes_{k[t]}k[t,t^{-1}]$ in $\QC_!(S)\otimes_{k}\Mod_{k[t,t^{-1}]}$ (i.e., the image of the periodic cyclic homology).
\end{Lemma}

\Proof
We first recall that the forgetful functor $u_S$ is $i^!\otimes  \Mod_{k[t,t^{-1}]}:\Ind(\Coh(S_{\DR}))\otimes_k \Mod_{k[t,t^{-1}]}\to \QC_!(S)\otimes_k \Mod_{k[t,t^{-1}]}$
where $i^!:\Ind(\Coh(S_{\DR}))\simeq \QC_!(S_{\DR})\to \QC_!(S)$ 
is
the (colimit-preserving) $!$-pullback along $i:S\to S_{\DR}$.
This $!$-pullback is a right adjoint (cf. \cite[Vol.II, Chap.4, 3.1.4, .3.2.2]{Gai2}). We write $i_*$ for the left adjoint.
The left adjoint $i_*$ sends objects of $\Coh(S)$ to compact objects of $\QC_!(S_{\DR})$.
Therefore, we can regard $i_*$ as a functor obtained from the restriction $\Coh(S)\to \Coh(S_{\DR})$
by the Ind-construction.
Note that $(i_*)^{S^1}|_{\Coh(S)^{S^1}}:\Coh(S)^{S^1}\to \Coh(S_{\DR})^{S^1}$ can naturally be identified with $i_*\otimes \Perf_{k[t]}:\Coh(S)\otimes_k\Perf_{k[t]} \to  \Coh(S_{\DR})\otimes_k\Perf_{k[t]}$.
It follows that $\Ind((i_*)^{S^1}|_{\Coh(S)^{S^1}}):\Ind(\Coh(S)^{S^1})\to \Ind(\Coh(S_{\DR})^{S^1})$
can be identified with 
$i_*\otimes \Mod_{k[t]}:\QC_!(S)\otimes_k\Mod_{k[t]} \to  \QC_!(S_{\DR})\otimes_k\Mod_{k[t]}$.
The right adjoint of $i_*\otimes \Mod_{k[t]}$ is given by
$i^!\otimes \Mod_{k[t]}:\QC_!(S_{\DR})\otimes_k\Mod_{k[t]}\to \QC_!(S)\otimes_k\Mod_{k[t]} $.
The base change to $k[t,t^{-1}]$ yields $u_S:\Ind(\Coh(S_{\DR}))\otimes_k \Mod_{k[t,t^{-1}]}\to \QC_!(S)\otimes_k \Mod_{k[t,t^{-1}]}$.
Namely, $u_S=i^!\otimes  \Mod_{k[t,t^{-1}]}$ is a right adjoint of 
$\Ind((i_*)^{S^1}|_{\Coh(S)^{S^1}})\otimes_{\Mod_{k[t]}}\Mod_{k[t.t^{-1}]}$.

Let $\Ind(\Coh(S)^{S^1})\to \Ind(\Coh(LS)^{S^1})$ be a colimit-preserving functor which is the $\Ind$-extension
of the pushforward $(\iota^{\textup{IndCoh}}_*)^{S^1}:\Coh(S)^{S^1}\to \Coh(LS)^{S^1}$.
Let $r_{\iota}:\Ind(\Coh(LS)^{S^1})\to \Ind(\Coh(S)^{S^1})$ be its right adjoint.
Let $(\iota^!)^{S^1}:\QC_!(LS)^{S^1}\to \QC_!(S)^{S^1}$ be a right adjoint functor of
$(\iota_*^{\textup{IndCoh}})^{S^1}$, which is induced by
the $!$-pullback right adjoint functor $\iota^!:\QC_!(LS)\to \QC_!(S)$ of $\iota_*^{\textup{IndCoh}}:\QC_!(S)\to \QC_!(LS)$.
Define
\[
\Psi_S:\Ind(\Coh(S)^{S^1})\to \QC_!(S)^{S^1} 
\]
to be the essentially unique colimit-preserving functor
which extends the inclusion $\Coh(S)^{S^1}\to  \QC_!(S)^{S^1} $.
Let $\Phi_S$ denote the left adjoint of $\Psi_S$. 
By construction, we have a commutative diagram 
\[
\xymatrix{
\QC_!(LS)^{S^1}   \ar[r]^(0.4){\Phi_{LS}} \ar[d]_{(\iota^{!})^{S^1}} & \Ind(\Coh(LS)^{S^1}) \ar[d]^{r_\iota}  \ar[r] &\Ind(\Coh(LS)^{S^1})\otimes_{k[t]}k[t,t^{-1}] \ar[d]^{r_\iota\otimes_{k[t]}k[t,t^{-1}]}  \\
\QC_!(S)^{S^1}  \ar[r]^(0.4){\Phi_S}   & \Ind(\Coh(S)^{S^1})  \ar[r] & \Ind(\Coh(S)^{S^1})\otimes_{k[t]}k[t,t^{-1}].
}
\]
Let $\textup{DL}_S:\QC_!(S_{\DR})\otimes_{k}k[t,t^{-1}]\to \Ind(\Coh(LS)^{S^1})\otimes_{k[t]}k[t,t^{-1}]$ be a right adjoint of
$\textup{LD}_S:\Ind(\Coh(LS)^{S^1})\otimes_{k[t]}k[t,t^{-1}] \to   \QC_!(S_{\DR})\otimes_{k}k[t,t^{-1}]$.
By the definitions of right adjoint functors, there is a canonical equivalence
$u_S\simeq r_\iota\otimes_{k[t]}k[t,t^{-1}] \circ \textup{DL}_S$.
By our assumption,
$(\iota^{!})^{S^1}(H)\simeq \varUpsilon_S(\HH_\bullet(\CCC/S))$.
Since $\textup{LD}_S$ is an equivalecne (in particular, it is fully faithful),  we are reduced to showing that 
$\Phi_S:\QC_!(S)^{S^1}  \to \Ind(\Coh(S)^{S^1})\simeq \QC_!(S)\otimes_k\Mod_{k[t]}$ sends
 $\varUpsilon_S(\HH_\bullet(\CCC/S))$ with the $S^1$-action to 
$\varUpsilon_S(\mathcal{HN}_\bullet(\CCC/S))$ (that is, the image of $\mathcal{HN}_\bullet(\CCC/S)$
under $\varUpsilon_S\otimes \textup{id}:\QC(S)\otimes_k \Mod_{k[t]}\stackrel{\sim}{\to} \QC_!(S)\otimes_k \Mod_{k[t]}$).
It is clear because the diagram
\[
\xymatrix{
\QC(S)^{S^1} \ar[r]^(0.4){\Phi_S^{\textup{QC}}} \ar[d]^{\varUpsilon_S} & \QC(S)\otimes_k \Mod_{k[t]} \ar[d]^(0.4){\varUpsilon_S\otimes \textup{id}} \\
\QC_!(S)^{S^1} \ar[r]^(0.4){\Phi_S}  &  \QC_!(S)\otimes_k \Mod_{k[t]}
}
\]
commutes.
\QED

By Construction~\ref{indequivariant} and the above Lemma, we have:

\begin{Theorem}
\label{maindmodule}
Suppose that $S$ is a smooth scheme of finite type over $k$.
We have constructed a right $\ZZ/2\ZZ$-periodic crystal (D-modules) 
as the object $\Omega^\circ(\CCC) \in \Ind(\Coh(S_{\DR}))\otimes_k \Mod_{k[t,t^{-1}]}$ which lies over $\mathcal{HP}_\bullet(\CCC/S)$. 
More precisely, there exists a canonically defined objects  $\Omega^\circ(\CCC)$
of $\textup{Crys}(S)\otimes_{k}k[t,t^{-1}]$
whose image in $\QC_!(S)\otimes_{k}k[t,t^{-1}]$ under $u_S$ is
equivalent to $\Upsilon_S(\mathcal{HN}_\bullet(\CCC/S))\otimes_{k[t]}k[t,t^{-1}]$.

Assume that $S$ is affine. The same statement holds
for $\Omega^\bullet(\CCC)$.
\end{Theorem}

\begin{Remark}
The image of $\Omega^{\circ,\textup{int}}(\CCC)$ and $\Omega^{\bullet,\textup{int}}(\CCC)$ under $\Ind(\Coh(LS)^{S^1})\to \Ind(\Coh(S)^{S^1})$
is canonically equivalent to the negative cyclic homology/complex $\mathcal{HN}_\bullet(\CCC/S)$. 
Moreover, $\Omega^{\circ}(\CCC)$ and $\Omega^{\bullet}(\CCC)$ are obtained from $\Omega^{\circ,\textup{int}}(\CCC)$ and $\Omega^{\bullet,\textup{int}}(\CCC)$
by inverting $t$.
Therefore, the pairs $(\Omega(\CCC)^\circ,\Omega^{\circ,\textup{int}}(\CCC))$ and $(\Omega^{\bullet}(\CCC),\Omega^{\bullet,\textup{int}}(\CCC))$ can be thought of as a categorical counterpart of  (the Rees construction of) the
Hodge filtered D-module.

\end{Remark}


\begin{thebibliography}{99}


\bibitem{And}
R. Andrade,
From manifolds to invariants of En-algebras,
 PhD thesis, Massachusetts Institute of
Technology, 2010




\bibitem{AF}
D. Ayala and J. Francis,
Factorization homology of topological manifolds,
J. Topol. Vol. 8 (2015)  1045--1084.



\bibitem{AFF}
D. Ayala and J. Francis,
Fibrations of $\infty$-categories,
Higher Structures 4 (1) 168--265, (2020).









\bibitem{AFT}
D. Ayala, J. Francis, and H.L. Tanaka,
Factorization homology of stratified manifolds,
Selecta Math. (N.S.) 23, no. 1, 293--362 (2017)

\bibitem{Bal}
E. Balzin,
Reedy model structures in Families, 
arXiv:1803.00681v2



\bibitem{Bar}
C. Barwick,
On left and right model categories and left and right Bousfield localizations,
Homology, Homotopy and Applications, vol. 12, 2010, 245–-320.




\bibitem{BN}
D. Ben-Zvi and D. Nadler,
Loop spaces and connections,
J. Topol., (2012) 1--54.


\bibitem{BGT1}
A. Blumberg, D. Gepner and G. Tabuada,
A universal characterization of higher K-theroy,
Geometry and Topology, 17 (2013), 733–-838.


\bibitem{BGT2}
A. Blumberg, D. Gepner and G. Tabuada,
Uniqueness of the multiplicative cyclotomic trace,
Adv. Math. 260 (2014) 191-232.


\bibitem{BM}
A. Blumberg and M. Mandell,
Localization theorems in topological Hochschild homology
and topological cyclic homology,
Geometry and Topology 16 (2012) 1053–-1120.




\bibitem{Cis}
D.-C. Cisinski,
Higher categories and Homotopical Algebra, 
Cambridge studies in advanced matematics 180,  2019.




\bibitem{Ind}
D. Gaitsgory,
Ind-coherent sheaves,
Mosc. Math. J., (2013), Vol. 13, 399-–528.

\bibitem{Gai2}
D. Gaitsgory and N. Rozenblyum,
A study in Derived Algebraic Geometry Volume I, II,
Mathematical Survey and Monographs, 22, American Math. Soc. 2017


\bibitem{Get}
Ezra Getzler,
Cartan homotopy formulas and the Gauss-Manin connection in cyclic homology,
In Quantum deformations of algebras and their representations (Ramat-Gan, 1991/1992; Rehovot,1991/1992)
volume 7 of Israel Math. Conf. Proc., pages 65–78. Bar-Ilan Univ., Ramat Gan, 1993.



\bibitem{Ginot}
G. Ginot,
Notes on factorization algebras, factorization homology and applications,



\bibitem{Good}
T. Goodwillie,
Cyclic homology, derivations, and the free loopspace,
Topology, 24 187–-215,
1985.


\bibitem{HP}
Y. Harpaz and M. Prasma,
Grothendieck construction for model
categories,
Adv. Math.
281, (2015), 1306--1363.



\bibitem{H}
B. Hennion,
Tangent Lie algebras of derived Artin stacks,
J. Reine Angew. Math. 741 (2018), 1435-5345.






\bibitem{Ho}
G. Horel,
Factorization homology and
calculus \'{a} la Kontsevich-Soibelman.
J. Noncomm. Geom. 11 (2017) 703--740.



\bibitem{HSS}
M. Hovey, B. Shilpley, and J. Smith,
Symmetric spectra,
J. Amer. Math. Soc. 13 (2000), 149--208.







\bibitem{I}
I. Iwanari,
Differnetial calculus of Hochschild pairs for infinity-categories,
SIGMA $\mathbf{16}$ (2020), 97  57 pages,
Special Issue on Primitive Forms and Related Topics in honor of Kyoji Saito for his 77th birthday.






\bibitem{IMA}
I. Iwanari,
Moduli theory associated to Hochschild pairs,
preprint available at the author's webpage


\bibitem{DC2}
I. Iwanari,
On D-modules of categories II,
preprint 

\bibitem{DC3}
I. Iwanari,
On D-modules of categories III,
forthcoming





\bibitem{L}
J. L. Loday,
Cyclic Homology,
Springer.



\bibitem{HTT}
J. Lurie,
Higher Topos Theory,
Annals Math. Studies, 170  2009.


\bibitem{HA}
J. Lurie,
Higher Algebra,
draft 2017.



\bibitem{DAG}
J. Lurie,
Derived Algebraic Geometry Series,
preprint



\bibitem{Maz}
A. Mazel-Gee,
Goerss-Hopkins obstruction theory via model $\infty$-categories,
thesis University of California, Berkeley (2016).


\bibitem{Muro}
F. Muro,
Dwyer-Kan homotopy theory of enriched categories,
J. Topol.  Vol. 8, (2015), 377–-413,


\bibitem{PS1}
D. Pavlov and J. Scholbach,
Admissibility and rectification of colored symmetric operads,
J. Topol. 11  (2018) 559--601.


\bibitem{PS2}
D. PAvlov and J. Scholbach,
Symmetric operads in abstract symmetric spectra,
J. Math. Ins. Jussiue 18 (2018), 707--758.





\bibitem{P}
A. Preygel,
Ind-coherent complexes on loop spaces and connections,
in ``Stacks and categories in geometry, topology, and algebra'' 643, 289-323, 2014.


\bibitem{PMF}
A. Preygel,
Thom-Sebastiani and Duality for Matrix Factorizations,
arXiv:1101.5834











\bibitem{S}
B. Shipley,
A convenient model structure for commutative ring spectra,
Contemp. Math. 346 (2004), 473--483.




\bibitem{HSSS}
M. Hoyois, P. Safronov, S. Scherotzke and N. Sibilla,
The categorified Grothendieck-Riemann-Roch theorem,
arXiv:1804.00879.




\bibitem{TV2}
B. Toen and G. Vezzosi,
Homotopical Algebraic Geometry II: Geometric stacks and Applications,
Memoirs Amer. Math. Soc.  902 (2008)






\end{thebibliography}
\end{document}